\documentclass{amsart}

\usepackage{graphicx}
\usepackage{graphics}
\usepackage{amssymb, amsmath, amsthm}
\usepackage{geometry}
\usepackage{lipsum}
\usepackage{enumerate}
\usepackage{dsfont}
\usepackage{mathrsfs}
\usepackage{xcolor}
\usepackage{multicol}
\usepackage{dblfloatfix}
\input xy
\xyoption{all}

\usepackage{hyperref}

\newtheorem{theorem}{Theorem}[section]
\newtheorem{lemma}[theorem]{Lemma}
\newtheorem{corollary}[theorem]{Corollary}
\newtheorem{proposition}[theorem]{Proposition}

\theoremstyle{definition}

\theoremstyle{remark}
\newtheorem{remark}[theorem]{Remark}
\newtheorem{question}{Question}

\numberwithin{equation}{section}
\renewcommand{\a}{\mathbf{a}}
\renewcommand{\b}{\mathbf{b}}
\renewcommand{\t}{\mathbf{t}}
\newcommand{\T}{\mathbf{T}}

\newcommand{\mc}{\mathcal}
\newcommand{\ms}{\mathscr}
\newcommand{\mb}{\mathbb}
\newcommand{\mf}{\mathfrak}
\newcommand{\on}{\operatorname}
\renewcommand{\v}{\mathbf{v}}
\newcommand{\w}{\mathbf{w}}
\newcommand{\x}{\mathbf{x}}
\newcommand{\y}{\mathbf{y}}

\newcommand{\s}{{\scriptscriptstyle \leqslant }}
\renewcommand{\leq}{\leqslant}
\newcommand{\n}{{\scriptscriptstyle\lhd}}

\newcommand{\ds}{\displaystyle}

\begin{document}
\title{Zeta functions of the 3-dimensional almost-Bieberbach groups}

\author{D. Sulca}
\address{Facultad de Matem\'aticas, Astronom\'ia y F\'isica, Universidad Nacional de C\'ordoba and Centro de Investigaci\'on y Estudios de Matem\'aticas, Ciudad Universitaria, C\'ordoba X5000HUA, Argentina}
\email{sulca@famaf.unc.edu.ar}









\begin{abstract}
The subgroup zeta function and the normal zeta function of a finitely generated virtually nilpotent group can be expressed as finite sums of Dirichlet series admitting Euler product factorization. We compute these series except for a finite number of local factors when the group is virtually nilpotent of Hirsch length 3. We deduce that they can be meromorphically  continued to the whole complex plane and that they satisfy local functional equations. The complete computation (with no exception of local factors) is presented for those groups that are also torsion-free, that is, for the 3-dimensional almost-Bieberbach groups.
\end{abstract}

\maketitle	

\section{Introduction}
\noindent The {\em subgroup zeta function} and the {\em normal zeta function} of a finitely generated group $G$ are, respectively, the Dirichlet series 
\begin{align*}
\zeta_G^{\s}(s)=\sum_{n=1}^\infty\frac{a_n^\s(G)}{n^s}\quad\quad \mbox{and}\quad\quad  \zeta_G^\n(s)=\sum_{n=1}^\infty \frac{a_n^\n(G)}{n^s},
\end{align*}
where $a_n^\s(G)$ denotes the number of subgroups and $a_n^\n(G)$  the number of normal subgroups of index $n$.
They were introduced  by Grunewald, Segal, and Smith in \cite{GSS} as a means to study groups of polynomial subgroup growth (PSG).
Indeed, these are precisely the groups for which these series are not only formal gadgets but also define analytic functions on some right-half plane of the complex plane. A natural problem is to understand how structural information of a group of PSG is encoded into arithmetical and analytic properties of its zeta functions and vice-versa. 

Lubtozky, Mann and Segal characterized the finitely generated residually finite groups of PSG: these are precisely the virtually solvable groups of finite rank \cite{LMS}. On the other hand, since they were introduced, zeta functions of groups have been extensively studied only in when the group is a finitely generated torsion-free nilpotent group: a $\mf{T}$-group. The theory for these groups is rich and continues growing.
We summarize some general results and refer the reader to the survey \cite{V2} and the references therein for more information and variations of the topic. Let $N$ be a $\mf{T}$-group, let $*\in\{\leq,\lhd\}$ and let $\alpha_N^*$ denote the abscissa of convergence of $\zeta_N^*(s)$. Then the following holds.
\begin{enumerate}[(N1)]
	\item $\alpha_N^*\leq h(N)$ (= the Hirsch length of $N$) \cite[Proposition 1]{GSS}.
	\item $\alpha_N^*\in\mb{Q}$ and there is $\delta>0$ such that $\zeta_N^*(s)$ admits a meromorphic continuation to the region $\{s\in\mathbb{C}: \on{Re}(s)>\alpha_N^*-\delta\}$  \cite[Theorem 1.1]{dSG}.
	\item $\zeta_N^*(s)$ has an Euler product factorization
	\begin{align*}
		\zeta_N^*(s)=\prod_{p\ \textrm{prime}} \zeta_{N,p}^*(s)
	\end{align*}
where $
	\zeta_{N,p}^*(s)=\sum_{k=0}^\infty a_{p^k}^*(N){p^{-ks}}.$
 In addition, each local factor $\zeta_{N,p}^*(s)$ is a rational function in $p^{-s}$ with coefficients in $\mb{Q}$ \cite[Proposition 1.3 and Theorem 1]{GSS}.
\item For almost all prime $p$, $\zeta_{N,p}^\s(s)$ satisfies the functional equation 
\begin{align*}
	\zeta_{N,p}^\s(s)|_{p\to p^{-1}}=(-1)^h p^{\binom{h}{2}-hs}\zeta_{N,p}^\s(s),
\end{align*}
where $h=h(N)$ \cite[Corollary 1.1]{V1}. If in addition $N$ is a $\mf{T}_2$-group (a $\mf{T}$-group of nilpotency class 2), then for almost all prime $p$, $\zeta_{N,p}^{\n}(s)$ satisfies the functional equation
\begin{align*}
	\zeta_{N,p}^\n(s)|_{p\to p^{-1}}=(-1)^hp^{\binom{h}{2}-(d+h)s}\zeta_{N,p}^\n(s),
\end{align*}
where $d$ is the rank of $N/Z(N)$ \cite[Theorem C]{V1}.
\end{enumerate}
\begin{remark} Properties (N1)-(N4) are, in fact, corollaries of the analogous properties established for {\em the subring zeta function} and {\em the ideal zeta function} of nilpotent Lie rings additively isomorphic to some $\mathbb{Z}^h$. The translation is done via the Mal'cev correspondence. The question of whether the ideal zeta function of a nilpotent Lie ring of nilpotency class $>2$ satisfies local functional equations led to the introduction and investigation of {\em the submodule zeta function}; see, e.g., \cite{Rossmann2016}, \cite{V3}, and \cite{LeeVoll2020}.
	On the other hand, it has long been known that, for instance, the ideal zeta function of the filiform nilpotent Lie ring $\on{Fil}_4$ fails to satisfy local functional equations; see \ \cite[Theorem 2.39]{dSW}. 
\end{remark}

Regarding the computation of zeta functions, the following are classical examples.
\begin{enumerate}
	\item For the free abelian group $\mb{Z}^h$, it holds that $\zeta_{\mb{Z}^h}^\s(s)=\zeta(s)\zeta(s-1)\cdots\zeta(s-h+1)$, where $\zeta(s)$ is the Riemann zeta function.
	In \cite[Chap.\ 15]{LS}, there are five different proofs of this elementary fact.
	\item Let $\mathbf{H}(\mb{Z})$ denote the integral points of the Heisenberg unipotent group scheme $\mathbf{H}$. Then
	\begin{align*}
		\zeta_{\mathbf{H}(\mb{Z})}^\s(s)=\frac{\zeta(s)\zeta(s-1)\zeta(2s-2)\zeta(2s-3)}{\zeta(3s-3)}\quad\quad  \mbox{and}\quad\quad  \zeta_{\mathbf{H}(\mb{Z})}^\n(s)=\zeta(s)\zeta(s-1)\zeta(3s-2)
	\end{align*}
	\cite[Section 8]{GSS}. 	A general theory and calculation of $\zeta_{\mathbf{H}_3(\mc{O}_K)}^\n(s)$, for $K$ a number field and $\mc{O}_K$ its ring of integers, are presented in \cite{SV1} and \cite{SV2}. 
\end{enumerate}
 \noindent Finding explicit formulae for the zeta functions of a $\mf{T}$-group is in general a difficult task. A substantial list of examples is recorded in \cite{dSW}, and new ones have emerged more recently; see, e.g., \cite{CSV}, \cite{SV1}, \cite{SV2}, \cite{V4}. In \cite{Rossmann2015} and \cite{Rossmann2018}, Rossmann developed a method for computing certain zeta functions associated with $\mf{T}$-groups and other algebraic structures. This algorithm was implemented in  \cite{Rossmann2017}.

\medskip

This work studies and presents formulae for the zeta functions of the virtually nilpotent groups of Hirsch length 3 (see the next paragraph for a detailed description of the content). This is not the first attempt in dealing with zeta functions of groups that are not nilpotent. In \cite{dS}, du Sautoy investigated zeta functions of compact $p$-adic analytic groups ($=$ virtually uniform pro-$p$-groups) and proved that they are rational functions in $p^{-s}$. Then, in \cite{dSMS}, the authors studied zeta functions of virtually abelian groups, proving for example, that they can be meromorphically continued to the whole complex plane (this property is not shared by zeta functions of $\mf{T}$-groups in general, as discussed in \cite[Chap. 7]{dSW}). The paper \cite{Su} contains general properties of zeta functions of virtually nilpotent groups that we now recall.
 Let $G$ be a finitely generated virtually nilpotent group, and let $N\lhd G$ be a finite index normal subgroup that is a $\mf{T}$-group. It is easy to check that
\begin{align*}
	\zeta_G^\s(s)=\sum_{N\leq H\leq  G}[G:H]^{-s}\zeta_{G}^{H,\s}(s)\quad \quad\mbox{and}\quad\quad  \zeta_G^\lhd(s)=\sum_{N\leq H\lhd G}[G:H]^{-s}\zeta_{G}^{H,\n}(s),
\end{align*}
where 
\begin{align*}
	\zeta_{G}^{H,\s}(s):=\sum_{\substack{A\leq G:\ AN=H}}[H:A]^{-s}\quad \quad \mbox{and}\quad \quad \zeta_{G}^{H,\n }(s):=\sum_{\substack{A\lhd G:\ AN=H}}[H:A]^{-s}.
\end{align*}
The series $\zeta_G^{H,\s}(s)$ and $\zeta_G^{H,\n}(s)$ will be referred to as {\em partial zeta functions} of $\zeta_G^\s(s)$ and $\zeta_G^\lhd(s)$ (with respect to $N$).
Each of them is an Euler product of {\em cone integrals} \cite[Theorem 1]{Su}, and hence it satisfies properties (N2) and (N3) above by the general theory of cone integrals developed in \cite{dSG}. 
Property (N1) also holds with the bound $h(N)=h(G)$ replaced by $h(G)+1$ \cite[Proposition 2.12]{Su}. 
The partial zeta functions were computed explicitly, except for a finite number of local factors, when the group $N$ is abelian; see  \cite[Section 2]{dSMS}. Computing the exceptional factors is, in general, more involved. The complete computation was done for the plane crystallographic groups. This is the main result in \cite{Mc}, and it is also summarized in \cite[Section 4]{dSMS}.  To the author's knowledge, these are the only existing examples of computations of zeta functions of infinite groups that are not nilpotent.

\medskip

The article is organized as follows.
In Section \ref{Section: Local factors at good primes}, we obtain formulae for the local factors of $\zeta_{G}^{G,\s}(s)$ and $\zeta_G^{N,\n}(s)$ when $p\nmid [G:N]$ for every virtually nilpotent group of Hirsch length 3 (see Table \ref{Table with zeta_G^G} and Table \ref{Table with zeta_G^N}). This suffices to conclude that the subgroup and normal zeta functions of each of these groups extend to meromorphic functions on the whole plane (Corollary \ref{Meromorphic continuation partial subgroup zeta function} and Corollary \ref{meromorphic contunuation normal zeta function}). We also deduce local functional equations, similar to but different from the ones presented in (N4) (Corollary \ref{main corollary counting subgroups} and Corollary \ref{main corollary normal zeta functions}).
The arguments in this section are of a group-theoretical nature, akin to \cite[Section 2]{dSMS}. 
In section \ref{Section: method of p-adic integration}, we review a method for expressing local zeta functions of groups in terms of $p$-adic integrals. This method is used later to compute local factors at primes $p|[G:N]$ (those not considered in Section \ref{Section: Local factors at good primes}). Section \ref{Section: The formulae} presents the formulae for the zeta functions of the torsion-free virtually nilpotent groups of Hirsch length 3. A group like this is either a 3-dimensional Bieberbach group or a 3-dimensional almost-Bieberbach group (the fundamental group of a 3-dimensional infra-nilmanifold). Section \ref{Section: Proofs for Bieberbach groups} contains proofs of the formulae for the 3-dimensional Bieberbach groups. It also contains formulae for the zeta functions of a large family of Bieberbach groups with holonomy $C_2$, to illustrate how involved the method for computing local factors at bad primes is, even in the case of virtually abelian groups. 
Finally, Section \ref{Section: proof for AB-groups} contains proofs of the formulae  for the 3-dimensional almost Bieberbach groups.

\medskip 

{\em Notation and conventions}

\smallskip
\noindent $|S|$ denotes the cardinality of a set $S$.  
For a group $G$, $Z(G)$ denotes its center; $A\leq G$ means that $A$ is a subgroup; $A\lhd G$ means that $A$ is a normal subgroup; $[G:A]$ denotes the index of $A$ in $G$; $\on{N}_G(A)$ denotes the normalizer of $A$ in $G$. In sums involving $[G:A]$, only subgroups of finite index are considered.
For a subset $S\subset G$, $\langle S\rangle$ denotes the subgroup generated by $S$. When $G$ is profinite, $\overline{\langle S\rangle}$ denotes the closed subgroup generated by $S$.
For $x,y\in G$, we denote
  $ ^xy=xyx^{-1}$ and $[x,y]=xyx^{-1}y^{-1}(= {}^xy\cdot y^{-1})$. For $S,T\subset G$ we denote $[S,T]=\langle [s,t]| s\in S, t\in T\rangle$. We use without mentioning the fact that if $G$ is a finitely generated profinite group, then every finite index subgroup is open, and if $H$ and $K$ are closed subgroups, then $[H:K]$ is also closed; see \cite{NS}.

The variable $p$ will be reserved for prime numbers. We denote by 
$v_p:\mb{Q}_p\to\mb{Z}\cup \{\infty\}$ the usual $p$-adic valuation, and by $|\cdot |_p=p^{-v_p(\cdot)}$ the $p$-adic norm.
$\zeta_p(s)$ denotes $\frac{1}{1-p^{-s}}$, the local factor of the Riemann zeta function $\zeta(s)$ at $p$.

\section{Local factors at good primes and meromorphic continuation}\label{Section: Local factors at good primes} 
\noindent Let $G$ be a finitely generated virtually nilpotent group, and let $N\lhd G$ be a finite index normal subgroup that is a $\mf{T}$-group. For a prime $p$, let $G_p$ be the completion of $G$ with respect to the family $\{B\lhd G_p: B\subseteq N_p,\ [N_p:B]\ \mbox{a power of } p\}$. For $X\leq G$, we denote by $X_p$ the closure of $X$ in $G$. For a subquotient $X/Y$ of $G$ we denote $(X/Y)_p=X_p/Y_p$.
 Note that $N_p$ is the pro-$p$ completion of $N$, and that the inclusion $G\hookrightarrow G_p$ induces an isomorphism $G/N\cong G_p/N_p$.

As noted in the introduction, to study $\zeta_G^\s(s)$, it is enough to focus on the partial zeta functions $\zeta_{G}^{H,\s}(s)=\zeta_H^{H,\s}(s)$ for $N\leq H\leq G$. There is no loss of generality if we only consider $\zeta_G^{G,\s}(s)$. Similarly, to study $\zeta_{G}^{\n}(s)$, it is enough to consider $\zeta_{G}^{H,\n}(s)$ for $N\leq H\lhd G$. 
The advantage is that we have 
\begin{align*}
	\zeta_{G}^{G,\s}(s)=\prod_{p\ \textrm{prime}} \zeta_{G_p}^{G_p,\s}(s)\quad\quad\mbox{and}\quad\quad \zeta_{G}^{H,\n }(s)=\prod_{p\ \textrm{prime}}\zeta_{G_p}^{H_p,\n}(s),
\end{align*}
where $\zeta_{G_p}^{G_p,\s}(s)$ and $\zeta_{G_p}^{H_p,\n}(s)$ are computed with respect to $N_p$ \cite[Proposition 2.2]{Su}. 

Along the section, unless otherwise specified,  $N$ denotes a $\mf{T}_2$-group. We also fix the following notation:
\begin{align*}
	Z=Z(N),\quad T=N/Z,\quad P=G/N,\quad \Lambda=\mb{Z}[P]\quad \mbox{and}\quad \Lambda_p=\mb{Z}_p[P].
\end{align*}
We will first obtain general expressions for $\zeta_{G_p}^{G_p,\s}(s)$ and $\zeta_{G_p}^{H_p,\n}(s)$ when $p\nmid |P|$. Then we shall specialize to the case of Hirsch length 3, giving explicit formulae in terms of well-known series. 
 The case $N\cong\mb{Z}^3$ (and more generally $N\cong\mb{Z}^h$ for any $h$) was settled in \cite{dSMS}. 
  

The action by conjugation of $G$ on $N$ ($(g,n)\mapsto {}^g n=g n g^{-1}$) induces structures of left $\Lambda$-modules on the abelian groups $T$ and $Z$ (and hence $T_p$ and $Z_p$ become $\Lambda_p$-modules).
 For this reason, we use additive notation when working with them.
The commutator operation $[x,y]=xyx^{-1}y^{-1}$ induces a bilinear map $[\cdot,\cdot]:T\times T\to Z$ (and hence also a $\mb{Z}_p$-bilinear map $[\cdot,\cdot]:T_p\times T_p\to Z_p$) that is compatible with the action of $P$, that is, $[\alpha\cdot x,\alpha\cdot y]=\alpha\cdot [x,y]$ for all $\alpha\in P$ and $x,y\in T$. 
If $U_1$ and $U_2$ are $\mb{Z}_p$-submodules of $T_p$, then $[U_1,U_2]$ denotes the $\mb{Z}_p$-submodule (=subgroup) of $Z_p$ generated by $[x_1,x_2]$ with $x_i\in U_i$.

For a ring $R$ and a left $R$-module $M$, we use the notation $X\leq_R M$ to indicate that $X$ is a left $R$-submodule of $M$. For $r\in R$,  $r_M:M\to M$ denotes the left multiplication by $r$. 
If $R=\Lambda$ or $\Lambda_p$, we denote by $\on{Der}(P,M)$ the set of derivations from $P$ to $M$.
   For a subset $S\subseteq P$, we denote $C_M(S)=\{x\in M: \alpha\cdot x=x,\ \forall\alpha\in S\}$. We will use the fact that if $p\nmid |P|$ and $M$ is a finitely generated $\Lambda_p$-module that is torsion-free as a $\mb{Z}_p$-module, then $M$ is a projective $\Lambda_p$-module, and all finite index $\Lambda_p$-submodules of $M$ are isomorphic. This follows from the fact that when $p\nmid |P|$,  $\Lambda_p$ is a maximal order of $\mb{Q}_p[P]$. 

\subsection{Counting subgroups}\label{counting subgroups}
We begin with a general formula for the local factors of $\zeta_G^{G,\s}(s)$ at ``good" primes.
\begin{proposition}\label{formula for the partial subgroup zeta function}
	If $p\nmid |P|$, then
	$$\zeta_{G_p}^{G_p,\s}(s)=\sum_{\substack{U\leq_{\Lambda_p} T_p,\ V\leq_{\Lambda_p} Z_p\\ [U,U]\subseteq V}} [T_p:U]^{-s}[T_p:C_{T_p}(P)+U][Z_p:V]^{-s}[Z_p:C_{Z_p}(P)+V]|\on{Hom}_{\Lambda_p}(T_p,Z_p/V)|.$$	
\end{proposition}
\begin{proof}
	For a subgroup $U\leq T_p$, we denote by $\tilde{U}$ its pre-image in $N_p$.	We set
	\begin{align*}
		\mc{A}&:=\{A\leq G_p:\ AN_p=G_p,\ [G_p:A]<\infty \},\\
		\mc{U}&:=\{(U\leq_{\Lambda_p} T_p, V\leq_{\Lambda_p} Z_p):\  [T_p:U]<\infty,\ [Z_p:V]<\infty,\ [U,U]\subseteq V\},\\
		\mc{V}&:=\{(U,V,C): (U,V)\in\mc{U},\  C\in\mc{A},\ C\cap N_p=\tilde{U}\},
	\end{align*}
	and define maps $\Phi:\mc{A}\to\mc{V}$, $\Psi:\mc{V}\to\mc{U}$ by 
	\begin{align*} 
		\Phi(A)=(((A\cap N_p)Z_p)/Z_p,A\cap Z_p, Z_pA)\quad\mbox{and}\quad \Psi(U,V,C)=(U,V).
	\end{align*}
	It is straightforward to check that $\Phi$ is well-defined. Note that for $A\in(\Psi\Phi)^{-1}(U,V)$, it holds that $[G_p:A]=[N_p:A\cap N_p]=[N_p:\tilde{U}][Z_p:V]=[T_p:U][Z_p:V]$; therefore,
	\begin{align*}
		\zeta_{G_p}^{G_p,\s}(s)&=\sum_{A\in \mc{A}} [G_p:A]^{-s}=\sum_{(U,V)\in\mc{U}}[T_p:U]^{-s}[Z_p:V]^{-s}|(\Psi\Phi)^{-1}(U,V)|.
	\end{align*}
	
	We now fix $(U,V)\in \mc{U}$ and show that $|\Psi^{-1}(U,V)|=|\on{Der}(P,T_p/U)|$. 
	Note first that $\Psi^{-1}(U,V)$ is in a bijection with the set of complements of $N_p/\tilde{U}$ in $G_p/\tilde{U}$. Since $N_p/\tilde{U}$ is a normal Hall subgroup of $G_p/\tilde{U}$, there is at least one complement by the Schur-Zassenhaus Theorem, and therefore the number of complements is $|\on{Der}(P,N_p/\tilde{U})|=|\on{Der}(P,T_p/U)|$ (cf.\ \cite[Proposition 1, Chap. 3]{Se}).
	
	We now fix $(U,V,C)\in\mc{V}$ and show that $ |\Phi^{-1}(U,V,C)|=|\on{Hom}_{\Lambda_p}(U,Z_p/V)||\on{Der}(P,Z_p/V)|$.
	The subgroups $\tilde{U}$ and $V$ are normal closed subgroups of $G_p$, and the condition  $[U,U]\subseteq V$ implies that the quotient $\tilde{U}/V$ is abelian. Thus, $\tilde{U}/V$ is a finitely generated $\mb{Z}_p$-module and $Z_p/V$ is clearly the torsion submodule.
	We use $C$ to give $\tilde{U}/V$ a structure of $\mb{Z}_p[P]$-module. Firstly, the action by conjugation of $C$ on $\tilde{U}$ induces a structure of $C/\tilde{U}$-module on $\tilde{U}/V$. Secondly, the inclusion $C\to G_p$ induces an isomorphism $C/\tilde{U}\cong G_p/N_p=P$. Therefore, $\tilde{U}/V$ becomes a $\mb{Z}_p[P]$-module. 
	Note that $Z_p/V$ has two structures of $\mb{Z}_p[P]$-module: one as a subobject of $\tilde{U}/V$ and one as a quotient of $Z_p$. It is easy to check that these two structures coincide.  Similarly, the structure of $\mb{Z}_p[P]$-module on $\tilde{U}/Z_p=U$ as a quotient of $\tilde{U}/V$ and the one as a subobject of $T_p$ are the same.
	
		Given $A\leq \Phi^{-1}(U,V,C)$, observe that $(A\cap N_p)/V$ is a complement of $Z_p/V$ in $\tilde{U}/V$. We claim that this complement is $P$-invariant. Indeed, $A\cap N_p$ is normal in $A$ and in $(A\cap N_p)Z_p$, so it is normal in $AZ_p=C$. Therefore, $(A\cap N_p)/V$ is $C/\tilde{U}$-invariant, and hence $P$-invariant. 
	Now, according to Lemma \ref{number of F-invariant complements} below, the number of $P$-invariant complements of $Z_p/V$ in $\tilde{U}/V$ is $|\on{Hom}_{\Lambda_p}(U,Z_p/V)|$. Fix one such complement, say $B/V$. 
	Note that $B$ is normal in $C$ (this follows from the condition of $B/V$ being $P$-invariant) and the set of those $A\in\Phi^{-1}(U,V,C)$ such that $A\cap N_p=B$ is in a bijection with the set of complements of $\tilde{U}/B$ in $C/B$.
	Since $\tilde{U}/B$ is a normal Hall subgroup of $C/B$, there is at least one such complement by the Schur-Zassenhaus Theorem, and then the number of complements is  $|\on{Der}(C/\tilde{U},\tilde{U}/B)|=|\on{Der}(P,\tilde{U}/B)|=|\on{Der}(P,Z_p/V)|$. We conclude that $|\Phi^{-1}(U,V,C)|=|\on{Hom}_{\Lambda_p}(U,Z_p/V)||\on{Der}(P,Z_p/V)|$. 
	
	To end the proof, we only need to show that $|\on{Der}(P,T_p/U)|=[T_p:U+C_{T_p}(P)]$, $|\on{Der}(P,Z_p/V)|=[Z_p:V+C_{Z_p}(P)]$ and  $|\on{Hom}_{\Lambda_p}(U,Z_p/V)|=|\on{Hom}_{\Lambda_p}(T,Z_p/V)|$. Here we use again that $p\nmid |P|$. The first two equalities follow from \cite[Lemma 2.4]{dSMS}, and for the last one we use the fact that $U$ and $T$ are isomorphic as $\Lambda_p$-modules.
\end{proof}
 The following lemma was used in the proof above.
\begin{lemma}\label{number of F-invariant complements}
	Let $R$ be a ring, let $M$ be a left $R$-module, and fix $X\leq_R M$. If there is a complement of $X$ in $M$ (i.e,. there is $Y\leq_R M$ such that $X+Y=M$ and $X\cap Y=0$), then the set of complements of $X$ in $M$ is in a bijection with $\on{Hom}_R(M/X,X)$. If $R=\mb{Z}_p[P]$ with $p\nmid |P|$, and if $M/X$ is finitely generated and torsion-free as a $\mb{Z}_p$-module, then the set of ($P$-invariant) complements of $X$ is in a bijection with $\on{Hom}_{\mb{Z}_p[P]}(M/X,X)$.
\end{lemma}
\begin{proof}
Assume that there is a complement of $X$, say $Y$. Given any other complement $K\leq_R M$, we define $\varphi_{K}:Y\to X$ as follows. If $y\in Y$, then we can write $y=k+x$ uniquely with $k\in K$ and $x\in X$. We set $\varphi_{K}(y)=x$. It is clear that $\varphi_{K}\in\on{Hom}_R(Y,X)$. Conversely, given $\varphi\in\on{Hom}_R(Y,X)$, we define $K_\varphi=\{y-\varphi(y): y\in Y\}$. It is easy to check that $K_\varphi$ is a complement of $X$, and that $\varphi_{K_\varphi}=\varphi$ and $K_{\varphi_{K}}=K$. Therefore, the set of complements of $X$ is in a bijection with $\on{Hom}_R(Y,X)$, which in turn is in a bijection with $\on{Hom}_R(M/X,X)$ since $M/X$ and $Y$ are isomorphic. This proves the first part of the lemma
 
We now show the second part, so we assume now that $R=\mb{Z}_p[P]$ with $p\nmid |P|$ and that $M/X$ is finitely generated and torsion-free as a $\mb{Z}_p$-module.  Then $M/X$ is a free $\mb{Z}_p$-module of finite rank, and hence it is a projective $\mb{Z}_p[P]$-module (here we use the fact that $p\nmid |P|$). This implies that there is at least one complement of $X$; therefore, the set of complements is in a bijection with $\on{Hom}_{\mb{Z}_p[P]}(M/X,X)$, as shown in the first part of the lemma. This completes the proof.
\end{proof}

The following lemma will be used to calculate $|\on{Hom}_{\Lambda_p}(T_p,Z_p/V)|$ in some particular cases. 
\begin{lemma}\label{Hom-trivial}
	Let $M$ be a $\Lambda_p$-module and assume that $p\nmid |P|$. If there is  $\alpha\in P$ such that $C_M(\alpha)=M$ and $C_{Z_p}(\alpha)=0$,  then $\on{Hom}_{\Lambda_p}(M,Z_p/V)=\{0\}$ for all $V\leq_{\Lambda_p} Z_p$. Similarly, if $C_M(\alpha)=0$ and $C_{Z_p}(\alpha)=Z_p$, then also $\on{Hom}_{\Lambda_p}(M,Z_p/V)=\{0\}$.
\end{lemma}
\begin{proof} 
Since $p\nmid |P|$ and since $Z_p$ is finitely generated and torsion-free as a $\mb{Z}_p$-module, $Z_p$ is a projective $\Lambda_p$-module. Thus, the lemma in both cases follows from the equality $\on{Hom}_{\Lambda_p}(M,Z_p)=0$, whose verification is straightforward.
\end{proof}

We now specialize to the case when $N$ has Hirsch length 3. Recall that given a finite subgroup $F\subset GL_2(\mathbb{Z})$, either $F$ is included in $SL_2(\mathbb{Z})$, in which case $F$ is isomorphic to one of the cyclic groups $C_1,C_2,C_3,C_4, C_6$, or else there is $\beta\in F$ of order 2 with determinant $-1$, in which case $F$ is isomorphic to one of the dihedral groups $D_1 (\cong C_2), D_2 (\cong C_2\times C_2), D_3, D_4, D_6$.
\begin{lemma}\label{main lemma}
	Assume that $N$ is a $\mf{T}_2$-group of Hirsch length 3 (so that $T\cong\mb{Z}^2$ and $Z\cong \mathbb{Z}$). Set $F:=\on{Im}(P\to GL(T))$ and  $\eta:=\on{rank}_{\mb{Z}}(C_T(P))$.
	 Then the following holds.
	\begin{enumerate}
					\item $C_Z(P)=Z$ if and only if $F\subset SL(T)$.
					\item $\eta=\left\{\begin{array}{ll}2&\mbox{if } F\ \mbox{is trivial}, \\ 1&\mbox{if } F\not\subset SL(T)\ \mbox{and } F\cong D_1, \\ 0 &\mbox{otherwise} \end{array}\right.$
					\item If $p\nmid |P|$, then for any $V\leq_P Z_p$ of finite index, we have $|\on{Hom}_{\Lambda_p}(T_p,Z_p/V)|=[Z_p:V]^\eta$.  \end{enumerate}
				\end{lemma}
			\begin{proof}  
				 Fix an ordered basis $\{x,y\}$ for $T$. Note that $[x,y]\in Z$ is non-zero.
				  Given $\alpha\in P$, let $\begin{pmatrix}
					a&b\\ c&d
				\end{pmatrix}$ be the matrix of $\alpha_T$ with respect to $\{x,y\}$. We have $\alpha\cdot [x,y]=[ax+cy,bx+dy]=(ad-bc) [x,y]$. It follows that $C_Z(P)=Z$ if and only if $F\subset SL(T)$. 
				 It also follows that  $Z$ is naturally an $F$-module. 
			 
	We now show (2) and (3). Fix a prime $p\nmid |P|$ and a subgroup $V\leq_{\Lambda_p} Z_p$ of finite index. Note that (2) and (3) hold clearly for $F=C_1$. 
	Assume next that $F\cong D_1$ and  $F\not\subset SL(T)$. Then the generator of $F$ has eigenvalues 1 and $-1$, whence the $\Lambda$-submodules					   $T_+:=C_T(P)$ and $T_-:=\{x\in T: \beta\cdot x=-x\}$ have additive rank 1. It follows that $\eta=1$. Next,
  since $p\neq 2$, the decomposition $x=\frac{x+\alpha\cdot x}{2}+\frac{x-\alpha\cdot x}{2}$ holds in $T_p$ and yields a decomposition $T_p=(T_+)_p\oplus (T_-)_p$. Since $C_Z(F)=0$ by (1), 
						  Lemma \ref{Hom-trivial} implies that $|\on{Hom}_{\Lambda_p}((T_+)_p,Z_p/V)|=1$. Thus,
						  $|\on{Hom}_{\Lambda_p}(T_p,Z_p/V)|
						  =|\on{Hom}_{\Lambda_p}((T_-)_p,Z_p/V)|=[Z_p:V]=[Z_p:V]^\eta$.
						  
						  Assume now that $C_1\neq F=\langle\alpha\rangle\subset SL(T)$.
						  It holds that $C_{T}(\alpha)=0$, and hence $\eta=0$, since otherwise 1 would be an eigenvalue of $\alpha$. However, since $\alpha\in SL(T)$ the latter  implies that $\alpha$ is the identity, which is not the case.
	Since $C_Z(\alpha)=Z$ by (1), Lemma \ref{Hom-trivial} implies that $|\on{Hom}_{\Lambda_p}(T_p,Z_p/V)|=1=|\on{Hom}_P(T_p,Z_p/V)|^\eta$.  
	
	 Assume finally that $F\cong D_d$ with $d\in\{2,3,4\}$. The intersection $F\cap SL(T)$ is not the trivial group. Indeed, the product of any two elements of $F$ is in $SL(T)$ and $F$ has more than two elements. It follows from the previous case that $\eta=0$ and $|\on{Hom}_{\Lambda_p}(T_p,Z_p/V)|=1=|\on{Hom}_P(T_p,Z_p/V)|^\eta$.  We have covered all the cases, so the proof is complete.
							\end{proof}
		\begin{lemma}\label{index of the commutator}
	Assume that $N$ is a $\mf{T}_2$-group of Hirsch length 3. If $U\leq T_p$ has finite index, then $[Z_p:[U,U]]=|[Z:[N,N]]|_p^{-1}[T_p:U]$		\end{lemma}
\begin{proof}
	Let $\{x,y\}$ be a basis for the $\mb{Z}_p$-module $T_p$, and let $U=\leq T_p$ be a $\mb{Z}_p$-submodule of finite index, say generated by $ax+by$ and $cx+dy$. Note that $[T_p,T_p]=\mb{Z}_p[x,y]$ and that $[U,U]=\mb{Z}_p [ax+by,cx+dy]=\mb{Z}_p (ad-bc)[x,y]$. Thus, $[Z_p:[U,U]]=[Z_p:[T_p,T_p]]|ad-bc|_p^{-1}$.	 On the other hand, clearly $[T_p:U]=|ad-bc|_p^{-1}$. Thus, $[Z_p:[U,U]]=[Z_p:[T_p,T_p]][T_p:U]$. Finally, it is clear that $[Z_p:[T_p,T_p]]=|[Z:[T,T]]|_p^{-1}=|[Z,[N,N]]|_p^{-1}$.
\end{proof}
\begin{theorem}\label{main theorem for zeta_G^G}
	Assume that $N$ is a $\mf{T}_2$-group of Hirsch length 3. Let $F:=\on{Im}(P\to GL(T))$ and $E=T\rtimes F$, which is a plane crystallographic group. If $p\nmid |P|$, then 
		\begin{align*}
		\zeta_{G_p}^{G_p,\s}(s)&=\zeta_p(s-\eta-\epsilon) \left(\zeta_{E_p}^{E_p,\s}(s) -p^{-s+\eta+\epsilon}|[Z:[N,N]]|_p^{s-\eta-\epsilon}\zeta_{E_p}^{E_p,\s}(2s-\eta-\epsilon) \right), 
	\end{align*}
where $\eta:=\on{rank}_{\mb{Z}}(C_T(F))$,  $\epsilon$ is $0$ or 1 according to whether $F$ is included in $SL(T)$ or not, and $\zeta_{E_p}^{E_p,\s}(s)$ is computed with respect to $T_p\lhd E_p$.
\end{theorem}
\begin{proof} 
	Fix $p\nmid |P|$. 	
	Any finite index subgroup $V\leq Z_p$ is $P$-invariant, and  $[Z_p:C_{Z_p}(P)+V]$ is $1$ or $[Z_p:V]$ according to whether $C_{Z}(P)$ is $Z$ or $0$, which in turn, by Lemma \ref{main lemma}(1), is translated into whether $F\subset SL(T)$ or not. Thus, $[Z_p:C_{Z_p}(P)+V]=[Z:V]^\epsilon$, where $\epsilon$ is as in the proposition.    
	 	 Next, by Lemma \ref{main lemma}(3), $|\on{Hom}_{\Lambda_p}(T_p,Z_p/V)|=[Z_p:V]^\eta$. 
	 	 Therefore, by Proposition \ref{formula for the partial subgroup zeta function},  
	 	   \begin{align*}
	 	 	\zeta_{G_p}^{G_p,\s}(s)&=\sum_{\substack{U\leq_{\Lambda_p} T_p,\ V\leq Z_p\\ [U,U]\subseteq V}} [T_p:U]^{-s}[T_p:C_{T_p}(P)+U][Z_p:V]^{-s+\epsilon+\eta}\\
	 	 	&=\sum_{\substack{U\leq_{\Lambda_p} T_p}} [T_p:U]^{-s}[T_p:C_{T_p}(P)+U]\left(\sum_{V\leq Z_p} [Z_p:V]^{-s+\epsilon+\eta}-\sum_{\substack{V\leq p[U,U]}}[Z_p:V]^{-s+\epsilon+\eta}\right)\\
	 	 	&=\sum_{\substack{U\leq_{\Lambda_p} T_p}} [T_p:U]^{-s}[T_p:C_{T_p}(P)+U]\left(\zeta_p(s-\epsilon-\eta)-[Z_p:p[U,U]]^{-s+\epsilon+\eta}\zeta_p(s-\epsilon-\eta) \right)\\
	 	 	&=\zeta_p(s-\epsilon-\eta)\left(\sum_{\substack{U\leq_{\Lambda_p} T_p}} [T_p:U]^{-s}[T_p:C_{T_p}(P)+U]\left(1-p^{-s+\epsilon+\eta}|[Z:[N,N]]|_p^{s-\eta-\epsilon}[T_p:U]^{-s+\epsilon+\eta} \right)\right),
	 	 \end{align*}
 	 where, in the last equality, we used Lemma \ref{index of the commutator}.
To complete the proof we have to show that $\sum_{\substack{U\leq_{\Lambda_p} T_p}} [T_p:U]^{-s}[T_p:C_{T_p}(F)+U]=\zeta_{E_p}^{E_p,\s}(s)$. However, this follows from  \cite[Proposition 2.3 and Lemma 2.4]{dSMS}.
	\end{proof}

\begin{corollary}\label{main corollary counting subgroups} 
 If $p\nmid |P|$, then $\zeta_{G_p}^{G_p}(s)$ is given as in Table \ref{Table with zeta_G^G}, where $k=[Z:[N,N]]$ and
  for $d\in\{3,4,6\}$,
	 $\chi_d:\mathbb{N}\to\mathbb{C}$ is the extended residue class character, $\displaystyle 
	\chi_d(n)=\left\{ \begin{array}{ll}
		1&\mbox{if } n\equiv 1\mod d\\
		-1&\mbox{if } n\equiv -1\mod d\\
		0&\mbox{otherwise},
	\end{array}\right. $,
	and  $L_p(s,\chi_d)=(1-\chi_d(p)p^{-s})^{-1}$ is the local factor at $p$ of the Dirichlet $L$-function of $\chi_d$, $\displaystyle 
	L(s,\chi_d)=\sum_{n=1}^\infty \chi_d(n){n^{-s}}.$
	
	Therefore, if  $p\nmid |P|[Z:[N,N]]$, then $\zeta_{G_p}^{G_p,\s}(s)$ satisfies the functional equation
	\begin{align*}
		&\zeta_{G_p}^{G_p,\s}(s)|_{p\to p^{-1}}=(-1)^3p^{-3s+3}\zeta_{G_p}^{G_p,\s}(s)&&\mbox{if } F\cong C_1,C_2,D_2;\\
		&\zeta_{G_p}^{G_p,\s}(s)|_{p\to p^{-1}}=(-1)^3p^{-3s+2}\chi_d(p)\zeta_{G_p}^{G_p,\s}(s)&&\mbox{if } F\cong C_d,\ d\in\{3,4,6\};\\
		&\zeta_{G_p}^{G_p,\s}(s)|_{p\to p^{-1}}=p^{-3s+3}\zeta_{G_p}^{G_p,\s}(s)&&\mbox{if } F\cong D_d,\ d\in\{3,4,6\}.
	\end{align*}

\begin{small} 
	\begin{table}[h]
	\begin{tabular}{|c|c|cc|}
		\hline 
		&$F$ &$\zeta_{G_p}^{G_p,\s}(s)$ &  \\
		\hline 
		& $C_1$& $\begin{array}{c} \\	\zeta_p(s-2)\left(\zeta_p(s)\zeta_p(s-1)-p^{-s+2}|k|_p^{s-2}\zeta_p(2s-2)\zeta_p(2s-3) \right)\\ \\
			\dfrac{\zeta_p(s)\zeta_p(s-1)\zeta_p(2s-2)\zeta_p(2s-3)}{\zeta_p(3s-3)} 
		\end{array}$ & $\begin{array}{l}p\nmid |P|\\ \\ p\nmid |P|k\end{array}$ \\
		\cline{2-4}
		$F\subset SL_2(T)$	& $C_2$& $\begin{array}{c} \\
			\zeta_p(s)\left(\zeta_p(s-1)\zeta_p(s-2)-p^{-s}|k|_p^{s}\zeta_p(2s-1)\zeta_p(2s-2)\right) \\ \\
			\dfrac{\zeta_p(s-1)\zeta_p(s-2)\zeta_p(2s-1)\zeta_p(2s-2)}{\zeta_p(3s-3)}
		\end{array}$ & $\begin{array}{l}p\nmid |P|\\ \\ p\nmid |P|k\end{array}$  \\
		\cline{2-4}
		& $\begin{array}{c}C_d\\ d=3,4,6\end{array}$& 	$\begin{array}{c}\\ 
			\zeta_p(s)\left(L_p(s-1,\chi_d)\zeta_p(s-1)-p^{-s}|k|_p^{s}L_p(2s-1,\chi_d)\zeta_p(2s-1)\right)\\ \\
			\dfrac{\zeta_p(s-1)\zeta_p(2s-1)L_p(s-1,\chi_d)L_p(2s-1,\chi_d)}{L_p(3s-2,\chi_d)}	
		\end{array}$&  $\begin{array}{l}p\nmid |P|\\ \\ p\nmid |P|k\end{array}$ \\
		\hline 
		& $D_1$& $\begin{array}{c}\\ 
			\zeta_p(s-2)\left(\zeta_p(s)\zeta_p(s-1)-p^{-s+2}|k|_p^{s-2}\zeta_p(2s-2)\zeta_p(2s-3) \right)\\ \\
			\dfrac{\zeta_p(s)\zeta_p(s-1)\zeta_p(2s-2)\zeta_p(2s-3)}{\zeta_p(3s-3)}
		\end{array}$ & $\begin{array}{l}p\nmid |P|\\ \\ p\nmid |P|k\end{array}$ \\
		\cline{2-4}
		$F\not\subset SL_2(T)$	& $D_2$& $\begin{array}{c} \\
			\zeta_p(s-1)\left(\zeta_p(s-1)^2-p^{-s+1}|k|_p^{s-1}\zeta_p(2s-2)^2\right) \\ \\
			\dfrac{\zeta_p(s-1)^2\zeta_p(2s-2)^2}{\zeta_p(3s-3)}
		\end{array}$ & $\begin{array}{l}p\nmid |P|\\ \\ p\nmid |P|k\end{array}$  \\
		\cline{2-4}
		& $\begin{array}{c}D_d\\ d=3,4,6\end{array}$& 	$\begin{array}{c}\\ 
			\zeta_p(s-1)\left(\zeta_p(2s-2)-p^{-s+1}|k|_p^{s-1}\zeta_p(4s-4)\right)\\
			\\
			\dfrac{\zeta_p(2s-2)\zeta_p(3s-3)\zeta_p(4s-4)}{\zeta_p(6s-6)}
		\end{array}$&  $\begin{array}{l}p\nmid |P|\\ \\ p\nmid |P|k\end{array}$ \\
		\hline 	
	\end{tabular}
\caption{Local factors of $\zeta_G^{G,\s}(s)$ at $p\nmid |P|$ ($F:=\on{Im}(P\to GL(T))$ and $k:=[Z:[N,N]]$).}
\label{Table with zeta_G^G}
\end{table} 
\end{small} 
\end{corollary}
\begin{proof}
	The formula for $\zeta_{G_p}^{G_p,\s}(s)$ in each case follows from Theorem  \ref{main theorem for zeta_G^G} and uses the formula for $\zeta_{E}^{E,\s}(s)$ obtained in \cite[Chap.\ 5]{Mc} as a step in the computation of the subgroup zeta functions of the plane crystallographic group $E$ (see also \cite[Section 4.1]{dSMS}). The computation of $\eta$ was done in Lemma \ref{main lemma}(2). If in addition $p\nmid |P|[Z:[N,N]]$, then $|[Z:[N,N]]|_p=1$ and the simplification of the formula is straightforward. The functional equation follows by inspection of formula. 
	\end{proof} 

\begin{corollary}\label{Meromorphic continuation partial subgroup zeta function}
 $\zeta_{G}^{G,\s}(s)$ has abscissa of convergence in the set $\{\frac{3}{2},2,3\}$ and admits a meromorphic continuation to the whole plane. The same holds for $\zeta_G^\s(s)$. 
\end{corollary}
\begin{proof}
	According to Table \ref{Table with zeta_G^G}, there is a Dirichlet series $Z(s)=\prod_{p}Z_p(s)$ such that: (a) it has abscissa of convergence in the set $\{\frac{3}{2},2,3\}$; (b) it admits a meromorphic continuation to the whole plane; (c) for any finite set of primes, say $S$, $Z(s)$ and $\prod_{p\notin S}Z_p(s)$ have the same abscissa of convergence; and
	(d) $Z(s)$ coincides with $\zeta_{G}^{G,\s}(s)$ except for a finite number of local factors.
	Since an exceptional local factor $\zeta_{G_p}^{G_p,\s}(s)$ is a rational function of $p^{-s}$ with rational coefficients \cite{dS}, it follows from (b) and (d) that $\zeta_G^{G,s}(s)$ also has a meromophic continuation to the whole plane.  
	Moreover, $\zeta_G^{G,\s}(s)$ and $Z(s)$ have the same abscissa of convergence. Indeed, by (c) and (d), it is enough to show that the abscissa of convergence of each local factor of $\zeta_{G}^{G,\s}(s)$ is strictly less than the abscissa of convergence of $\zeta_{G}^{G,\s}(s)$. Now, this follows from the fact that $\zeta_{G}^{G,\s}(s)$ is an Euler product of cone integrals
 \cite[Theorem 1]{Su}. This important property about cone integrals was established in \cite[Section 4]{dSG}, and it was a key point in obtaining analytic properties of global zeta functions. 
 
 To prove the last assertion of the corollary, we apply the first part to each of the partial zeta functions of $\zeta_G^\s(s)$.
\end{proof}

\begin{remark}
	If $F=C_1$ in Table \ref{Table with zeta_G^G}, we recover \cite[Proposition 8.1]{GSS}. See also \cite{KlopschVoll2007} for a generalization.
\end{remark}

\subsection{Counting normal subgroups}
We begin with some preliminaries from elementary group theory.
Let $G$ be for the moment any group, and let $N\lhd G$ be a normal subgroup. We define inductively a series of normal subgroups $\gamma_1(G,N)\supseteq\gamma_2(G,N)\supseteq\ldots$ by setting $$\gamma_1(G,N):=N\ \mbox{and}\ \gamma_{i}(G,N):=[\gamma_{i-1}(G,N),G],\  \mbox{for}\  i\geq 2.$$ 
\begin{lemma}\label{Hall nilpotent subgroups} Assume that $G$ is finite, and that $N$ is a normal Hall subgroup. If $\gamma_{c+1}(G,N)=1$ for some $c$, then $N$ has a unique complement, say $C$. In addition, $G=N\times C$ and  $\gamma_{c+1}(G)\cap N=\{1\}$.
\end{lemma}
\begin{proof}
	By the Schur-Zassenhaus theorem, there is at least one complement of $N$ in $G$, and they are all conjugated. 
	We prove that there is only one by induction on $c$.
	
	 Assume that $c=1$, and let $C$ be a complement of $N$ in $G$. Since $[N,C]\subseteq [N,G]=\gamma_2(G,N)=1$, it follows that $C$ is normal in $NC=G$, and therefore, $C$ is the unique complement of $N$.

	Assume now that $c>1$. By the inductive hypothesis, there is a unique complement, say $C'/\gamma_c(G,N)$, of $N/\gamma_c(G,N)$ in $G/\gamma_c(G,N)$. 
	Now, given a complement $C$ of $N$ in $G$, clearly $(C\gamma_c(G,N))/\gamma_c(G,N)$ is a complement of $N/\gamma_c(G,N)$ in $G/\gamma_c(G,N)$, thus $C\gamma_c(G,N)=C'$. This implies that $C$ is a complement of $\gamma_c(G,N)$ in $C'$. However, $[\gamma_c(G,N),C']\subseteq\gamma_{c+1}(G,N)=1$, so by the case $c=1$, there is only one possibility for $C$. This completes the induction. 
	
	The final part of the lemma is clear.
\end{proof}
\begin{lemma}\label{Gamma_c(G) is included in A}
	Assume that $N$ is nilpotent, say of class $c$, and let $A\lhd G$ such that $AN=G$. Then $\gamma_{c+1}(G)\subseteq A$; in particular,  $\gamma_{c+1}(G,N)\subseteq A\cap N$.
\end{lemma}
\begin{proof}
	Note that $B:=A\cap N$ is also normal in $G$, and there is an identification  $G/B=N/B\times A/B$. It follows that $\gamma_{c+1}(G/B)=\gamma_{c+1}(N/B)\times\gamma_{c+1}(A/B)=1\times\gamma_{c+1}(A/B)\subseteq A/B$. Since $\gamma_{c+1}(G/B)=(\gamma_{c+1}(G)B)/B$, we deduce that $\gamma_{c+1}(G)\subset A$. In particular, $\gamma_{c+1}(G,N)\subseteq \gamma_{c+1}(G)\cap N\subseteq A\cap N$.
	\end{proof} 

\smallskip 

We return to the setting introduced at the beginning of the section except that we do not assume that the nilpotency class is 2 yet. We fix an intermediate normal subgroup $N\leq H\lhd G$ and consider the series $\zeta_{G_p}^{H_p,\n}(s)$. 
  Note that $\gamma_i(H,N)$ is normal in $G$ for all $i$, and that $(\gamma_i(H,N))_p=\gamma_i(H_p,N_p)$ and $(\gamma_i(H))_p=\gamma_i(H_p)$ (cf. \cite[Theorem 1.4]{NS}).
 
\begin{proposition}\label{simplification for normal extension zeta function}
	 Let $c$ be the nilpotency class of $N$. Set $G':=G/\gamma_{c+1}(H,N)$ and $N':=N/\gamma_{c+1}(H,N)$. If $p\nmid |H/N|$, then $\zeta_{G_p}^{H_p,\n}(s)=\zeta_{G_p'}^{N_p',\n}(s)$. In particular, if $N_p'\cong\mb{Z}_p/k\mb{Z}_p$ for some $k\in\mb{Z}$, then $\zeta_{G_p}^{H_p,\n}(s)=\frac{1-p^{-s}|k|_p^s}{1-p^{-s}}$.
\end{proposition}
\begin{proof}
Let $H'=H/\gamma_{c+1}(H,N)$. By Lemma \ref{Gamma_c(G) is included in A}, we have $\zeta_{G_p}^{H_p,\n}(s)=\zeta_{G'_p}^{H'_p,\n}(s)$, where the series on the right is computed with respect to $N_p'\lhd G_p'$ (this holds for all $p$). Assume now that $p\nmid |H/N|$.
Given $A'\leq H_p'$ of finite index and normal in $G_p'$ such that $A'N'_p=H_p'$, the intersection $B':=N_p'\cap A'$ is normal in $G_p'$ and $[N_p':B']=[H_p':A']$. Conversely, given $B'\leq N_p'$ of finite index and normal in $G_p'$, we have $\gamma_{c+1}(H_p'/B',N_p'/B')=(\gamma_{c+1}(H_p',N_p') B')/B'=1$ since $\gamma_{c+1}(H_p',N_p')=1$. Therefore, we can apply Lemma \ref{Hall nilpotent subgroups} to $N'_p/B'\lhd H'_p/B'$. It follows that
	 there is a unique $A'\leq H_p'$ such that $A'N_p=H_p'$ and $A'\cap N_p'=B'$. By the uniqueness, $A'/B'$ is normal in $G_p'/B'$, and hence $A'$ is normal in $G_p'$.
	 We deduce that $\zeta_{G'_p}^{H'_p,\n }(s)=\zeta_{G_p'}^{N_p',\n}(s)$. This proves the first part of the lemma.
	 
	  If in addition $N_p'\cong\mb{Z}_p/k\mb{Z}_p$, then every finite index subgroup of $N_p'$ is characteristic; therefore, $\zeta_{G_p'}^{N_p',\n}(s)=\zeta_{\mb{Z}_p/k\mb{Z}_p}^\s(s)$, which is clearly equal to $\frac{1-p^{-s}|k|_p^s}{1-p^{-s}}$.
	\end{proof}
\begin{corollary}\label{the abscissa of convergence of the normal zeta function}
	 $\zeta_G^\n(s)$ and $\zeta_G^{N,\n}(s)$ have the same abscissa of convergence.
\end{corollary}
\begin{proof}
We have to show that for each intermediate normal subgroup $N\leq H\lhd G$, the abscissa of convergence of $\zeta_G^{H,\n}(s)$ is bounded by that of $\zeta_G^{N,\n}(s)$. 
 We set $G':=G/\sqrt{\gamma_{c+1}(H,N)}$ and $N':=N/\sqrt{\gamma_{c+1}(H,N)}$, where $\sqrt{\gamma_{c+1}(H,N)}:=\{x\in N: x^n\in \gamma_{c+1}(H,N)\ \mbox{for some } n\in\mathbb{N}\}$. By Proposition \ref{simplification for normal extension zeta function}, $\zeta_{G}^{H,\n}(s)$ and $\zeta_{G'}^{N',\n}(s)$ have the same but a finite number of local factors.
 	By \cite[Theorem 1]{Su}, both series are Euler products of cone integrals. As explained in the proof of Corollary \ref{Meromorphic continuation partial subgroup zeta function}, this suffices to ensure that $\zeta_{G}^{H,\n}(s)$ and $\zeta_{G'}^{N',\n}(s)$ have the same abscissa of convergence. Finally, it is clear that the abscissa of convergence of $\zeta_{G'}^{N',\n}(s)$ is bounded by that of $\zeta_G^{N,\n}(s)$. 
\end{proof}

\begin{remark}
	We claim that if $p\nmid |H/N|$, then $\gamma_{c+1}(H_p,N_p)=N_p\cap\gamma_{c+1}(H_p)$; in particular, if $H_p/N_p$ is nilpotent of class $\leq c$, then $\gamma_{c+1}(H_p,N_p)=\gamma_{c+1}(H_p)$. Indeed,
	 let $B\lhd G_p$ be a finite index normal subgroup such that $\gamma_{c+1}(H_p,N_p)\subseteq B\subseteq N_p$. Since $N_p/B$ is a normal Hall subgroup of $H_p/B$, and since $\gamma_{c+1}(H_p/B,N_p/B)=(\gamma_{c+1}(H_p,N_p) B)/B=1$, we obtain by Lemma \ref{Hall nilpotent subgroups} that $\gamma_{c+1}(H_p/B)\cap N_p/B=1$. It follows that $(\gamma_{c+1}(H_p)B)\cap N_p\subseteq B$, and hence $\gamma_{c+1}(H_p)\cap N_p\subseteq B$. Since $N_p/\gamma_{c+1}(H_p,N_p)$ is residually finite, the obvious inclusion $\gamma_{c+1}(H_p,N_p)\subseteq  \gamma_{c+1}(H_p)\cap N_p$ must be an equality.
\end{remark}

We now return to the case when $N$ is a $\mf{T}_2$-group. Recall the notation introduced at the beginning of the section.
Given $V\leq Z_p$ of finite index, we denote by $X(V)\leq T_p$ the largest subgroup such that $[T_p,X(V)]\subseteq V$. It has finite index in $T_p$, and if $V$ is in addition $P$-invariant, then  $X(V)$ is also $P$-invariant. 
 Proposition \ref{simplification for normal extension zeta function} allows us to focus only on the first partial zeta function  $ \zeta_G^{N,\n}(s)=\sum_{\substack{B\lhd G: B\subseteq N}}[N:B]^{-s}$, the other ones being of the form $\zeta_{N'}^{G',\n}(s)$ except for a finite number of local factors, where $G'$ is a quotient of $G$.
\begin{proposition}\label{formula for the partial normal zeta function}
	Fix a prime $p$. Assume that all  finite index $\Lambda_p$-submodules of $T_p$ are isomorphic. Assume also that for every $V\leq_{\Lambda_p} Z_p$ and $U\leq_{\Lambda_p}T_p$ of finite index such that $U\subseteq X(V)$, there exists $B\lhd G_p$ of finite index such that $B\subseteq N_p$, $(BZ_p)/Z_p=U$ and $B\cap Z_p=V$. Then it holds that  	
	$$
	\zeta_{G_p}^{N_p,\n}(s)=\zeta_{E_p}^{ T_p,\n}(s) \left(
	\sum_{V\leq_{\Lambda_p} Z_p} [T_p:X(V)]^{-s}[Z_p:V]^{-s}|\on{Hom}_{\Lambda_p}(T_p,Z_p/V)|\right).
	$$
	The above assumptions hold if $p\nmid |P|$.
\end{proposition}
\begin{proof}
	The proof is similar to that of \cite[Lemma 6.1]{GSS}, so we omit it. We just point out that given $U\leq_{\Lambda_p} T_p$ and $V\leq_{\Lambda_p} Z_p$ of finite index with $U\subseteq X(V)$, the number of normal subgroups $B\lhd G_p$ such that $(BZ_p)/Z_p=U$ and $B\cap Z_p=V$ is equal to $\on{Hom}_{\Lambda_p}(U,Z_p/V)$. This uses the second assumption and Lemma \ref{number of F-invariant complements}. In addition, by the first assumption, $\on{Hom}_{\Lambda_p}(U,_p/V)=\on{Hom}_{\Lambda_p}(T_p/V)$. 
%
%
 
The first assumption is satisfied when $p\nmid |P|$ since in this case $\Lambda_p$ is a maximal order of $\mb{Q}_p[P$], and the second assumption is satisfied by Lemma \ref{number of F-invariant complements}.
\end{proof}

We now specialize to the case of Hirsch length 3.
\begin{lemma}\label{index of the centralizer of V}
	Assume that $N$ is a $\mf{T}_2$-group of Hirsch length 3. Given $V\leq Z_p$ of finite index, the following holds:
	\begin{enumerate}
		\item If $[T_p,T_p]\subseteq V$, then $[T_p:X(V)]=1$.
		\item If $V\subseteq [T_p,T_p]$, then $[T_p:X(V)]=|[Z:[N,N]]|_p^2[Z_p:V]^2$.
	\end{enumerate}
\end{lemma}
\begin{proof}
The first assertion is trivial since in this case $X(V)=T_p$. Assume now that $V\subseteq [T_p,T_p]$.	
Let $\{x,y\}$ be a basis for $T_p$.	Given a $\mb{Z}_p$-submodule $X\leq T_p$ of finite index, say generated by $ax+by$ and $cx+dy$, we have $[X,T_p]=\langle a[x,y], b[x,y],c[x,y],d[x,y]\rangle$. Thus, $[X,T_p]\subseteq V$ if and only if $a[x,y], b[x,y], c[x,y], d[x,y]\in V$. Since $[x,y]$ generates $[T_p,T_p]$, the previous condition holds if and only if $a,b,c,d\in \frac{[Z_p:V]}{[Z_p:[T_p,T_p]]}\mb{Z}_p$, which in turn is equivalent to saying that $X\subseteq \frac{[Z_p:V]}{[Z_p:[T_p,T_p]]} T_p$. Thus, $X(V)=\frac{[Z_p:V]}{[Z_p:[T_p,T_p]]} T_p$. The index of this subgroup is $[Z_p:V]^2|[Z_p:[T_p,T_p]]|_p^2=[Z_p:V]^2|[Z:[N,N]]|_p^2$. 
\end{proof}

\begin{theorem}\label{main theorem for zeta_N^G}
	Assume that $N$ is a $\mf{T}_2$-group of Hirsch length 3. Set $F:=\on{Im}(P\to GL(T))$ and $E:=T\rtimes F$, which is a plane crystallographic group. If $p\nmid |P|$, then
	$$
	\zeta_{G_p}^{N_p,\n}(s)=\zeta_{E_p}^{T_p,\n}(s) \left(
	\frac{1-|[Z:[N,N]]|_p^{s-\eta}}{1-p^{-s+\eta}}+ |[Z:[N,N]]|_p^{s-\eta}\zeta_p(3s-\eta)\right),
	$$
	where $\eta:=\on{rank}_{\mb{Z}}(C_T(P))$. 
	In particular, if $p\nmid |P|[Z:[N,N]]$, then
		$$
	\zeta_{G_p}^{N_p,\n}(s)=\zeta_{E_p}^{T_p,\n}(s) \zeta_p(3s-\eta).
	$$ 
\end{theorem}
\begin{proof}
By Lemma \ref{main lemma}, $|\on{Hom}_{\Lambda_p}(T_p,Z_p/V)|=[Z_p:V]^\eta$ for any $V\leq_{\Lambda_p} Z_p$, where $\eta=\on{rank}_{\mb{Z}}(C_T(P))$.
	It follows from Proposition \ref{formula for the partial normal zeta function} and Lemma \ref{index of the centralizer of V} that
		\begin{align*} 
	\zeta_{G_p}^{ N_p,\n}(s)&=\zeta_{E_p}^{T_p,\n}(s) \left(	\sum_{\substack{ [T_p,T_p]\subsetneqq V\leq Z_p}} [Z_p:V]^{-s+\eta}+
	\sum_{V\leq [T_p,T_p]} |[Z:[N,N]]|_p^{-2s}[Z_p:V]^{-3s+\eta}\right)\\
	&=\zeta_{E_p}^{T_p,\n}(s) \left(\sum_{k=0}^{v_p([T_p,T_p])-1} p^{(-s+\eta)k}+
	\sum_{V\leq [T_p,T_p]} |[Z:[N,N]]|_p^{2s}[Z_p:[T_p,T_p]]^{-3s+\eta}[[T_p,T_p]:V]^{-3s+\eta}\right)\\
	&=\zeta_{E_p}^{T_p,\n}(s) \left(\frac{1-|[Z:[N,N]]|_p^{s-\eta}}{1-p^{-s+\eta}}+ |[Z:[N,N]]|_p^{s-\eta}\zeta_p(3s-\eta)\right).
	\end{align*}
If in addition $p\nmid [Z:[N,N]]$, then $|[Z:[N,N]]|_p=1$, and the proof follows.
\end{proof}
\begin{corollary}\label{main corollary normal zeta functions}
If $p\nmid |P|$, then $\zeta_{G_p}^{N_p,\n}(s)$ is given as in Table \ref{Table with zeta_G^N}, where $k:=[Z:[N,N]]$. Therefore, if $p\nmid |P|[Z:[N,N]]$, then $\zeta_{G_p}^{N_p,\n}(s)$ satisfies the following functional equation: 	\begin{align*}
 	&\zeta_{G_p}^{N_p,\n}(s)|_{p\to p^{-1}}=(-1)^3p^{-5s+3}\zeta_{G_p}^{N_p,\n}(s)&&\mbox{if } F=C_1,\\
 	&\zeta_{G_p}^{N_p,\n}(s)|_{p\to p^{-1}}=(-1)^3 p^{-5s+1}\zeta_{G_p}^{N_p,\n}(s)&&\mbox{if } F\cong C_2= D_1, \\
 	&\zeta_{G_p}^{N_p,\n}(s)|_{p\to p^{-1}}=(-1)^3\chi_d(p) p^{-5s}\zeta_{G_p}^{N_p,\n}(s)&&\mbox{if } F\cong C_d,\ d\in\{3,4,6\},\\
 	&\zeta_{G_p}^{N_p,\n}(s)|_{p\to p^{-1}}=(-1)^3p^{-5s}\zeta_{G_p}^{N_p,\n}(s) &&\mbox{if  } F\cong D_2,\\
 	&\zeta_{G_p}^{N_p,\n}(s)|_{p\to p^{-1}}=p^{-5s}\zeta_{G_p}^{N_p,\n}(s) &&\mbox{if  } F\cong D_d,\ d\in\{3,4,6\}.
 \end{align*}
	\begin{small} 
		\begin{table}[h]
			\begin{tabular}{|c|c|cc|}
				\hline 
				&$F$ &$\zeta_{G_p}^{N_p,\n}(s)$ &  \\
				\hline 
				& $C_1$& $\begin{array}{c} \\	\ds \zeta_p(s)\zeta_p(s-1) \left(
					\frac{1-|k|_p^{s-2}}{1-p^{-s+2}}+ |k|_p^{s-2}\zeta_p(3s-2)\right)\\ \\
					\ds\zeta_p(s)\zeta_p(s-1)\zeta_p(3s-2) 
				\end{array}$ & $\begin{array}{l}\\ p\nmid |P|\\ \\ p\nmid |P|k\end{array}$ \\
				\cline{2-4}
				$F\subset SL_2(T)$	& $C_2$& $\begin{array}{c} \\
					\ds \zeta_p(s)\zeta_p(s-1) \left(
					\frac{1-|k|_p^{s}}{1-p^{-s}}+ |k|_p^{s}\zeta_p(3s)\right) \\ \\
					\ds \zeta_p(s)\zeta_p(s-1) \zeta_p(3s)
				\end{array}$ & $\begin{array}{l}\\ p\nmid |P|\\ \\ p\nmid |P|k\end{array}$  \\
				\cline{2-4}
				& $\begin{array}{c}C_d\\ d=3,4,6\end{array}$& 	$\begin{array}{c}\\ 
					\ds	\zeta_p(s)L_p(s,\chi_d) \left(
					\frac{1-|k|_p^{s}}{1-p^{-s}}+ |k|_p^{s}\zeta_p(3s)\right)\\ \\
					\ds \zeta_p(s) L_p(s,\chi_d) \zeta_p(3s)
				\end{array}$&  $\begin{array}{l}\\ p\nmid |P|\\ \\ p\nmid |P|k\end{array}$ \\
				\hline 
				& $D_1$& $\begin{array}{c}\\ 
					\ds \zeta_p(s)^2 \left(
					\frac{1-|k|_p^{s-1}}{1-p^{-s+1}}+ |k|_p^{s-1}\zeta_p(3s-1)\right)\\ \\
					\ds \zeta_p(s)^2 \zeta_p(3s-1)
				\end{array}$ & $\begin{array}{l}\\ p\nmid |P|\\ \\ p\nmid |P|k\end{array}$ \\
				\cline{2-4}
				$F\not\subset SL_2(T)$	& $D_2$& $\begin{array}{c} \\
					\ds \zeta_p(s)^2 \left(
					\frac{1-|k|_p^{s}}{1-p^{-s}}+ |k|_p^{s}\zeta_p(3s)\right) \\ \\
					\ds \zeta_p(s)^2  \zeta_p(3s)
				\end{array}$ & $\begin{array}{l}\\ p\nmid |P|\\ \\ p\nmid |P|k\end{array}$  \\
				\cline{2-4}
				& $\begin{array}{c}D_d\\ d=3,4,6\end{array}$& 	$\begin{array}{c}\\ 
					\ds \zeta_p(2s) \left(
					\frac{1-|k|_p^{s}}{1-p^{-s}}+ |k|_p^{s}\zeta_p(3s)\right)\\
					\\
					\ds\zeta_p(2s)\zeta_p(3s)
				\end{array}$&  $\begin{array}{l}\\ p\nmid |P|\\ \\ p\nmid |P|k\end{array}$ \\
				\hline 	
			\end{tabular}
			\caption{Local factors of $\zeta_G^{N,\n}(s)$ at $p\nmid |P|$ ($F:=\on{Im}(P\to GL(T))$ and $k:=[Z:[N,N]]$).}
			\label{Table with zeta_G^N}
		\end{table} 
	\end{small} 
\end{corollary}
\begin{proof}
	This follows from Theorem \ref{formula for the partial normal zeta function} and uses the formula for $\zeta_{E}^{T,\n}(s)$ obtained in \cite[Chap. 6]{Mc} as a step in the computation of the normal zeta function of the plane crystallographic group $E$ (see also \cite[Section 4.2]{dSMS}). The calculation of $\eta$ was done in Lemma \ref{main lemma}. The local functional equation follows by inspection of the formula.
\end{proof}

\begin{corollary}\label{meromorphic contunuation normal zeta function}
	$\zeta_{G}^{N,\n}(s)$ has abscissa of convergence in the set $\{\frac{1}{2},1,2\}$ and admits a meromorphic continuation to the whole plane. The same holds for $\zeta_G^\n(s)$.
\end{corollary}
\begin{proof}
	The assertion about $\zeta_{G}^{N,\n}(s)$ follows by inspection of Table \ref{Table with zeta_G^N} and the same argument used in the proof of Corollary \ref{Meromorphic continuation partial subgroup zeta function}. Next, by Corollary \ref{the abscissa of convergence of the normal zeta function},  $\zeta_G^\n(s)$ has the same abscissa of convergence as $\zeta_G^{N,\n}(s)$. To complete the proof, it is enough to show that for each intermediate normal subgroup $N\leq H\lhd G$ different from $N$, the series $\zeta_{G}^{H,\n}(s)$ admits a meromorphic continuation to the whole plane. We use the notation from the proof of Corollary \ref{the abscissa of convergence of the normal zeta function} with $c=2$. Since the local factor of $\zeta_{G}^{H,\n}(s)$ at a prime $p$ is a rational function in $p^{-s}$ by \cite{dS}, it is enough to show that $\zeta_{G'}^{N',\n}(s)$ admits a meromorphic continuation to the whole plane. If $N'=N$, then we are in the same situation as in the case of $\zeta_G^N(s)$, so we are done. If $N'\neq N$, then either $N'\cong\mb{Z}^2$, in which case the result follows from the formulae of the normal zeta functions of the plane crystallographic groups given in \cite[Chap.\ 6]{Mc}, or else $N'\cong\mb{Z}$ or $N'=0$, in which case the result is trivial.
\end{proof}

\subsection{Remark on the local functional equations}
One natural problem is to decide whether local functional equations hold for $\zeta_G^{G,\s}(s)$ and $\zeta_G^{N,\n}(s)$ when the Hirsch length of $N$ is higher than 3. When $N$ is abelian, this follows from the explicit formulae obtained in \cite[Section 2]{dSMS}. In contrast to the situation of zeta functions of nilpotent groups, 
 here the local functional equations are not uniform in $p$. They depend on how $p$ ramifies in certain number fields that arise in the decomposition of $\mb{Q}[G/N]$ into simple algebras.
 Propositions \ref{formula for the partial subgroup zeta function} and \ref{formula for the partial normal zeta function} enable us to linearize the problem when the nilpotency class is 2 (indeed, the linearization for $\zeta_G^{N,\n}(s)$ holds more generally by the Mal'cev correspondence). We will not go into details.  
Briefly, the problem can be stated as follows. Let $L$ be a 2-step nilpotent Lie ring additively isomorphic to $\mb{Z}^h$, and let $Z$ denote its center. Let $P$ be a finite group acting on $L$ by Lie ring automorphisms.
We consider the Dirichlet series
\begin{align*}
\zeta_{P\curvearrowright L}^{\s}(s)=\sum [L:B]^{-s}|\on{Der}(P,L/B)|\quad\mbox{and}\quad 	\zeta_{P\curvearrowright L}^{\n}(s)=\sum [L:B]^{-s}, 
\end{align*}
where the first sum runs over the $P$-invariant finite index subrings, and the second one only over the $P$-invarian finite index ideals. Corollaries \ref{main corollary counting subgroups} and \ref{main corollary normal zeta functions} imply that these series satisfy local functional equations when $L$ has additive rank 3. In analogy with the results of \cite{V1}, we might ask:
\begin{question}
	Do $\zeta_{P\curvearrowright L}^{\s}(s)$ and $\zeta_{P\curvearrowright L}^{\n}(s)$ satisfy local functional equations when $L$ is a 2-step nilpotent Lie ring of rank $>3$? Does $\zeta_{P\curvearrowright L}^{\s}(s)$ satisfy local functional equations for arbitrary rings  (not necessarily nilpotent) that are additively isomorphic to some $\mb{Z}^h$?
\end{question}

 We return to the case when $L$ is a 2-step nilpotent Lie ring additively isomorphic to $\mb{Z}^3$. It is not difficult to check, say by inspection of Table \ref{Table with zeta_G^G} and the results of \cite[Section 2]{dSMS}, that if we forget the structure of Lie ring on $L$ (obtaining $\mb{Z}^3$) and consider just $\zeta_{P\curvearrowright \mb{Z}^3}^{\s}(s)$, then $\zeta_{P\curvearrowright L}^{\s}(s)$ and $\zeta_{P\curvearrowright \mb{Z}^3}^{\s}(s)$ satisfy the same local functional equations. By this, and again in analogy with (N4), we can also ask:
\begin{question}
	Let $L$ be a ring additively isomorphic to $\mb{Z}^h$. Let $P$ be a finite group acting on $L$.  If $\zeta_{P\curvearrowright L}^{\s}(s)$  satisfies local functional equations, are these equations the same as those satisfied by $\zeta_{P\curvearrowright \mb{Z}^h}^{\s}(s)$?   
\end{question}
In a forthcoming paper, we show that Question 1 and Question 2 have positive answer for various nilpotent Lie rings of rank 4.

\section{Local zeta functions as $p$-adic integrals}\label{Section: method of p-adic integration}
\noindent 
Let $N_p$ be a torsion-free finitely generated nilpotent pro-$p$ group and let $G_p$ be a profinite group that includes $N_p$ as an open normal subgroup. Fix also an intermediate normal subgroup $N_p\leq H_p\lhd G_p$.
We review a method, developed in  \cite[Section 2]{GSS} for $\mf{T}$-groups and extended to virtually nilpotent groups in \cite{Su}, to express $\zeta_{H_p}^{H_p,\s}(s)$ and $\zeta_{G_p}^{H_p,\n}(s)$ as $p$-adic integrals.  This method is used in the next sections to calculate local factors at ``bad" primes in some cases. 

\subsection{Expressing $\zeta_{N_p}^\s(s)$ and $\zeta_{N_p}^\n(s)$ as $p$-adic integrals.} 
An additional reference for this part is \cite[Chap. 15]{LS}, where the notation is more adapted to ours.
 Fix a Mal'cev basis $\x=(x_1,\ldots,x_h)$ for $N_p$. By definition, the series of subgroups
$$N_p=\overline{\langle x_1,\ldots,x_h\rangle}\supset \overline{\langle x_2,\ldots,x_h\rangle}\supset\cdots \supset\overline{\langle x_h\rangle}$$
is central, and each $x\in N_p$ determines a unique vector $\a=(a_1,\ldots,a_h)\in\mb{Z}_p^h$ such that  $x=\x^\a:=x_1^{a_1}\cdots x_h^{a_h}.$

Let $\on{T}_h(\mb{Z}_p)$ denote the set of $h\times h$-upper-triangular matrices with entries in $\mb{Z}_p$, and let $\on{T}_h^+(\mb{Z}_p)\subset\on{T}_h(\mb{Z}_p)$ denote the subset of those matrices with non-zero determinant. For $\t\in\on{T}_h(\mb{Z}_p)$, we set
$$B_\t:=\overline{\langle\x^{\t_1},\ldots,\x^{\t_h}\rangle}\leq N_p,$$
where $\t_i$ denotes the $i$-th row of $\t$.
The subgroup $B_\t$ is open if and only if $\t\in\on{T}_h^+(\mb{Z}_p)$.

Given an open subgroup $B\leq N_p$, we say that $\t\in\on{T}_h^+(\mb{Z}_p)$ {\em represents a good basis for $B$} (with respect to $\x$) if $B=B_\t$ and 
$B=\{(\x^{\t_1})^{\lambda_1}\ldots (\x^{\t_h})^{\lambda_h}:\lambda_1,\ldots,\lambda_h\in\mb{Z}_p\}.$
In this case, $(\x^{\t_1},\ldots,\x^{\t_h})$ is a Mal'cev basis for $B$ (also called a {\em good basis for $B$}).
We set 
$$\mc{M}(B):= \{\t\in \on{T}_h(\mb{Z}_p):\t\ \mbox{represents a good basis for } B\}$$
Note that given $\t=(t_{ij})\in\mc{M}(B)$, the value $|t_{ii}|_p$ depends only on $B$ (and on the fixed Mal'cev basis $\x$). Indeed, we have $|t_{ii}|_p^{-1}= [\overline{\langle x_i,\ldots,x_h\rangle}:B\cap \overline{\langle x_i,\ldots,x_h\rangle}]$. We also see from this that
\begin{equation}\label{inverse of the index}
	[N_p:B]^{-1}=\prod_{i=1}^h|t_{ii}|_p.
\end{equation} 
We collect results from \cite[Section 2]{GSS} in the next lemma. The topological group $\on{T}_h(\mb{Z}_p)\cong\mb{Z}_p^{h(h+1)/2}$ has a normalized Haar measure, which we denote by $\mu$.
\begin{lemma}\label{properties of M(B)}
	$\mc{M}(B)$ is an open subset of $\on{T}_h(\mb{Z}_p)$, and for any $\t=(t_{ij})\in\mc{M}(B)$   
	\begin{align}\label{Haar measure of M(B)} \mu(\mc{M}(B))=(1-p^{-1})^h\prod_{i=1}^h|t_{ii}|_p^i.
	\end{align} 
	Therefore, for a complex variable $s$,
	\begin{equation}\label{expression of the index as a p-adic integral} [N:B]^{-s}=\frac{1}{(1-p^{-1})^h}\int_{\mc{M}(B)} \prod_{i=1}^h|t_{ii}|_p^{s-i}d\mu. 
	\end{equation} 
\end{lemma}
\begin{proof}
	The first part of the lemma is proved in \cite[Lemma 2.5]{GSS}; see also \cite[Lemma 15.1.1]{LS}. We recall here the proof of (\ref{expression of the index as a p-adic integral}). We start from the right-hand side and use (\ref{inverse of the index}) and (\ref{Haar measure of M(B)}):
	$$\frac{1}{(1-p^{-1})^h}\int_{\mc{M}(B)} \prod_{i=1}^h|t_{ii}|_p^{s-i}d\mu=\frac{1}{(1-p^{-1})^h}\left( \prod_{i=1}^h|t_{ii}|_p^{s-i}\right) \mu({\mc{M}(B)})=\prod_{i=1}^h|t_{ii}|_p^{s-i}\prod_{i=1}^h|t_{ii}|_p^{i}=[N_p:B]^{-s}.$$
\end{proof}

We now set $\mc{M}_N^\s:=\cup_{B\leq N}\mc{M}(B)$ and $\mc{M}_N^\n:=\cup_{B\lhd N}\mc{M}(B)$, where the unions run only over open subgroups. These are open subsets of $\on{T}_h^+(\mb{Z}_p)$, by Lemma \ref{properties of M(B)}. They both coincide with $\on{T}_h^+(\mb{Z}_p)$ if $N$ is abelian. In the general case, the following descriptions were given in \cite[Lemmas 2.3 and 2.4]{GSS}:
\begin{align}
	\label{condition to represent a good basis} \mc{M}_{N_p}^\s&=\{\t\in\on{T}_h^+(\mb{Z}_p):[\x^{\t_i},\x^{\t_j}]\in \overline{\langle \x^{\t_{j+1}},\ldots,\x^{\t_h}\rangle},\ \ 1\leq i<j\leq h\},\\
	\nonumber  \mc{M}_{N_p}^\n&=\{\t\in\on{T}_h^+(\mb{Z}_p):[x_i,\x^{\t_j}]\in \overline{\langle \x^{\t_{j+1}},\ldots,\x^{\t_h}\rangle},\ \ 1\leq i,j\leq h\}.
\end{align}
As an immediate consequence of formula (\ref{expression of the index as a p-adic integral}), we obtain the following corollary.
\begin{corollary}[Proposition 2.7, \cite{GSS}]\label{zeta functions of N}
	For $*\in\{\leq,\lhd\}$,
	\begin{align*}
		\zeta_{{N}_p}^*(s)=\frac{1}{(1-p^{-1})^{h}}\int_{\mc{M}_{N_p}^*}\prod_{i=1}^h|t_{ii}|_p^{s-i}d\mu.
	\end{align*}
\end{corollary}

\subsection{Expressing $\zeta_{H_p}^{H_p,\s}(s)$ and $\zeta_{G_p}^{H_p,\n}(s)$ as $p$-adic integrals} We begin with the following particular case, which is an immediate consequence of  (\ref{expression of the index as a p-adic integral}).
\begin{corollary}\label{integral expression for zeta_G^N}
 Let $\alpha_1,\ldots,\alpha_r\in G_p$ such that their classes modulo $N_p$ generate $G_p/N_p$. Then
		\begin{align*}
		\zeta_{{G}_p}^{N_p,\n}(s)=\frac{1}{(1-p^{-1})^{h}}\int_{\substack{\t\in\mc{M}_{N_p}^\n\\ {}^{\alpha_j}(\x^{\t_i})\in \overline{\langle \x^{\t_{1}},\ldots,\x^{\t_h}\rangle},\\ i=1,\ldots,h,\ j=1,\ldots,r}}\prod_{i=1}^h|t_{ii}|_p^{s-i}d\mu.
	\end{align*}
\end{corollary}
It remains to consider the case $H_p\neq N_p$. 
We set $F:=H_p/N_p$ and denote by $e$ its identity element.
The following is fixed in the rest of this section.
\begin{enumerate}[-]
	\item A presentation $\langle f_1,\ldots,f_r| R_\lambda(f_1,\ldots,f_r)=e,\ \lambda\in J\rangle$ for the group $F$ ($f_j\neq e,\ \forall j $).
	\item  Group words $w_f(X_1,\ldots,X_r)$, $f\in F$, such that $f=w_f(f_1,\ldots,f_r)$. In the case $f=f_j$ ($j=1,\ldots,r$), we just set $w_{f_j}=X_j$, and in the case $f=e$, we set $w_e=e$, the empty word.
	\item A transversal $\{\beta_f: f\in F\}$ to the cosets of $N_p$ in $H_p$ such that $\beta_f N_p=f$ for all $f\in F$ and $\beta_e=1$ (the identity of $N_p$). We denote $\beta_i:=\beta_{f_i}$ (this is why we denoted the identity of $F$ by $e$ and not by 1).
		\item Elements $\alpha_1,\ldots,\alpha_s\in G_p$ whose classes modulo $N_p$ generate $G_p/N_p$. Assume that the first $t$ of them ($t\leq s$), $\alpha_1,\ldots,\alpha_t$, generate $G_p$ modulo $H_p$.
\end{enumerate}

\noindent Given $A\leq H_p$ of finite index such that $AN_p=H_p$, the intersection $A\cap N_p\subseteq N_p$ is open and $A$ can be written as 
$$A=(A\cap N_p)\bigcup\left(\bigcup_{f\in F\setminus\{e\}} \beta_f n_f(A\cap N_p)\right)$$
for some $n_f\in N_p$. Now, given an open subgroup $B\leq N_p$ and  $n_f\in N_p$ for $f\in F\setminus\{e\}$,  Lemma \ref{conditions for A=A(B,nf) to be a subgroup} below establishes necessary and sufficient conditions on the elements $n_f$ for the set
\begin{align*}
	A(B,(n_f)_{f\in F\setminus\{e\}}):=B\bigcup \left(\bigcup_{f\in F\setminus\{e\}} \beta_f n_f B\right)
\end{align*}
to be a subgroup of $G_p$. In this case, necessarily $A(B,(n_f)_{f\in F\setminus\{e\}})N_p=H_p$.
\begin{lemma}\label{conditions for A=A(B,nf) to be a subgroup}
	Fix an open subgroup $B\leq N_p$ and $n_f\in N_p$ for each  $f\in F\setminus\{e\}$. Let $n_j:=n_{f_j}$. 
	 Then the set $A:=A(B,(n_f)_{f\in F\setminus\{e\}})$ is a subgroup of $G_p$ if and only if the following three conditions are satisfied.
	\begin{enumerate}
		\item  $ {}^{\beta_jn_{j}}B\subseteq B$ for $1\leq j\leq r$. 
		\item $R_\lambda(\beta_{1}n_{1},\ldots,\beta_{r}n_{r})\in B$ for all $\lambda\in J$.
		\item $\beta_fn_f\in w_f(\beta_{1}n_{1},\ldots,\beta_{r}n_{r})B$ for all $f\in F\setminus\{e,f_1,\ldots,f_r\}.$
	\end{enumerate}
	Moreover,	$A$ is a normal subgroup of $G_p$ if and only if  (2), (3), and the following four conditions are satisfied.
	\begin{enumerate}
		\item[(4)] $B$ is normal in $N_p$.
		\item[(5)] $ {}^{\alpha_k}B\subseteq B$ for $1\leq k\leq s$.
		\item[(6)] $[\beta_{j}n_{j},x_i]\in B$ for $1\leq j\leq r$ and $1\leq i\leq h$.
		\item[(7)] $[\alpha_k,\beta_j n_j]\in w_{[\bar{\alpha}_k,f_j]}(\beta_1 n_1,\ldots,\beta_r n_r)B$ for $1\leq k\leq t$ and $1\leq j\leq r$, where $\bar{\alpha}_k$ denotes the class of $\alpha_k$ at $G_p/N_p$.
	\end{enumerate}
\end{lemma}
\begin{proof}
	Assume first that $A$ is a subgroup of $G_p$. Then clearly $B=A\cap N_p$, so $B$ is normal in $A$; in particular, (1) holds. Note that
	$
	R_\lambda(\beta_{1}n_{1},\ldots,\beta_{r}n_{r})N_p=R_\lambda(f_1,\ldots,f_r)=e\ (\in F),
	$ and
	hence $R_\lambda(\beta_{1}n_{1},\ldots,\beta_{r}n_{r})\in A\cap N_p=B$; thus, (2) holds. Similarly, $\beta_fn_f$ and $w_f(\beta_{1}n_{1},\ldots,\beta_{r}n_{r})$ have the same image at $F$, and hence (3) also holds.
	
	Conversely, assume that the conditions (1), (2), and (3) are satisfied for the collection $(n_f)_{f\in F\setminus\{e\}}$.
	Condition (1) implies that  $\beta_{j}n_{j}\in \on{N}_{H_p}(B)$, the  normalizer of $B$ in $H_p$, and then (3) implies that $\beta_f n_f\in \on{N}_{H_p}(B)$ for all $f\in F\setminus\{e\}$. It follows that  $A\subseteq \on{N}_{H_p}(B)$. Let $A'/B$ denote the subgroup of $\on{N}_{H_p}(B)/B$ generated by $A/B$. By (3), $A'/B$ is generated by the classes of $\beta_{1}n_{1},\ldots,\beta_{r}n_{r}$, and by (2), these generators satisfy the relations defining the presentation of $F$. It follows that,  $|A'/B|\leq |F|$. On the other hand, $|A/B|$ has exactly $|F|$ elements. Thus, $A/B=A'/B$, whence $A=A'$. We conclude that $A$ is a subgroup of $G_p$. This completes the proof of the first part of the lemma.
	
	Assume now that $A$ is a normal subgroup of $G_p$. According to the first part of the lemma, (2) and (3) are satisfied. Note that (4) and (5) also hold since $B=A\cap N_p$ is also normal in $G_p$. Now, the normality of $A$ implies that $ [\beta_j n_j,x_i]\in A$, and the normality of $N_p$ that $[\beta_j n_j,x_i]\in N_p$; thus $[\beta_j n_j,x_i]\in B=A\cap N_p$, whence (6) is also satisfied.  Finally, $(w_{[\bar{\alpha}_k,f_j]}(\beta_1 n_1,\ldots,\beta_r n_r))^{-1}[\alpha_k,\beta_j n_j]$ is in $A$ since $A$ is normal, and it is in $N_p$ since its image at $G_p/N_p$ is the identity. This shows that (7) also holds.
	
	Conversely, assume that (2) to (7) are satisfied. Note that (4) and (5) imply  (1); hence, according to the first part of the lemma, $A$ is a subgroup.  They also imply that $B$ is normal in $G_p$. We saw in the second paragraph of the proof that $A/B$ is generated by the classes of the elements $\beta_{j}n_{j}$; thus, by (6),  $N_p/B$ is included in $\on{N}_{G_p/B}(A/B)$. It follows that $A/B$ is normal in $H_p/B$.  Finally, by (7),  $A/B$ is normal in $G_p/B$, and hence $A$ is normal in $G_p$. This completes the proof of the lemma.
\end{proof}

We set more notation for the next proposition.
We denote by $\on{M}_{k\times h}(\mb{Z}_p)$ the set of $k\times h$-matrices with entries in $\mb{Z}_p$. Given $\v=(v_{ij})\in \on{M}_{k\times h}(\mb{Z}_p)$, its $i$-th row vector $(v_{i1},\ldots,v_{ih})\in\mb{Z}_p^h$ is denoted by $\v_i$. The normalized Haar measure of  $\on{T}_h(\mb{Z}_p)\times\on{M}_{k\times h}(\mb{Z}_p)\cong \mb{Z}_p^{kh+h(h+1)/2}$ is denoted by $\mu$. Elements of $\on{T}_h(\mb{Z}_p)\times\on{M}_{k\times h}(\mb{Z}_p)$ are denoted $(\t,\v)$. 

\begin{proposition}\label{integral expression for zeta_H^H and zeta_G^H} It holds:
	\begin{align*}\zeta_{H_p}^{H_p,\s}(s)=\frac{1}{(1-p^{-1})^h}\int_{\mc{T}_{H_p}^{H_p,\s}} \prod_{i=1}^h|t_{ii}|_p^{s-i-r}d\mu&&\mbox{and}&& \zeta_{G_p}^{H_p,\n}(s)=\frac{1}{(1-p^{-1})^h}\int_{\mc{T}_{G_p}^{H_p,\n}} \prod_{i=1}^h|t_{ii}|_p^{s-i-r}d\mu,
	\end{align*}
	where $\mc{T}_{H_p}^{H_p,\s}\subset \on{T}_h(\mb{Z}_p)\times \on{M}_{r,h}(\mb{Z}_p)$ is the set of pairs $(\t,\v)$ such that $\t\in\mc{M}_{N_p}^\s$ and 
	\begin{align}\label{domain of integration for subgroup}	 \ ^{\beta_{j}}[\x^{\v_j},\x^{\t_i}] [\beta_{j},\x^{\t_i}],\ R_\lambda(\beta_{1}\x^{\v_1},\ldots,\beta_{r}\x^{\v_r})\in \overline{\langle \x^{\t_{1}},\ldots,\x^{\t_h}\rangle},\ 1\leq i\leq h,\ 1\leq j\leq r,\ \lambda\in J,
	\end{align}
	and $\mc{T}_{G_p}^{H_p,\n}\subset \on{T}_h(\mb{Z}_p)\times \on{M}_{r,h}(\mb{Z}_p)$ is the set of pairs $(\t,\v)$ such that $\t\in\mc{M}_{N_p}^\n$ and 
	\begin{align}
	\label{domain of integration for normal subgroup} ^{\alpha_{k}}(\x^{\t_i}),\ R_\lambda(\beta_{1}\x^{\v_1},\ldots,\beta_{r}\x^{\v_r}),\  {}^{\beta_{j}}[\x^{\v_j},x_i][\beta_{j},x_i],\   (w_{[\bar{\alpha}_l,f_j]}(\beta_1\x^{\v_1},\ldots,\beta_r\x^{\v_r}))^{-1} [\alpha_l,\beta_j \x^{\v_j}] \in \overline{\langle \x^{\t_{1}},\ldots,\x^{\t_h}\rangle}\\
\nonumber	 1\leq k\leq s,\ 1\leq i\leq h,\  \lambda\in J,\ 1\leq j\leq r,\ 1\leq l\leq t. 
	\end{align}
\end{proposition}
\begin{proof}
	Throughout the proof, we shall refer repeatedly to conditions (1)-(7) of Lemma \ref{conditions for A=A(B,nf) to be a subgroup}. We will also use the same notation. The first part of this lemma establishes a bijection between the family $\ms{F}_{H_p}^{H_p,\s}$ of open subgroups $A\leq G_p$ such that $AN_p=H_p$ and the family of sequences $(n_fB)_{f\in F}$, where $B$ is an open subgroup of $N_p$ and the $n_f$'s are elements of $N$ with $n_e=1$ such that  (1), (2), and (3) are fulfilled. 
	As (3) simply expresses $n_fB$ in terms of $n_{1}B,\ldots,n_{r}B$ (recall that $n_j$ denotes $n_{f_j}$), it follows that $\ms{F}_{H_p}^{H_p,\s}$ is in a bijection with the set, say $\ms{T}$, of sequences $(B,n_{1} B,\ldots,n_{r} B)$ for which conditions (1) and (2) are satisfied. Thus,
	\begin{equation}\label{eq1:main prop} \zeta_{H_p}^{H_p,\s}(s):=\sum_{A\in \ms{F}_{H_p}^{H_p,\s}}[H_p:A]^{-s}=\sum_{(B,n_{1}B,\ldots,n_{r}B)\in\ms{T}}[N_p:B]^{-s}.
	\end{equation} 
	We now fix $(B,n_{1}B,\ldots,n_{r}B)\in\mathscr{T}$ and set  
	\begin{align*}
		\mc{S}(B,n_{1}B,\ldots, n_{r}B):=\{\v\in \on{M}_{r\times h}(\mb{Z}_p):(B,\x^{\v_1}B,\ldots,\x^{\v_r} B)= (B,n_{1} B,\ldots, n_{r} B)\}.
	\end{align*}
	This is an open subset of $\on{M}_{r\times h}(\mb{Z}_p)$ of Haar measure $[N_p:B]^{-r}$. In fact, the mapping $\mathbb{Z}_p^h\rightarrow N_p$ given by $\a\mapsto \x^\a$ is a homeomorphism that preserves the (normalized) Haar measure (cf. \cite[Lemma 2.4]{Su}), and  the measure of each coset $n_{j}B$ is $[N_p:B]^{-1}$. Combining this with Lemma \ref{properties of M(B)} and  formula (\ref{inverse of the index}), we deduce that $\mc{M}(B)\times \mc{S}(B,n_{1}B,\ldots, n_{r}B)$ is an open subset of $\on{T}_h(\mb{Z}_p)\times\on{M}_{r\times h}(\mb{Z}_p)$, and that for any $\t=(t_{ij})\in\mc{M}(B)$
	$$
	\mu(\mc{M}(B)\times \mc{S}(B,n_{f_1}B,\ldots, n_{f_r}B))= (1-p^{-1})^h\prod_{i=1}^h|t_{ii}|_p^i \prod_{i=1}^h|t_{ii}|_p^r =(1-p^{-1})^h\prod_{i=1}^h|t_{ii}|_p^{r+i}.
	$$
	It follows (again, using formula (\ref{inverse of the index})) that
	\begin{align*}
		[N_p:B]^{-s}&=\prod_{i=1}^h|t_{ii}|_p^s=(1-p^{-1})^{-h}\prod_{i=1}^h|t_{ii}|_p^{s-r-i} \mu(\mc{M}(B)\times \mc{S}(B,n_{f_1}B,\ldots, n_{f_r}B))\\
		&=\frac{1}{(1-p^{-1})^h}\int_{\mc{M}(B)\times \mc{S}(B,n_{f_1}B,\ldots, n_{f_r}B)}\prod_{i=1}^h|t_{ii}|_p^{s-r-i}d\mu.
	\end{align*}
 This and (\ref{eq1:main prop}) imply that
	\begin{align*}
		\zeta_{H_p}^{H_p,\s}(s)=\frac{1}{(1-p^{-1})^h}\int_{\mc{T}_{H_p}^{H_p,\s}}\prod_{i=1}^h|t_{ii}|_p^{s-r-i}d\mu,
	\end{align*}
	where $\mc{T}_{H_p}^{H_p,\s}:= \bigcup_{(B,n_{1}B,\ldots,n_{r}B)\in\ms{T}}\mc{M}(B)\times \mc{S}(B,n_{1}B,\ldots,n_{r}B)$.
	 Note that this is the set of pairs $(\t,\v)\in\on{T}_h(\mb{Z}_p)\times\on{M}_{r\times h}(\mb{Z}_p)$ such that (a) $\t\in\mathcal{M}_{N_p}^\s$; 
and (b) the elements $n_{j}:=\x^{\v_j}$, for $j=1,\ldots,r$, satisfy conditions (1) and (2)  with $B=\overline{\langle \x^{\t_{1}},\ldots,\x^{\t_h}\rangle}$. These are precisely the conditions listed in  (\ref{domain of integration for subgroup}) since $^{\beta_{j}}[\x^{\v_j},\x^{\t_i}] [\beta_{j},\x^{\t_i}]=[\beta_{j}\x^{\v_j}, \x^{\t_i}]$.  This proves the integral  expression for $\zeta_{H_p}^{H_p,\s}(s)$.

	The expression for $\zeta_{G_p}^{H_p,\n}(s)$ is obtained similarly. This time we have to consider sequences $(n_fB)_{f\in F}$ such that conditions (2) to (7) are satisfied. Working as in the previous case, we arrive at the expression
	\begin{align*}
		\zeta_{G_p}^{H_p,\n }(s)=\frac{1}{(1-p^{-1})^h}\int_{\mc{T}_{G_p}^{H_p,\n}}\prod_{i=1}^h|t_{ii}|_p^{s-r-i}d\mu	\end{align*}
	where $\mc{T}_{G_p}^{H_p,\n}=\cup \mc{M}(B)\times\mc{S}(B,n_{1}B,\ldots,n_{r}B)$, and the union runs over those $(B,n_{1}B,\ldots,n_{r}B)$ with $B$ open and normal in $N_p$ and the $n_j$'s satisfying (2), (5), (6), (7). 
	Therefore, $\mc{T}_{G_p}^{H_p,\n}$ consists of the pairs $(\t,\v)$ such that (a) $\t\in\mc{M}_{N_p}^\n$; (b) $^{\alpha_k}(\x^{\t_i})\in \overline{\langle \x^{\t_{1}},\ldots,\x^{\t_h}\rangle}$ for $i=1,\ldots,h$ and $k=1,\ldots,s$; and (c) the elements $n_{j}:=\x^{\v_j}$  ($j=1,\ldots,r$) satisfy conditions (2),  (6) and (7). Clearly
(b) and (c) are equivalent to (\ref{domain of integration for normal subgroup}) since
	$^{\beta_{j}}[\x^{\v_j},x_i][\beta_{j},x_i]=[\beta_{j}\x^{\v_j},x_i]$. This completes the proof of the proposition.
\end{proof}

\subsection{A method to simplify cone integrals}
We end this section with an elementary observation that we will use frequently to simplify the calculation of cone integrals.
 	
Given $m, n\in\mb{N}$, we consider subsets of $\mb{Z}_p^m\times \mb{Z}_p^n$ obtained as follows.
\begin{enumerate}[(i)]
	\item Fix a measurable subset $D_0\subset\mb{Z}_p^m$ and rational functions $g_i(\T)\in\mb{Q}_p(\T)=\mb{Q}_p(T_1,\ldots,T_m)$, $i=1,\ldots,n$, such that $g_i(\t)\in\mb{Z}_p$ whenever $\t\in D_0$.
	\item Given an integer $i$ such that $1\leq i\leq n$, assume that we have defined a measurable subset $D_{i-1}\subset\mb{Z}_p^m\times\mb{Z}_p^{i-1}$. Let $k_i(\T,V_1,\ldots,V_{i-1}),\ \lambda_i(\T,V_1,\ldots,V_{i-1})\in\mb{Q}_p(\T,V_1,\ldots,V_{i-1})$ be rational functions such that $k_i(\t,v_1,\ldots,v_{i-1})\in\mb{Z}_p$ and $\lambda_i(\t,v_1,\ldots,v_{i-1})\in\mb{Z}_p^*$ whenever $(\t,v_1,\ldots,v_{i-1})\in D_{i-1}$. We define $D_i$ in one of the following two ways:
	\begin{align*}
		\mbox{I:}\ \ & D_i=\left\{(\t,v_1,\ldots,v_i)\in\mb{Z}_p^m\times\mb{Z}_p^i: \begin{matrix}&(\t,v_1,\ldots,v_{i-1})\in D_{i-1},\\   &g_i(\t)\ |\ k_i(\t,v_1,\ldots,v_{i-1})+\lambda_i(\t,v_1,\ldots,v_{i-1}) v_i\end{matrix}\right\}\\
		\mbox{II:}\ \ & D_i=\left\{(\t,v_1,\ldots,v_i)\in\mb{Z}_p^m\times\mb{Z}_p^i: \begin{matrix}(\t,v_1,\ldots,v_{i-1})\in D_{i-1},\\ |g_i(\t)|_p=|k_i(\t,v_1,\ldots,v_{i-1})+\lambda_i(\t,v_1,\ldots,v_{i-1}) v_i|_p\end{matrix}\right\}.
	\end{align*}
	\item We finally set $D=D_n\subset\mb{Z}_p^m\times\mb{Z}_p^n$.  
\end{enumerate}
We say that $(D_0,g_1(\T),\ldots,g_n(\T))$ is the initial datum of definition of $D$ and call $(v_1,\ldots,v_n)$ the sequence of \emph{pivots}. A pivot $v_i$ is said to be of type I or type II according to the way we choose to define $D_i$.

\begin{proposition}\label{proceso de simplifiacion de integrales}
	Let $f_0,g_0\in\mb{Q}_p[\T]=\mb{Q}_p[T_1,\ldots,T_m]$ be non-zero polynomials, and consider the integral
	\begin{align}\label{int11}
		\int_{D\subseteq\mb{Z}_p^m\times\mb{Z}_p^n}|f_0(\t)|_p^s|g_0(\t)|_pd\mu(\t),\ \ s\in\mb{C},
	\end{align}
	where $D\subseteq \mb{Z}_p^m\times\mb{Z}_p^n$ is a measurable subset defined from a datum $(D_0,g_1(\T),\ldots,g_n(\T))$ and pivots $(v_1,\ldots,v_n)$. If the number of pivots of type II is $r$, then 
	\begin{align*}
		\int_{D}|f_0(\t)|_p^s|g_0(\t)|_pd\mu(\t)=
		(1-p^{-1})^r\int_{D_0}|f_0(\t)|_p^s\prod_{i=0}^n|g_i(\t)|_pd\mu(\t). 
	\end{align*}
\end{proposition} 
\begin{proof}
	The integral (\ref{int11}) can be performed as follows. We start integrating with respect to the variable $v_n$. Note that $g_n(\t)|k_n(\t,v_1,\ldots,v_{n-1})+\lambda_n(v_1,\ldots,v_{n-1})v_n$ if and only if $v_n\in g_n(\t)\mb{Z}_p-\frac{k_n(\t,v_1,\ldots,v_{n-1})}{\lambda_n(\t,v_1,\ldots,v_{n-1})}$, and the Haar measure of this set is  $|g_n(\t)|_p$. Similarly, $|g_i(\t)|_p=|k_n(\t,v_1,\ldots,v_{n-1})+\lambda_n(\t,v_1,\ldots,v_{n-1})v_n|_p$ if and only if $v_n\in g_n(\t)\mb{Z}_p^*-\frac{k_n(\t,v_1,\ldots,v_{n-1})}{\lambda_n(\t,v_1,\ldots,v_{n-1})}$ and the Haar measure of the this set is $(1-p^{-1})|g_n(\t)|_p$. To sum up, after integrating with respect to $v_n$, the integrand is multiplied by $|g_n(\t)|_p$ or by $(1-p^{-1})|g_n(\t)|_p$ according to whether $v_n$ is a pivot of type I or II, and the new domain of integration is clearly $D_{n-1}$. We next integrate with respect to $v_{n-1}$ and so on. After $n$ steps, we arrive at the desired form of the integral.
\end{proof} 

\section{The formulae}\label{Section: The formulae}

\noindent We now present the complete formulae for the subgroup and normal zeta functions of the 3-dimensional almost-Bieberbach groups, whose definition we recall below. Proofs of these formulae will be left for the next two and last sections. For details of the next discussion, we refer to \cite{DeK}.

Let $\mc{N}$ be a connected and simply connected nilpotent Lie group, and let $\on{Aut}(\mc{N})$ be the topological group of Lie group automorphisms of $\mc{N}$.
The semi-direct product $\mc{N}\rtimes \on{Aut}(\mc{N})$ acts on $\mc{N}$ in a canonical way by 
$$^{(n,\alpha)}x=n\alpha(x),\ \ \forall x,n\in\mc{N},\ \ \forall \alpha\in\on{Aut}(\mc{N}).$$
All maximal compact subgroups of $\on{Aut}(\mc{N})$ are conjugate. We fix one of them, say $\mc{C}\subset \on{Aut}(\mc{N})$.  A torsion-free uniform discrete subgroup $G$ of $\mc{N}\rtimes \mc{C}$ is called an almost-Bieberbach group (abbreviated as AB-group). Note that the quotient space $_{G}\backslash{\mc{N}}$ is a compact manifold of the same dimension as $\mc{N}$ whose fundamental group is identified with $G$. The dimension of $G$ is defined as the dimension of $\mc{N}$.	
When $\mc{N}=\mb{R}^n$ and $\mc{C}=O(\mb{R}^n)$, the orthogonal group, we recover the definition of Bieberbach group. The quotient space $_G\backslash\mb{R}^n$, with the metric induced from the Euclidean space $\mb{R}^n$, is a compact flat manifold.

The AB-groups were characterized algebraically as the finitely generated torsion-free virtually nilpotent groups  (cf.\ \cite[Theorem 3.1.3]{DeK}). If $G$ is an AB-group, then its Fitting subgroup $\on{Fitt}(G)$ (the maximal normal nilpotent subgroup of $G$) is indeed maximal nilpotent. It follows that any intermediate subgroup $\on{Fitt}(G)\leq H\leq G$ is also an AB-group with $\on{Fitt}(H)=N$.
An AB-group $G$ is a Bieberbach group if and only if $\on{Fitt}(G)$ is abelian.

\subsection{The 3-dimensional Bieberbach groups and their zeta functions}\label{3-dimensional Bieberbach groups and their zeta functions}
There are only ten 3-dimensional Bieberbach groups up to isomorphism (cf. \cite[Chap. 3]{Wolf1972Spaces}). They are listed below with their corresponding zeta functions expressed in terms of the partial zeta functions with respect to the Fitting subgroup (which is always $\langle  x_1,x_2,x_3\rangle$). 
 The first six are the fundamental groups of the orientable compact flat manifolds, and the last four are the fundamental groups of the non-orientable ones. 
We have arranged the presentation so that $\zeta_G^{G\s}(s)$ appears as the first term in the formula of $\zeta_G^\s(s)$ and $[G:N]^{-s}\zeta_N^\s(s)$ as the last one. The same holds for $\zeta_G^\n(s)$. 

\medskip 
$\ms{G}_1=\langle x_1,x_2,x_3: [x_1,x_2]=[x_1,x_3]=[x_2,x_3]=1\rangle$ has zeta functions
\begin{small} 
\begin{align*}
	\zeta_{\ms{G}_1}^\s(s)&=\zeta_{\ms{G}_1}^\lhd(s)=\zeta(s)\zeta(s-1)\zeta(s-2).
\end{align*}
\end{small} 
\hrule 
\smallskip 
$\ms{G}_2=\langle \alpha,x_1,x_3,x_3:[x_i,x_j]=1\ \forall i,j,\ \alpha^2=x_3,\ ^\alpha x_1=x_1^{-1},\ ^\alpha x_2=x_2^{-1} \rangle$ has zeta functions
\begin{small} 
\begin{align*} 
\zeta_{\ms{G}_2}^\s(s)&=\zeta(s)\zeta(s-1)\zeta(s-2)(1-2^{-s})+	2^{-s}\zeta(s)\zeta(s-1)\zeta(s-2)=\zeta(s)\zeta(s-1)\zeta(s-2),\\
	\zeta_{\ms{G}_2}^\lhd(s)&=(1+6\cdot 2^{-s}+4\cdot 4^{-s})(1-2^{-s})\zeta(s)+2^{-s}\zeta(s)^2\zeta(s-1)(1+3\cdot 2^{-s}).
\end{align*}
\end{small} 
\hrule 
\smallskip 
 $\ms{G}_3=\langle \alpha,x_1,x_2,x_3:[x_i,x_j]=1\ \forall i,j,\ \alpha^{3}=x_1,\  ^\alpha x_2=x_3,\  ^\alpha x_3=x_2^{-1}x_3^{-1}\rangle$ has zeta functions 
\begin{small}  
\begin{align*}
	\zeta_{\ms{G}_3}^\s(s)&=(1-3^{-s})\zeta(s)\zeta(s-1)L(s-1,\chi_3)+3^{-s}\zeta(s)\zeta(s-1)\zeta(s-2),\\
	\zeta_{\ms{G}_3}^\lhd(s)&=(1+3\cdot 3^{-s})(1-3^{-s})\zeta(s)+3^{-s}\zeta(s)^2L(s,\chi_3)(1+2\cdot 3^{-s}).
\end{align*}
\end{small} 
\hrule 
\smallskip 
 $\ms{G}_4=\langle \alpha,x_1,x_2,x_3:[x_i,x_j]=1\ \forall i,j,\ \alpha^4=x_1,\ ^\alpha x_2=x_3,\ ^\alpha x_3=x_2^{-1} \rangle$ has zeta functions
 \begin{small} 
\begin{align*} 
	\zeta_{\ms{G}_4}^\s(s)&=\zeta(s)\zeta(s-1)L(s-1,\chi_4)(1-2^{-s})+2^{-s}\zeta(s)\zeta(s-1)\zeta(s-2)(1-2^{-s})+4^{-s}\zeta(s)\zeta(s-1)\zeta(s-2),
	\\ 
	\zeta_{\ms{G}_4}^\lhd(s)&=(1+2\cdot 2^{-s})(1-2^{-s})\zeta(s)+2^{-s}(1+2\cdot 2^{-s}+2\cdot 4^{-s})(1-2^{-s})\zeta(s)+4^{-s}\zeta(s)^2L(s,\chi_4)(1+2^{-s}).
\end{align*}
\end{small} 
\hrule 
\smallskip 
$\ms{G}_5=\langle \alpha,x_1,x_2,x_3:[x_i,x_j]=1\ \forall i,j,\ \alpha^6=x_1,\ ^\alpha x_2=x_3,\ ^\alpha x_3=x_2^{-1}x_3 \rangle$ has zeta functions
\begin{small} 
\begin{align*}
	\zeta_{\ms{G}_5}^\s(s)=&\zeta(s)\zeta(s-1)L(s-1,\chi_3)(1-2^{-s})(1-3^{-s})+2^{-s}\zeta(s)\zeta(s-1)L(s-1,\chi_3)(1-3^{-s})\\
	&+3^{-s}\zeta(s)\zeta(s-1)\zeta(s-2)(1-2^{-s}) +6^{-s}\zeta(s)\zeta(s-1)\zeta(s-2),\\
	\zeta_{\ms{G}_5}^\lhd(s)=&(1-2^{-s})(1-3^{-s})\zeta(s)+2^{-s}(1+3^{-s})(1-3^{-s})\zeta(s)+3^{-s}(1+4^{-s})(1-2^{-s})\zeta(s)+6^{-s} \zeta(s)^2L(s,\chi_3).
\end{align*} 
\end{small} 
\hrule 
\smallskip 
 $\ms{G}_6=\left\langle
\alpha,\beta,x_1,x_2,x_3: \begin{array}{l}
	[x_i,x_j]=1\ \forall i,j,\ 
	\alpha^2=x_1,\ ^\alpha x_2=x_2^{-1},\ ^\alpha x_3=x_3^{-1},\\  \beta^2=x_2,\ ^\beta x_1=x_1^{-1},\ ^\beta x_3=x_3^{-1},\ (\alpha\beta)^2=x_3 \end{array} \right\rangle
$ has zeta functions
\begin{small} 
\begin{align*}
	\zeta_{\ms{G}_6}^\s(s)&=\zeta(s-1)^3(1-2\cdot 2^{-s})^3+3\cdot 2^{-s}\zeta(s)\zeta(s-1)\zeta(s-2)(1-2^{-s})+4^{-s}\zeta(s)\zeta(s-1)\zeta(s-2),\\
	\zeta_{\ms{G}_6}^\n(s)&=1+3\cdot 2^{-s}(1+2\cdot 2^{-s})(1-2^{-s})\zeta(s)+4^{-s}\zeta(s)^3(1+4\cdot 2^{-s}+4^{-s}).
\end{align*}
\end{small} 
\hrule 
\smallskip 
$\ms{B}_1=\langle \varepsilon, x_1,x_2,x_3:[x_i,x_j]=1\ \forall i,j,\  \varepsilon^2=x_1,\ ^\varepsilon x_2=x_2,\ ^\varepsilon x_3=x_3^{-1}\rangle$ has zeta functions
\begin{small} 
 \begin{align*}
	\zeta_{\ms{B}_1}^\s(s)&=\zeta(s)\zeta(s-1)^2(1-2^{-s})(1+2^{-(s-1)})+2^{-s}\zeta(s)\zeta(s-1)\zeta(s-2),\\
	\zeta_{\ms{B}_1}^\lhd(s)&=(1-2^{-s})(1+2^{-(s-2)})\zeta(s)\zeta(s-1)+2^{-s}(1+3\cdot 2^{-s})\zeta(s)^2\zeta(s-1).
\end{align*}
\end{small} 
\hrule 
\smallskip 
$\ms{B}_2=\langle \varepsilon,x_1,x_2,x_3: [x_i,x_j]=1\ \forall i,j,\ \varepsilon^2=x_1,\ ^\varepsilon x_2=x_2,\ ^\varepsilon x_3=x_1x_2x_3^{-1}\rangle$ has zeta functions
\begin{small} 
\begin{align*}
	\zeta_{\ms{B}_2}^\s(s)&=\zeta(s)\zeta(s-1)^2(1-2^{-s})(1-2\cdot 2^{-s}+8\cdot 4^{-s})+2^{-s}\zeta(s)\zeta(s-1)\zeta(s-2),\\
	\zeta_{\ms{B}_2}^\lhd(s)&=\zeta(s)\zeta(s-1)(1-2^{-s})+2^{-s}\zeta(s)^2\zeta(s-1)(1-2^{-s}+4\cdot 4^{-s}).
\end{align*}
\end{small} 
\hrule 
\smallskip 
$\ms{B}_3=\left\langle
\alpha,\varepsilon, x_1,x_2,x_3:\begin{array}{l}[x_i,x_j]=1\ \forall i,j,\
	\alpha^2=x_1,\  ^\alpha x_2=x_2^{-1},\ ^\alpha x_3=x_3^{-1},\\ 
	\varepsilon^2=x_2,\ ^\varepsilon x_1=x_1,\ ^\varepsilon x_3=x_3^{-1},\ [\varepsilon,\alpha]=x_2\end{array}
\right\rangle$ has zeta functions
\begin{small} 
\begin{align*}
	\zeta_{\ms{B}_3}^\s(s)=&\zeta(s)\zeta(s-1)^2(1-2^{-s})(1-4\cdot 4^{-s})+2\cdot 2^{-s} \zeta(s)\zeta(s-1)^2(1-2^{-s})(1+2\cdot 2^{-s})\\
	&+2^{-s}\zeta(s)\zeta(s-2)\zeta(s-3)(1-2^{-s})+4^{-s}\zeta(s)\zeta(s-1)\zeta(s-2),\\
	\zeta_{\ms{B}_3}^\lhd(s)=&\zeta(s)(1-2^{-s})(1+3\cdot 2^{-s})+2^{-s}\zeta(s)^2(1-2^{-s})(1+2^{-s}+2\cdot 4^{-s})+2^{-s}\zeta(s)(1-2^{-s})(1+2\cdot 2^{-s})\\
	&+2^{-s}\zeta(s)^2(1-2^{-s})(1+5\cdot 2^{-s}+2\cdot 4^{-s})+4^{-s}\zeta(s)^3(1+4\cdot 2^{-s}+ 4^{-s}).
	\end{align*}
\end{small} 
\hrule 
\smallskip 
$\ms{B}_{4}=\left\langle
\alpha,\varepsilon,x_1,x_2,x_3:\begin{array}{l} [x_i,x_j]=1\ \forall i,j,\
	\alpha^2=x_1,\ ^\alpha x_2=x_2^{-1},\ 	^\alpha x_3=x_3^{-1},\\
	 \varepsilon^2=x_2,\ ^\varepsilon x_1=x_1,\ ^\varepsilon x_3=x_3^{-1},\ [\varepsilon,\alpha]=x_2x_3 \end{array}\right\rangle$ has zeta functions
 \begin{small} 
\begin{align*}
	\zeta_{\ms{B}_4}^\s(s)&=\zeta(s)\zeta(s-1)^2(1-2^{-s})(1-2^{-(s-1)})^2+2^{-s}\zeta(s)\zeta(s-1)\zeta(s-2)(1-2^{-s})\\
&+2\cdot 2^{-s} \zeta(s)\zeta(s-1)^2(1-2^{-s})(1+2^{-(s-1)})+	4^{-s}\zeta(s)\zeta(s-1)\zeta(s-2),\\
	\zeta_{\ms{B}_4}^\lhd(s)&=4^{-s}\zeta(s)^3(1+4\cdot 2^{-s}+ 4^{-s}) +2^{-s}\zeta(s)(1-2^{-s})(1+2\cdot 2^{-s})\\
	&+2^{-s}\zeta(s)^2(1-2^{-s})(1+ 2^{-s}+2\cdot 4^{-s})+2^{-s}\zeta(s)^2(1-2^{-s})^2(1+2\cdot 2^{-s})+\zeta(s)(1-2^{-s}).
\end{align*}
\end{small}

\subsection{The 3-dimensional almost-Bieberbach groups and their zeta functions}\label{3-dimensional AB groups and their zeta functions}

The 3-dimensional AB-groups with non-abelian Fitting subgroup were classified in \cite{DIKL}; see also \cite[Chap. 7]{DeK}. They are arranged into seven families according to the isomorphism type of  $G/Z(N)$ ($G$ is the AB-group and $N$ the Fitting subgroup), which is one of the plane crystallographic groups $\mathbf{p1}$, $\mathbf{p2}$, $\mathbf{pg}$, $\mathbf{p3}$, $\mathbf{p4}$, $\mathbf{p6}$, $\mathbf{p2gg}$. 
They are listed below with their zeta functions expressed in terms of the partial zeta functions with respect to the Fitting subgroup, which is always $\langle x_1,x_2,x_3\rangle$. This time the normal zeta function has been arranged so that $[G:N]^{-s}\zeta_G^{N,\n}(s)$ appears as the first term in the formula and $\zeta_{G}^{G,\n}(s)$ as the last one.

\medskip 
	\noindent $N_k=\left\langle x_1,x_2,x_3: [x_2,x_1]=x_3^k,\  [x_1,x_3]=[x_2,x_3]=1\right\rangle$, $k\in\mb{N}$, has zeta functions
	\begin{small} 
		\begin{flalign*}
		\zeta_{N_k}^\s(s)=&\prod_{p\nmid k}\frac{\zeta_p(s)\zeta_p(s-1)\zeta_p(2s-2)\zeta_p(2s-3)}{\zeta_p(3s-3)}\times\prod_{p|k}\left(\zeta_p(s-2)\left(\zeta_p(s)\zeta_p(s-1)-p^{-s+2}|k|_p^{s-2}\zeta_p(2s-2)\zeta_p(2s-3) \right)\right),&
	\end{flalign*}
\begin{flalign*} 
		\zeta_{N_k}^\n(s)=&\prod_{p\nmid k}\zeta_p(s)\zeta_p(s-1)\zeta_p(3s-2) \times \prod_{p\mid k}\zeta_p(s)\zeta_p(s-1) \left(
		\frac{1-|k|_p^{s-2}}{1-p^{-s+2}}+ |k|_p^{s-2}\zeta_p(3s-\eta)\right).&
	\end{flalign*}
\end{small} 
\hrule 
\smallskip 
\noindent  $G_{\mathbf{p2},2k}=\left\langle \alpha, x_1,x_2,x_3:\begin{array}{l} [x_2,x_1]=x_3^{2k},\ [x_1,x_3]=[x_2,x_3]=1,\\
		 \alpha^2=x_3,\ ^{\alpha}x_1=x_1^{-1},\ ^{\alpha}x_2=x_2^{-1}\end{array} \right\rangle$, $k\in\mb{N}$, has zeta functions
	 \begin{small} 
	\begin{flalign*}
		\zeta_{G_{\mathbf{p2},2k}}^\s(s)=&\prod_{p\nmid 2k} \frac{\zeta_p(s-1)\zeta_p(s-2)\zeta_p(2s-1)\zeta_p(2s-2)}{\zeta_p(3s-3)}&\\
		&\times\prod_{\substack{p|k, p\neq 2}} \zeta_p(s)\left(\zeta_p(s-1)\zeta_p(s-2)-p^{-s}|k|_p^{s}\zeta_p(2s-1)\zeta_p(2s-2)\right)\times \zeta_2(s)\zeta_2(s-1)+2^{-s}\zeta_{N_{2k}}^\s(s),&\end{flalign*}
	\begin{flalign*} 
	\zeta_{G_{\mathbf{p2},2k}}^\n(s)=&2^{-s}\prod_{p\nmid 2k}\zeta_p(s)\zeta_p(s-1)\zeta_p(3s) \times \prod_{\substack{p\mid k, p\neq 2}} \zeta_p(s)\zeta_p(s-1) \left(
	\frac{1-|k|_p^{s}}{1-p^{-s}}+ |k|_p^{s}\zeta_p(3s)\right)&\\
	&\times \zeta_2(s)\zeta_2(s-1)\left(1+\frac{ 2^{2-s}(1-|k|_2^s)}{1-2^{-s}}+\zeta_2(3s)2^{2-s}|k|_2^s\right)+ \prod_{\substack{p\mid k, p\neq 2}}\left(\frac{1-p^{-s}|k|_p^s}{1-p^{-s}}\right)\times \left( 1+6\cdot 2^{-s}+4\cdot 4^{-s}\right). &
\end{flalign*}
\end{small} 
\hrule 
\smallskip 
\noindent  $G_{\mathbf{pg},2k}=\left\langle\beta,x_1,x_2,x_3:\begin{array}{l}[x_2,x_1]=x_3^{2k},\ [x_1,x_3]=[x_2,x_3]=1,\\
		 \beta^2=x_2,\  ^\beta x_1=x_1^{-1}x_3^{-k},\ ^\beta x_3=x_3^{-1}\end{array} \right\rangle$, $k\in\mb{N}$, has zeta functions
	 \begin{small} 
		\begin{flalign*}
			\zeta_{G_{\mathbf{pg},2k}}^\s(s)=&\prod_{p\nmid 2k}\frac{\zeta_p(s)\zeta_p(s-1)\zeta_p(2s-2)\zeta_p(2s-3)}{\zeta_p(3s-3)}\times \prod_{\substack{p\mid k, p\neq 2}} \zeta_p(s-2)\left(\zeta_p(s)\zeta_p(s-1)-p^{-s+2}|k|_p^{s-2}\zeta_p(2s-2)\zeta_p(2s-3) \right)&\\
			&\times \zeta_2(s-1)\left(\zeta_2(s-2)(1-|k|_2^{s-2})+\zeta_2(2s-3)|k|_2^{s-2}\right)+2^{-s}\zeta_{N_{2k}}^\s(s),&
			\end{flalign*}
		\begin{flalign*} 
	\zeta_{G_{\mathbf{pg},2k}}^\n(s)=&2^{-s}\prod_{p\nmid 2k} \zeta_p(s)^2\zeta_p(3s-1)\times \prod_{\substack{p\mid k, p\neq 2}}\zeta_p(s)^2 \left(
	\frac{1-|k|_p^{s-1}}{1-p^{-s+1}}+ |k|_p^{s-1}\zeta_p(3s-1)\right)\times\\
	&\times \zeta_2(s)^2\left(
	\zeta_2(s-1)(1+3\cdot 2^{-s})(1-|k|_2^{s-1})+\zeta_2(3s-1)(1+2^{-s}+6\cdot 4^{-s}-2\cdot 8^{-s}-2\cdot 16^{-s})|k|_2^{s-1} \right)&\\
	&+(1+(4+(-1)^k2)\cdot 2^{-s}+4\cdot 4^{-s})(1-2^{-s})\zeta(s).&
\end{flalign*}
\end{small} 
\hrule 
\smallskip 
\noindent 	$\displaystyle 
	G_{\mathbf{p3},k,\epsilon}=\left\langle \gamma,x_1,x_2,x_3:\begin{array}{l} [x_2,x_1]=x_3^k,\ [x_1,x_3]=[x_2,x_3]=1,\\
		 \gamma^3=x_3^{\epsilon},\ ^\gamma x_1=x_2,\ ^\gamma x_2=x_1^{-1}x_2^{-1}\end{array}\right\rangle$, with  $k\in\mb{N}$ and $\epsilon\in\{1,-1\}$ verifying $k(k+\epsilon)\equiv 0\mod 3$,
has zeta functions
\begin{small} 
\begin{flalign*}
	\zeta_{G_{\mathbf{p3},k,\epsilon}}^\s(s)=&\prod_{p\nmid 3k} \frac{\zeta_p(s-1)L_p(s-1,\chi_3)\zeta_p(2s-1)L_p(2s-1,\chi_3)}{L_p(3s-2,\chi_3)}\times&\\
	&\times  \prod_{\substack{p\mid k, p\neq 3}} \zeta_p(s)\left(L_p(s-1,\chi_3)\zeta_p(s-1)-p^{-s}|k|_p^{s}L_p(2s-1,\chi_3)\zeta_p(2s-1)\right)\times \zeta_3(s-1)+3^{-s}\zeta_{N_k}^\s(s),&
\end{flalign*}
\begin{flalign*}
	\zeta_{G_{\mathbf{p3},k,\epsilon}}^\n(s)=&3^{-s}\prod_{p\nmid 3k}\zeta_p(s)\zeta_p(3s)L_p(s,\chi_3)\times \prod_{p\mid k, p\neq 3} \zeta_p(s)L_p(s,\chi_3) \left(
	\frac{1-|k|_p^{s}}{1-p^{-s}}+ |k|_p^{s}\zeta_p(3s)\right)\times&\\
	&\times  \left(\zeta_{3}(s)+3\zeta_{3}(s)\left(3^{-s}\zeta_{3}(s)(1-| k|_{3}^{s})+3^{-3s}\zeta_{3}(3s)| k|_{3}^{s} \right)\right)+ \prod_{\substack{p\mid k, p\neq 3}} \frac{1-p^{-s}|k|_p^s}{1-p^{-s}}\times (1+3\cdot 3^{-s}).&
\end{flalign*}
\end{small} 
\hrule 
\smallskip 
\noindent 	 $G_{\mathbf{p4},2k,\epsilon}=\left\langle \gamma,x_1,x_2,x_3: \begin{array}{l}[x_2,x_1]=x_3^{2k},\ [x_1,x_3]=[x_2,x_3]=1,\\ \gamma^4=x_3^\epsilon,\ ^{\gamma}x_1=x_2,\ ^{\gamma}x_2=x_1^{-1}\end{array} \right\rangle
	$, with $k\in\mb{N}$ and $\epsilon\in\{1,-1\}$, has zeta functions
	\begin{small} 
		\begin{flalign*}
	\zeta_{G_{\mathbf{p4},2k,\epsilon}}^\s(s)=&\prod_{p\nmid 2k}\frac{\zeta_p(s-1)L_p(s-1,\chi_4)\zeta_p(2s-1)L_p(2s-1,\chi_{4})}{L_p(3s-2,\chi_{4})}\times&\\
	&\times \prod_{\substack{p\mid k, p\neq 2}}	\zeta_p(s)\left(L_p(s-1,\chi_d)\zeta_p(s-1)-p^{-s}|k|_p^{s}L_p(2s-1,\chi_d)\zeta_p(2s-1)\right)\times \zeta_2(s-1)+2^{-s}\zeta_{G_{\mathbf{p2},2k}}^\s(s),&
	\end{flalign*}
 \begin{flalign*}
 		\zeta_{G_{\mathbf{p4},2k,\epsilon}}^\n(s)=&4^{-s}\prod_{p\nmid 2k} \zeta_p(s)\zeta_p(3s)L_p(s,\chi_4)\times \prod_{\substack{p\mid k, p\neq 2}}	\zeta_p(s)L_p(s,\chi_d) \left(
 		\frac{1-|k|_p^{s}}{1-p^{-s}}+ |k|_p^{s}\zeta_p(3s)\right)\times &\\
 		&\times \left(\zeta_{2}(s)+2\zeta_{2}(s)\left(2^{-s}\zeta_{2}(s)(1-|2 k|_{2}^{s})+2^{-3s}\zeta_{2}(3s)|2 k|_{2}^{s} \right)\right)\\
 		& 	+2^{-s}\prod_{\substack{p\mid k, p\neq 2}}\frac{1-|k|_p^sp^{-s}}{1-p^{-s}}\times (1+2\cdot 2^{-s}+2\cdot 4^{-s}) +\prod_{\substack{p|k, p\neq 2}}\frac{1-|k|_p^sp^{-s}}{1-p^{-s}}\times (1+2\cdot 2^{-s}).&
 \end{flalign*}
\end{small} 
\hrule 
\smallskip 
\noindent $G_{\mathbf{p6},2k,\epsilon}=\displaystyle  \left\langle \gamma,x_1,x_2,x_3:\begin{array}{l}[x_2,x_1]=x_3^{2k},\ [x_1,x_3]=[x_2,x_3]=1,\\  \gamma^6=x_3^\epsilon,\ ^\gamma x_1=x_2,\ ^\gamma x_2=x_1^{-1}x_2\end{array}\right\rangle
$, $k\in\mb{N}$ and $\epsilon\in\{1,-1\}$ verifying $k(\epsilon+k)\equiv 0\mod 3$, has zeta functions
\begin{small} 
\begin{flalign*}
\zeta_{G_{\mathbf{p6},2k,\epsilon}}^\s(s)=&\prod_{p \nmid 6k} \frac{\zeta_p(s-1)\zeta_p(2s-1) L_p(s-1,\chi_3)L_p(2s-1,\chi_3)}{L_p(3s-2,\chi_3)}\cdot&\\
&\times\prod_{\substack{p\mid k, p\nmid 6}} \zeta_p(s)\left(L_p(s-1,\chi_6)\zeta_p(s-1)-p^{-s}|k|_p^{s}L_p(2s-1,\chi_6)\zeta_p(2s-1)\right)\times  \zeta_3(s-1)\zeta_2(2s-2)&\\
&+ 2^{-s}\zeta_{G_{\mathbf{p3},2k,\epsilon}}^{G_{\mathbf{p3},2k,\epsilon},\s}(s)+3^{-s}\zeta_{G_{\mathbf{p2},2k}}^{G_{\mathbf{p2},2k},\s}(s)+6^{-s}\zeta_{N_{2k}}^\s(s),&	
\end{flalign*}
\begin{flalign*}
		\zeta_{G_{\mathbf{p6},2k,\epsilon}}^\lhd(s)=&6^{-s}\prod_{p \nmid 6k} \zeta_p(s)\zeta_p(3s)L_p(s,\chi_3)\times \prod_{\substack{p\mid 2k, p\neq 3}} \zeta_p(s)L_p(s,\chi_3) \left(
		\frac{1-|2k|_p^{s}}{1-p^{-s}}+ |2k|_p^{s}\zeta_p(3s)\right)\times&\\
		&\times\left(\zeta_{3}(s)+\zeta_{3}(s)\left(3^{-s}\zeta_{3}(s)(1-| k|_{3}^{s})+3^{-3s}\zeta_{3}(3s)| k|_{3}^{s} \right)\right)&\\
		&+3^{-s}\prod_{\substack{p\mid k, p\neq 2}}\frac{1-p^{-s}|k|_p^s}{1-p^{-s}}\times (1+4^{-s})+2^{-s}\prod_{\substack{p\mid k, p\neq 3}}\frac{1-p^{-s}|2k|_p^s}{1-p^{-s}}\times (1+3^{-s})+\prod_{\substack{p\mid k, p\nmid 6}}\frac{1-p^{-s}|k|_p^s}{1-p^{-s}}.&
\end{flalign*}
\end{small} 
\hrule 
\smallskip 
	\noindent  $G_{\mathbf{p2gg},4k}=\left\langle 
	\alpha,\beta,x_1,x_2,x_3: \begin{array}{l} [x_2,x_1]=x_3^{4k},\ [x_1,x_3]=[x_2,x_3]=1,\\
		 \alpha^2=x_3,\ ^{\alpha}x_1=x_1^{-1}x_3^{2k}, ^{\alpha}x_2=x_2^{-1} x_3^{-2k},\\ 
		\beta^2=x_1,\ ^{\beta}x_2=x_2^{-1}x_3^{2k},\ ^{\beta}x_3=x_3^{-1},\ (\beta\alpha)^2=x_2
	\end{array}
	\right\rangle 
	$, $k\in\mb{N}$, has zeta functions
	\begin{small} 
	\begin{flalign*}
	\zeta_{G_{\mathbf{p2gg},4k}}^\s(s)=&	\prod_{p\nmid 2k}\frac{\zeta_p(s-1)^2\zeta_p(2s-2)^2}{\zeta_p(3s-3)}\times  \prod_{\substack{p\mid k, p\neq 2}} \zeta_p(s-1)\left((\zeta_p(s-1))^2-p^{-s+1}|k|_p^{s-1}(\zeta_p(2s-2))^2\right)&\\
	&+2\cdot 2^{-s}\zeta_{G_{\mathbf{pg},4k}}^{G_{\mathbf{pg},4k},\s}(s)+2^{-s}\zeta_{G_{\mathbf{p2},4k}}^{G_{\mathbf{p2},4k},\s}(s)+4^{-s}\zeta_{N_{4k}}^\s(s),&
	\end{flalign*}
	\begin{flalign*}
		\zeta_{G_{\mathbf{p2gg},4k}}^\n(s)=&4^{-s}\prod_{p\nmid 2k}\zeta_p(s)^2\zeta_p(3s)\times\prod_{\substack{p|k, p\neq 2}}
		 \zeta_p(s)^2 \left(
		 \frac{1-|k|_p^{s}}{1-p^{-s}}+ |k|_p^{s}\zeta_p(3s)\right)\times & \\
		 &\times \zeta_2(s)^2\left(\zeta_2(s)(1+4\cdot 2^{-s}+4^{-s})(1-|k|_2^{s})+\zeta_2(3s)(1+5\cdot 2^{-s}+2\cdot 4^{-s}+8^{-s}-16^{-s}-2\cdot 32^{-s})|k|_2^s \right)& \\
		&+2^{-s}\prod_{\substack{p\mid k, p\neq 2}}\frac{1-p^{-s}|k|_p^s}{1-p^{-s}}\times(1+2\cdot 2^{-s})+2\cdot 2^{-s}(1+2\cdot 2^{-s})(1-2^{-s})\zeta(s)+1.&  
	\end{flalign*}
\end{small} 

\begin{small} 
\section{Computing the zeta functions of the 3-dimensional Bieberbach groups}\label{Section: Proofs for Bieberbach groups}

\noindent We prove the formulae presented in \ref{3-dimensional Bieberbach groups and their zeta functions}. Throughout this section, we denote by $G$ the Bieberbach group under consideration and by $N$ the Fitting subgroup. We keep the notation introduced at the beginning of Section \ref{Section: Local factors at good primes} except that here $N$ is not a $\mf{T}_2$-group but rather an abelian group (and it is indeed an $\mb{Z}[P]$-module, where $P=G/N$). Instead, we denote
$$Z=C_N(P),\quad E=G/Z,\quad T=N/Z.
$$

The formula for the zeta functions of $\ms{G}_1$ are already known (see Introduction), so we focus here on the other groups.
We mentioned that an intermediate subgroup $N\leq H\leq G$ is again a Biebarbach group with Fitting subgroup $N$.
Therefore, to prove the formula for $\zeta_G^{G,\s}(s)$, it is enough to prove the one for $\zeta_G^{G,\s}(s)$.  The isomorphism classes of the intermediate subgroups $N\subsetneqq H\subsetneqq G$ (if there is any) will be identified when computing the partial zeta functions $\zeta_G^{H,\n}(s)$. 
The following will be useful. If $p\nmid |G/N|$, then
\begin{align}\label{general remark}
	\zeta_{G_p}^{G_p,\s}(s)&=\zeta_{E_p}^{E_p,\s}(s)\zeta_{Z_p}^\s(s)=\zeta_{E_p}^{T_p,\n}(s-1)\zeta_{Z_p}^\s(s),\\
\nonumber	 \zeta_{G_p}^{N_p,\n}(s)&=\zeta_{E_p}^{T_p,\n}(s)\zeta_{Z_p}^\s(s),\\
\nonumber	 \zeta_{G_p}^{H_p,\n}(s)&=\zeta_{E_p}^{H_p/Z_p,\n}(s)\zeta_{Z_p}^\s(s).
\end{align}
These follow immediately from the analysis in \cite[Section 2]{dSMS}. When the rank of $Z$ is 1, then $E$ is a plane crystallographic group, and we will be able to apply the results of \cite{Mc} (summarized in \cite[Section 4]{dSMS}). 

When computing local factors at primes $p|[G:N]$, we will sometimes use the method of $p$-adic integration with respect to the Mal'cev basis $\{x_1,x_2,x_3\}$; specifically, Corollary \ref{integral expression for zeta_G^N} and Proposition \ref{integral expression for zeta_H^H and zeta_G^H}. In this case, $\mc{M}_p^\s=\mc{M}_p^\n=\on{T}_3^+(\mb{Z}_p)$. 
The following lemma, whose verification is straightforward, will be used to translate the conditions defining the domains of integration  into {\em cone conditions}. 
\begin{lemma}\label{obtaining cone conditions abelian cas}
 An element $\x^\v\in N_p$ is in $\overline{\langle \x^{\t_1},\x^{\t_2},\x^{\t_3}\rangle}$ if and only if the following holds:
\begin{enumerate}
	\item $t_{11}|v_1$
	\item $t_{22}|-\frac{v_1}{t_{11}}t_{12}+v_2$
	\item $t_{33}|-\frac{-\frac{v_1}{t_{11}}t_{12}+v_2}{t_{22}}t_{23}-\frac{v_1}{t_{11}}t_{13}+v_3$ 
\end{enumerate}
\end{lemma} 

\subsection{Computing zeta functions of a family of groups including $\ms{G}_2$ and $\ms{B}_1$.} Given positive integers $a,b$, 
we shall compute the zeta functions of
\begin{align*}
 G=\left\langle \alpha,x_1,\ldots,x_a,x_{a+1},\ldots,x_{a+b}:\ \begin{array}{l} [x_i,x_j]=1\ \forall i,j,\ \alpha^2=x_{a+b},\\
 	 ^\alpha x_i=x_i^{-1}\ \mbox{for } i=1,\ldots,a,\\
 	 ^\alpha x_i=x_i\ \mbox{for } i=a+1,\ldots,a+b\end{array} \right\rangle. 
\end{align*}
This, of course, includes the cases $G=\ms{G}_2$ and $G=\ms{B}_1$.
Here $N=\langle x_1,\ldots,x_{a+b}\rangle$, $P\cong C_2$ is generated by the class of $\alpha$, and $Z=C_N(P)=\langle x_{a+1},\ldots,x_{a+b}\rangle$. Note that if $U\leq Z\leq V\leq N$, then $U$ and $V$ are $P$-submodules. 
We first prove a series of lemmas. When working inside $T$ or $Z$, we shall use additive notation.
\begin{lemma}\label{Hom from C_2}
	For $V\leq Z$ of finite index,  $\on{Hom}_P(T,Z/V)=\on{Hom}(C_2^a,Z/V)$, and this set has $[Z:2Z+V]$ elements.
	  \end{lemma} 
\begin{proof}
A group homomorphism $\varphi:T\to Z/V$ is of $P$-modules if and only if $\varphi(-x)=\varphi(x)$ for all $x\in T$, that is, $\varphi(2x)=0$ for all $x\in T$, or $\varphi(2T)=0$. This implies the first equality.	
Next, note that $\on{Hom}(C_2^a,Z/V)=\on{Hom}(C_2^a,((\frac{1}{2}V)\cap Z)/V)$, whose size is $[(\frac{1}{2} V) \cap Z:V]^a$. Thus, the second equality follows from the following calculation:
\begin{align}\label{calculus with indices} 
[(\frac{1}{2} V) \cap Z:V]=[ V\cap 2Z: 2V]=\frac{[V:2V]}{[V:V\cap 2Z]}=\frac{[Z:2Z]}{[(V+2Z):2Z]}=[Z:V+2Z].
\end{align} 
\end{proof}

\begin{lemma}\label{lemma on number of As with fixed B}
 Given $A\leq G$ of finite index such that $AN=G$, the subgroup $B:=A\cap N$ is normal of finite index in $G$ and $x_{a+b}\in (B\cap Z)+2Z$. If in addition $A$ is normal, then also $x_1^2,\ldots,x_a^2\in B$.

		 Conversely, given $B\leq N$ of finite index and normal in $G$ such that $x_{a+b}\in (B\cap Z)+2Z$, the set $\{A\leq G: AN=G,\ A\cap N=B\}$ has $[N:BZ][Z:2Z+(B\cap Z)]$ elements. If in addition $x_1^2,\ldots,x_a^2\in B$, then every $A$ in the latter set is normal in $G$.
	 \end{lemma}
 \begin{proof}
 	A subgroup $A\leq G$ such that $AN=G$ must contain $\alpha n$ for some $n\in N$. Then $(\alpha n)^2\in A$. It is easy to check that $(\alpha n)^2=x_{a+b}+2z$ for some $z\in Z$; thus $x_{a+b}+2z\in B\cap Z$, where $B=A\cap N$, or $x_{a+b}\in (B\cap Z)+2Z$. Note that $B$ is clearly normal, and if $A$ is also normal, then $[N,G]=\langle x_1^2,\ldots,x_a^2\rangle$ is included in $B$ by Lemma \ref{Gamma_c(G) is included in A}. This proves the first part of the lemma.
 	
 	We now fix $B\leq N$ of finite index and normal in $G$, and assume that $x_{a+b}+2z\in B\cap Z$ for some $z\in Z$, which will be fixed. 
 	Let $\{w_1,\ldots,w_s\}\subset \langle x_1,\ldots,x_a\rangle$ and $\{z_1,\ldots,z_t\}\subset Z$ be transversals to the cosets of $BZ$ in $N$ and to the cosets of $B\cap Z$ in $Z$, respectively. Then $\{w_iz_j: 1\leq i\leq s,\ 1\leq j\leq t\}$ is a transversal to the cosets of $B$ in $N$. A subgroup $A\leq G$ such that $AN=G$ and $A\cap N=B$ is of the form
 	$A=B\cup\alpha w_iz_j B$ for uniquely determined $w_i$ and $z_j$. Now $A=B\cup\alpha w_iz_j B$ is a subgroup if and only if $(\alpha w_iz_j)^2\in B$, that is $x_{a+b}z_j^2\in B$, or $x_{a+b}+2z_j\in B\cap Z$ in additive notation. In turn, this is equivalent to $2z_j-2z\in B\cap Z$, that is, $z_j-z\in (\frac{1}{2}(B\cap Z))\cap Z$. Therefore, only $[(\frac{1}{2}(B\cap Z))\cap Z:B\cap Z]$ of the $z_j$'s are allowed.
 	Since there is no restriction on $w_i$'s, the number of possibilities for $A$ is $[N:BZ][(\frac{1}{2}(B\cap Z))\cap Z:B\cap Z]$, which is equal to $[N:BZ][Z:2Z+(B\cap Z)]$ by (\ref{calculus with indices}). If in addition $[N,G]=\langle x_1^2,\ldots,x_a^2\rangle$ is included in $B$, then $G/B$ is abelian and hence $A$ is normal. This completes the proof of the lemma. 
 \end{proof} 
\begin{lemma}\label{lemma on the number of Bs for fixed U and v}
	Fix $U\leq T$ and $V\leq Z$ of finite index. Then $$
		|\{B\lhd G: (BZ)/Z=U,\ B\cap Z=V\}|=[Z: V+2Z]^a.$$
  If in addition $2T\subseteq U$, then 
  $$
		|\{B\lhd G: (BZ)/Z=U,\ B\cap Z=V,\ x_1^2,\ldots,x_a^2\in B\}|=[Z: V+2Z]^d,\quad \mbox{where } d=\dim_{\mb{F}_2}(U/2T).$$
\end{lemma}
\begin{proof}
 Let $\tilde{U}$ be the pre-image of $U$ in $N$.
	We are interested in the number of $P$-invariant complements of $Z/V$ in $\tilde{U}/V$. One of them is clearly $(U'+V)/V$, where $U'$ is the projection of $\tilde{U}$ onto $\langle x_1,\ldots,x_a\rangle$, and hence by Lemma \ref{number of F-invariant complements}  the number of them is $|\on{Hom}_{P}(\tilde{U}/Z,Z/V)|=|\on{Hom}_{P}(U,Z/V)|=|\on{Hom}_P(T,Z/V)|$. This, in turn, is equal to $[Z:V+2Z]^a$ by Lemma \ref{Hom from C_2}. This proves the first part of the lemma.
	
	 To show the second part, we use additive notation in $\tilde{U}/V$. There is no loss of generality if we assume that $U/2T$ is generated by the classes of $x_1,\ldots,x_d$. Let $\bar{x}_i$ denote the class of $x_i$ modulo $V$, $i=1,\ldots,a$.
	 We are interested in the number of $P$-invariant complements of $Z/V$ in $\tilde{U}/V$ that contain $2\bar{x}_1,\ldots,2\bar{x}_a$.  One of them is $B_0/V$, where $B_0=\langle \bar{x}_1,\ldots,\bar{x}_d,2\bar{x}_{d+1},\ldots,2\bar{x}_a\rangle$. To count how many there are, we follow the proof of Lemma \ref{number of F-invariant complements}.
	 A $P$-invariant complement is of the form $B_{\varphi}/V=\{x-\varphi(x): x\in B_0/V\}=\langle \bar{x}_1-\varphi(\bar{x}_1),\ldots,\bar{x}_d-\varphi(\bar{x}_d),2\bar{x}_{d+1}-\varphi(2\bar{x}_{d+1}),\ldots,2\bar{x}_a-\varphi(2\bar{x}_a)\rangle$ for a uniquely determined $\varphi\in \on{Hom}_P(B_0/V,Z/V)=\on{Hom}(B_0/(2B_0+V),Z/V)$. The last equality uses Lemma \ref{Hom from C_2}, which is possible since clearly $B_0/V\cong T$ as $P$-modules. Given $\varphi$ as above, we have $2\bar{x}_i=2\bar{x}_i- 2\varphi(\bar{x}_i)\in B_\varphi/V$ for $i=1,\ldots,d$. If it is also required that $2\bar{x}_{d+1},\ldots,2\bar{x}_a\in B_{\varphi}/V$, we need $\varphi(2\bar{x}_i)=0$ for all $i=d+1,\ldots,a$. Thus, the set of $P$-invariant complements of $Z/V$ in $\tilde{U}/V$ that contain $2\bar{x}_1,\ldots,2\bar{x}_a$ is in a bijection with $\on{Hom}(B_0/(2B_0+\langle 2x_{d+1},\ldots,2x_a\rangle+V),Z/V)$, that is, with $\on{Hom}(\frac{\langle x_1,\ldots,x_d\rangle}{2\langle x_1,\ldots,x_d\rangle},Z/V)$. By Lemma \ref{Hom from C_2}, this set has $[Z:V+2Z]^d$ elements. This completes the proof of the lemma.
\end{proof}

We now show a combinatorial lemma. 
Let $M$ be a finite dimensional vector space over $\mb{F}_2$. For a flag of subspaces $\beta: 0=F_0\subsetneqq F_1\subsetneqq\cdots\subsetneqq F_l=M$, the length $l$ is denoted by $l(\beta)$. The set of all flags is denoted by $\mc{F}(M)$.
\begin{lemma}
	Let $M$ be a non-zero finite dimensional vector space over $\mb{F}_2$. Then it holds that
 \begin{align}\label{formula 4}
	\sum_{W\leq M}\sum_{\beta\in\mc{F}(M/W)} (-1)^{l(\beta)}=0.
\end{align}
Let $k\geq 0$ and let $M$ be a $k$-dimensional vector space over $\mb{F}_2$. Then it holds that
\begin{align}\label{formula 2}
	\sum_{\beta\in\mc{F}(M)} (-1)^{l(\beta)}=(-1)^{k}2^{\binom{k}{2}}.
\end{align}
\end{lemma}
\begin{proof}
If $M$ has positive dimension, then any flag of $M$ of length $l$ produces a flag of length $l-1$ in $M/W$ for some  $0\neq W\leq M$, and conversely, any flag of $M/W$ of length $l-1$ extends to a flag of $M$ of length $l$. This implies that $\sum_{\beta\in\mc{F}(M)} (-1)^{l(\beta)}=	-\sum_{0\neq W\leq M}\sum_{\beta\in\mc{F}(M/W)} (-1)^{l(\beta)}$, which yields (\ref{formula 4}).

We now prove (\ref{formula 2}) by induction on $k$. The case $k=0$ is obvious, so we assume that $k>0$.
Let $A_k$ be the left-hand side and $B_k$ the right-hand side. Note that (\ref{formula 4}) implies that
$\sum_{i=0}^{k}A_i\binom{k}{i}_2=0$. Here $\binom{k}{i}_2$ denotes the 2-binomial coefficient, which expresses the number of $i$-dimensional subspaces in a $k$-dimensional vector space over $\mb{F}_2$. 
On the other hand, if we replace $t$ by $-1$ in the 2-binomial theorem
\begin{align}\label{formula 3}
	(1+t)(1+2t)\cdots (1+2^{k-1}t)=\sum_{i=0}^k 2^{\binom{i}{2}}\binom{k}{i}_2 t^i
\end{align}
we obtain that $0=\sum_{i=0}^k (-1)^i 2^{\binom{i}{2}}\binom{k}{i}_2=\sum_{i=0}^k B_k\binom{k}{i}_2$. 
Since $A_i=B_i$ for $i<k$ by the inductive hypothesis, and since $\binom{k}{k}_2=1$, it follows that $A_k=B_k$. This completes the induction.
 \end{proof}
\begin{lemma}\label{combinatorial lemma}
	Let $Z_1$ be a free abelian group of rank $b$ and let $Z_0\leq Z_1$ be a subgroup such that $2Z_1\subseteq Z_0$. Then
	\begin{align}\label{formula6} \sum_{\substack{V\leq Z_1\\ V+Z_0=Z_1}}[Z_1:V]^{-s}=\left(\prod_{i=1}^{b} \zeta(s-i+1)\right)  \left(\prod_{k=0}^{d-1} (1-2^{k-s}) \right), \quad\mbox{where } d=\dim_{\mb{F}_2} Z_1/Z_0.
		\end{align} 
\end{lemma}
\begin{proof}
 By  (\ref{formula 4}), the left-hand side of (\ref{formula6}) becomes
 \begin{align*}
 \sum_{V\leq  Z_1}[Z_1:V]^{-s}\left(\sum_{Z_0+V\leq W\leq Z_1}\sum_{\beta\in\mc{F}(Z_1/W)} (-1)^{l(\beta)}\right)
=\sum_{Z_0\leq W\leq Z_1}\sum_{V\leq W}[Z:W]^{-s}[W:V]^{-s}\left(\sum_{\beta\in\mc{F}(Z_1/W)} (-1)^{l(\beta)}\right).	
 \end{align*}
By (\ref{formula 2}), the latter becomes
\begin{align*}
		\prod_{i=1}^{b}\zeta(s-i+1)\left(\sum_{Z_0\leq W\leq Z_1}[Z:W]^{-s}(-1)^{\dim Z_1/W} 2^{\binom{\dim Z_1/W}{2}}\right)=\prod_{i=1}^{b}\zeta(s-i+1)\left(\sum_{j=0}^d2^{-js}\binom{d}{j}_2 (-1)^j 2^{\binom{j}{2}}\right).
\end{align*}
Finally, by using (\ref{formula 3}) with $k=d$ and $t=-2^{-s}$, the last expression becomes the right-hand side of (\ref{formula6}).
\end{proof}
 
\begin{proposition} 
	The partial zeta functions of $G$ with respect to $N$ are as follows:
	\begin{align*}
		\zeta_{G}^{G,\s}(s)&=\prod_{i=1}^a\zeta(s-i)\times\prod_{j=1}^{b}\zeta(s-j+1)\times \left(\sum_{k=1}^{b}  2^{(b-k)(a+1-s)} \binom{b-1}{k-1}_2 (1-2^{-s})\cdots (1-2^{k-1-s})\right),\\
	 \zeta_G^{G,\n}(s)&=\prod_{i=1}^{a}\zeta(s-i+1)\times\prod_{j=1}^{b}\zeta(s-j+1)
	 \times\left(\sum_{k=0}^b 2^{(b-k)(s-a)}\binom{b}{k}_2(1-2^{-s})\cdots (1-2^{k-1-s}),\right)\\
	 \zeta_{G}^{N,\n }(s)&=\prod_{j=1}^{b}\zeta(s-j+1)\times \left(\sum_{i=0}^a2^{(a-i)(1-s)}\binom{a}{i}_2\left(\sum_{k=1}^{b}2^{(b-k)(i+1-s)}\binom{b-1}{k-1}_2(1-2^{-s})\cdots (1-2^{k-1-s})\right) \right)
	\end{align*} 
\end{proposition}
\begin{proof} The first equality below uses Lemma \ref{lemma on number of As with fixed B}, the second one Lemma \ref{lemma on the number of Bs for fixed U and v}, and the last one Lemma \ref{combinatorial lemma}.
	\begin{align*}
		\zeta_{G}^{G,\s}(s)&=\sum_{\substack{B\lhd G: B\subseteq N\\ x_{a+b}\in (B\cap Z)+2Z}}[N:B]^{-s}[N:BZ][Z:(B\cap Z)+2Z] \\
		&=\sum_{\substack{U\leq T, V\leq Z:\\ x_{a+b}\in V+2Z}}[T:U]^{-s}[Z:V]^{-s} [Z:V+2Z]^a[T:U][Z:V+2Z]\\
		&=\zeta_{T}^\s(s-1)\left(\sum_{2Z+\langle x_{a+b}\rangle \leq Z_1\leq Z} [Z:Z_1]^{a+1-s} \sum_{\substack{V\leq Z_1\\ V+2Z=Z_1}}[Z_1:V]^{-s} \right)\\
        &=\prod_{i=1}^a\zeta(s-i)\sum_{k=1}^b\left(\sum_{\substack{2Z+\langle x_{a+b}\rangle \leq Z_1\leq Z\\ [Z:Z_1]=2^{b-k}}} 2^{(b-k)(a+1-s)} \sum_{\substack{V\leq Z_1\\ V+2Z=Z_1}}[Z_1:V]^{-s} \right)\\
		&=\prod_{i=1}^a\zeta(s-i)\times\prod_{j=1}^{b}\zeta(s-j+1) \times\left(\sum_{k=1}^{b} 2^{(b-1-k)(a+1-s)}\binom{b-1}{k-1}_2(1-2^{-s})\cdots (1-2^{k-1-s})\right). 
	\end{align*}
The proofs of the expressions for $\zeta_{G}^{N,\n}(s)$ and $\zeta_G^{G,\n}(s)$ are similar. In the first case we do not have to add the condition $x_{a+b}\in V+2Z$ or the factor $[T:U][Z:V+2Z]$ in the second line. In the second case, we have to add the condition $2T\leq U$ and replace $[Z:V+2Z]^a$ by $[Z:V+2Z]^{\dim_{\mb{F}_2}(U/T)}$ in the second line.
\end{proof}

\subsection{Computing the zeta functions of $\ms{G}_3$, $\ms{G}_4$, and $\ms{G}_5$} 
Let $G$ be any of the groups $\ms{G}_3$, $\ms{G}_4$, or $\ms{G}_5$. We set $d=|G/N|$. Note that $Z=\langle x_1\rangle$ and that $E$ is the plane crystallographic group $\mathbf{pd}$.

\subsubsection{Local factors of $\zeta_{G}^{G,\s}(s)$} If $p\nmid d$, then  by (\ref{general remark}) we have $\zeta_{G_p}^{G_p,\s}(s)=\zeta_{E_p}^{E_p,\s}(s)\zeta_{Z_p}^\s(s)$, and
according to \cite[5.10, 5.13, 5.16]{Mc}, the latter is equal to $\zeta_p(s-1)L_p(s-1,\chi_d)\zeta_p(s)$ if $d\in\{3,4\}$, and to $\zeta_p(s-1)L_p(s-1,\chi_3)\zeta_p(s)$ if $d=6$.

 Assume now that $p\mid d$. We claim that if $A\leq G_p$ is open and satisfies $AN_p=G_p$, then $Z_p\subset A$.
Indeed, $A$ contains $\alpha \x^\v$ for some $\v\in\mb{Z}_p^3$. Then $(\alpha \x^\v)^d=x_1^{dv_{1}+1}\in A$, and hence $x_1\in A$ since $dv_{1}+1\in\mb{Z}_p^*$. This proves the claim.
It follows now that $\zeta_{G_p}^{G_p,\s}(s)=\zeta_{E_p}^{E_p,\s}(s)$. According to  \cite[5.10, 5.13, 5.16]{Mc}, this series is equal to $\zeta_p(s-1)L_p(s-1,\chi_q)$, where $q=d$ if $d\in\{3,4\}$ and $q=3$ if $d=6$.

\subsubsection{Local factors of $\zeta_{G}^{N,\n}(s)$} 
If $p\nmid d$, then by (\ref{general remark}) we have $\zeta_{G_p}^{N_p,\n}(s)=\zeta_p(s)\zeta_{E_p}^{T_p,\n}(s)$, and  according to \cite[6.10, 6.13, 6.16]{Mc} the latter is equal to  $\zeta_p(s)L_p(s,\chi_d)\zeta_p(s)$ if $d\in\{3,4\}$ and to $\zeta_p(s)L_p(s,\chi_3)\zeta_p(s)$ if $d=6$.

Assume now that $p\mid d$. By Corollary \ref{integral expression for zeta_G^N}, we have 
\begin{align*}
	\zeta_{G_p}^{N_p,\n}(s)=\frac{1}{(1-p^{-1})^3}\int_{\mc{T}}|t_{11}|_p^{s-1}|t_{22}|_p^{s-2}|t_{33}|_p^{s-3}d\mu, 
\end{align*}
where $\mc{T}\subset \on{T}_3^+(\mb{Z}_p)$ is defined by the following equivalent conditions:  
\begin{align*} 
&	{}^\alpha(\x^{\t_1}) \x^{-\t_1},\ ^\alpha(\x^{\t_2}),\  ^\alpha(\x^{\t_3}) \in\overline{\langle \x^{\t_1},\x^{\t_2},\x^{\t_3}\rangle}\\
\Leftrightarrow\quad & 
x_2^{-t_{13}-t_{12}} x_3^{t_{12}-(r+1)t_{13}},\ x_2^{-t_{23}}x_3^{t_{22}-r t_{23}} \in\overline{\langle \x^{\t_1},\x^{\t_2},\x^{\t_3}\rangle},
	\end{align*}  
 where $r=1$ if $d=3$, $r=0$ if $d=4$, and $r=-1$ if $d=6$. Therefore, by Lemma \ref{obtaining cone conditions abelian cas}, the conditions defining $\mc{T}$ are
\begin{align*}
\begin{array}{lll} (1)\quad t_{22}|-t_{12}-t_{13},& (2)\quad  t_{33}|\frac{t_{12}+t_{13}}{t_{22}}t_{23}+t_{12}-(r+1)t_{13},& (3)\quad  t_{22}|t_{23},\\
(4)\quad  t_{33}|\frac{t_{23}}{t_{22}}t_{23}+t_{22}-r t_{23},& (5)\quad t_{22}|-t_{33},& (6)\quad t_{33}|\frac{t_{33}}{t_{22}}t_{23}. \end{array} 
\end{align*}
Note that (3) implies (6), and (4) can be written as $\frac{t_{33}}{t_{22}}|(\frac{t_{23}}{t_{22}})^2-r\frac{t_{23}}{t_{22}} +1$. 

{\em Assume first that $d\in\{3,4\}$.}  One easily checks that the equation $T^2-r T+1\equiv 0\mod p$, in $\mb{Z}_p$, has $-1$ as a unique solution modulo $p$, and there are no solutions in $\mb{Z}_p$ for $T^2-r T+1\equiv 0\mod p^2$. Therefore, (4) splits the domain of integration $\mc{T}$ as a disjoint union $\mc{T}^a\cup\mc{T}^b$ according to the cases (4a) $|t_{22}|_p=|t_{33}|_p$ and (4b) $|pt_{22}|_p=|t_{33}|_p$ and $p|\frac{t_{23}}{t_{22}}+1$. 

If we assume (4a), then (1) and (3) reduce (2) to $t_{33}|t_{12}-(r+1)t_{13}$, and then (1) reduces (2) again to $t_{33}|(r+2)t_{13}$. Note that $r+2=p$. The conditions defining $\mc{T}^a$ are, therefore, $
	t_{22}|-t_{12}-t_{13}$,  $t_{22}|pt_{13}$, $t_{22}|t_{23}$ and $|t_{22}|_p=|t_{33}|_p$, and the integral over $\mc{T}^a$ becomes
\begin{align*}
	\frac{1}{(1-p^{-1})^3}\int_{\mc{T}^a} |t_{11}|_p^{s-1}|t_{22}|_p^{2s-5}d\mu=\frac{\zeta_p(s)}{1-p^{-1}}\int_{t_{22}|pt_{13}}|t_{22}|_p^{2s-2}d\mu=\zeta_p(s)\left(1+p^{1-2s}\zeta_p(2s)\right),
\end{align*}
where in the first equality we integrated with respect to $t_{11}$ and applied Proposition \ref{proceso de simplifiacion de integrales} with the pivots $t_{12}, t_{23}, t_{33}$.

If we now assume (4b), then (3) and (5) are redundant, and since $r+2=p$, we see that (2) can be written as $t_{33}|(\frac{t_{23}}{t_{22}}+1)t_{12}+(\frac{t_{23}}{t_{22}}+1-p)t_{13}$. Using (1) and the assumption that $|t_{33}|_p=|pt_{22}|_p$ and $p|\frac{t_{23}}{t_{22}}+1$, we can reduce (2) to $t_{33}|pt_{13}$, which is equivalent to $t_{22}|t_{13}$. Thus, $\mc{T}^b$ is defined by the conditions
$t_{22}|-t_{12}-t_{13}$, $t_{22}|t_{13}$, $ pt_{22}|t_{23}+t_{22}$ and $|pt_{22}|_p=|t_{33}|_p
$, and the integral over $\mc{T}^b$ becomes 
\begin{align*}
\frac{1}{(1-p^{-1})^3}\int_{\mc{T}^b}|t_{11}|_p^{s-1}|t_{22}|_p^{2s-5}p^{3-s}d\mu=\frac{1}{(1-p^{-1})^2}\int |t_{11}|_p^{s-1}|t_{22}|_p^{2s-1}p^{1-s}d\mu=p^{1-s}\zeta_p(s)\zeta_p(2s),
\end{align*}
where, in the first equality, we used Proposition \ref{proceso de simplifiacion de integrales} with the pivots $t_{12}, t_{13}, t_{23}, t_{33}$.

 We conclude that
$$
	\zeta_{G_p}^{N_p,\n}(s)=\zeta_p(s)\left(1+p^{1-2s}\zeta_p(2s)+p^{1-s}\zeta_p(2s)\right)=\zeta_p^2(s)(1+(p-1)p^{-s})\quad \mbox{for } d\in\{3,4\}\ \mbox{and } p|d.
$$

{\em Assume now that $d=6$ (and hence $r=-1$).} If $p=3$, then the analysis is the same as above except that the condition $t_{22}|3t_{13}$ in the description of $\mc{T}^a$ must be replaced by $t_{22}|t_{13}$ since this time $r+2=1$, and the condition (4b) is now $|3t_{22}|_3=|t_{33}|_3$ and $3|\frac{t_{23}}{t_{33}}-1$ since $1\in\mb{Z}_3$ is the unique solution modulo $3$ of $T^2+T+1\equiv 0\mod 3$. Thus, $\mc{T}^a$ is defined by $
	t_{22}|-t_{12}-t_{13}$, $t_{22}|t_{13}$, $t_{22}|t_{23}$ and $|t_{22}|_3=|t_{33}|_3$,
and the integral over $\mc{T}^a$ becomes
\begin{align*}
	\frac{1}{(1-3^{-1})^3}\int_{\mc{T}^a} |t_{11}|_3^{s-1}|t_{22}|_3^{2s-5}d\mu&=\frac{1}{(1-3^{-1})^2}\int|t_{11}|_3^{s-1}|t_{22}|_3^{2s-1}d\mu=\zeta_3(s)\zeta_3(2s),
\end{align*}
where in the first equality we used Proposition \ref{proceso de simplifiacion de integrales} with the pivots $t_{12},t_{13},t_{23},t_{33}$. 

If we assume (4b), then (3) and (6) are redundant, and (2) can be written as $t_{33}|(\frac{t_{23}}{t_{22}}+1)t_{12}+\frac{t_{23}}{t_{22}}t_{13}$. Since $\frac{t_{23}}{t_{22}}$ is a unit, we can use the previous condition and $t_{22}|t_{33}$ to reduce (1) to $t_{22}|t_{12}$. Thus, the conditions defining $\mc{T}^b$ are $
t_{22}|t_{12}$, $3t_{22}|(\frac{t_{23}}{t_{22}}+1)t_{12}+\frac{t_{23}}{t_{22}}t_{13}$, $3t_{33}|t_{23}-t_{33}$ and|$3t_{22}|_3=|t_{33}|_3$,
and the integral over $\mc{T}^b$ becomes
\begin{align*}
\frac{1}{(1-3^{-1})^3}\int_{\mc{T}^b}|t_{11}|_3^{s-1}|t_{22}|_3^{2s-5}3^{3-s}d\mu=\frac{1}{(1-3^{-1})^2}\int |t_{11}|_3^{s-1}|t_{22}|_3^{2s-1}3^{-s}d\mu=3^{1-s}\zeta_3(s)\zeta_3(2s),
\end{align*}
where in the first equality we used Proposition \ref{proceso de simplifiacion de integrales} with the pivots $t_{12}, t_{23}, t_{13}, t_{33}$. 

We conclude that 
$$
	\zeta_{G_3}^{N_3,\n}(s)=\zeta_3(s)\zeta_3(2s)+3^{-s}\zeta_3(s)\zeta_3(2s)=\zeta_3(s)^2.$$
	
Assume finally that $p=2$. Since there are no solutions in $\mb{Z}_2$ for $T^2+T+1\equiv 0\mod 2$, (4) is equivalent to $|t_{33}|_2=|t_{22}|_2$ and (2) is reduced to $t_{22}|t_{12}$. Thus, the conditions defining $\mc{T}$ are
	$t_{22}|-t_{12}-t_{13}$, $t_{22}|t_{12}$, $t_{22}|t_{23}$ and $|t_{22}|_2=|t_{33}|_2$,
and the integral becomes
\begin{align*}
	\zeta_{G_2}^{N_2,\n}(s)=\frac{1}{(1-2^{-1})^3}\int_{\mc{T}}|t_{11}|_2^{s-1}|t_{22}|_2^{2s-5}d\mu=\frac{1}{(1-2^{-1})^2}\int |t_{11}|_2^{s-1}|t_{22}|_2^{2s-1}d\mu=\zeta_2(s)\zeta_2(2s),
\end{align*}
where in the first equality we used Proposition \ref{proceso de simplifiacion de integrales} with the pivots $(t_{13},t_{12}, t_{23},t_{33})$.

\subsubsection{Local factors of $\zeta_{G}^{G,\n}(s)$} If $p\nmid d$, then
$[N_p,G_p]=\overline{\langle x_2^{-1}x_3, x_2^{-1}x_3^{-r-1}\rangle}=\overline{\langle x_2^{-1}x_3, x_3^{-r-2}\rangle}=\overline{\langle x_2, x_3\rangle}$, where the last equality holds since $r+2\in\mb{Z}_p^*$. As $N_p/\overline{\langle x_2, x_3\rangle}\cong\mb{Z}_p$, we have $\zeta_{G_p}^{G_p,\n}(s)=\zeta_p(s)$ by Proposition \ref{simplification for normal extension zeta function}.

Assume now that $p\mid d$. Arguing as in the computation of $\zeta_{G_p}^{G_p,\s}(s)$ (with $p\mid d$), we find that $\zeta_{G_p}^{G_p,\n}(s)=\zeta_{E_p}^{E_p,\n}(s)$, and according to the results of \cite[6.10, 6.13, 6.16]{Mc}, this is equal to $1+p\cdot p^{-s}$ if $d=3$ or $4$, and to $1$ if $d=6$.

\subsubsection{Local factors of $\zeta_{G}^{H,\n}(s)$ for $d=4$ and $H=\langle \alpha^2,x_1,x_2,x_3\rangle$}\label{G_4, zeta_G^H, normal} 
Note that $H=\ms{G}_2$.  If $p\neq 2$, then $[N_p,H_p]=\overline{\langle x_2^2,x_2^2\rangle}=\overline{\langle x_2,x_2\rangle}$, and hence $N_p/\overline{\langle x_2,x_2\rangle}\cong\mb{Z}_p$. Therefore, by Proposition \ref{simplification for normal extension zeta function},  $\zeta_{G_p}^{H_p,\lhd}(s)=\zeta_p(s)$. 

Assume that $p=2$. An open subgroup $A\lhd G_2$ such that $AN_2=H_2$ must include $Z_2$. Indeed, $A$ contains $\alpha^2\x^\v$ for some $\v\in\mb{Z}_2^3$, and hence $(\alpha^2 \x^\v)^2=  x_1^{2v_1+1}\in A$, or $x_1\in A$ since $2t_1+1\in \mb{Z}_2^*$. It follows that $\zeta_{G_2}^{H_2,\n}(s)=\zeta_{E_2}^{H_2/Z_2,\n}(s)$, and this is equal to $1+2\cdot 2^{-s}+2\cdot 4^{-s}$, according to the results of \cite[6.10]{Mc}.

\subsubsection{Local factors of $\zeta_{G}^{H,\n}(s)$ for $d=6$ and $H=\langle \alpha^3,x_1,x_2,x_3\rangle$}
 Note that $H=\ms{G}_2$. Arguing as in the previous paragraph we find that $\zeta_{G_p}^{H_p,\n}(s)=\zeta_p(s)$ if $p\neq 2$, and that  $\zeta_{G_2}^{H_2,\n}(s)=\zeta_{E_2}^{H_2/Z_2,\n}(s)$, which is equal to $1+4^{-s}$, according to the results of \cite[6.16]{Mc}.

\subsubsection{Local factors of $\zeta_{G}^{K,\n}(s)$ for $d=6$ and $K=\langle \alpha^2,x_1,x_2,x_3\rangle$}
 Note that $K=\ms{G}_3$. Arguing as in the case $d=3$ and $K=G$, we obtain that $\zeta_{G_p}^{K_p,\n}(s)=\zeta_p(s)$ if $p\neq 3$, and that $\zeta_{G_3}^{K_3,\n}(s)=\zeta_{E_3}^{K_3/Z_3,\n}(s)$, which is equal to $1+3^{-s}$, according to the results of \cite[6.16]{Mc}.

\medskip
 
\subsection{Computing the zeta functions of $G=\ms{G}_6$.}
Note that $Z=0$ and that $N=\langle x_1\rangle\oplus\langle x_2\rangle\oplus\langle x_3\rangle$ is a decomposition of $N$ into a  direct sum of irreducible $\mb{Z}[P]$-modules that are pairwise non-isomorphic.

\subsubsection{Local factors of $\zeta_{G}^{G,\s}(s)$}\label{G_6, zeta_G^G, subgroup} If $p\neq 2$, then by using the above decomposition of $N$ and the results of \cite[Section 2]{dSMS}, we obtain that $\zeta_{G_p}^{G_p,\s}(s)=\zeta_p(s-1)\zeta_p(s-1)\zeta_p(s-1)$.  

We now compute $\zeta_{G_2}^{G_2,\s}(s)$. An open subgroup $A\leq G_2$ such that $AN_2=G_2$ must contain $\alpha\x^\mathbf{u}$,  $\beta \x^\v$ and  $\alpha\beta \x^\w$ for some $\mathbf{u},\v,\w\in\mb{Z}_2^3$, whence $A$ contains the squares of these elements, which are $x_1^{2u_1+1}$, $x_2^{2v_2+1}$ and $x_3^{2w_3+1}$. Since $2u_1+1, 2v_2+1, 2w_3+1\in\mb{Z}_2^*$, it follows that $A$ includes $x_1,x_2,x_3$ and hence also $N_2$. This implies that $A=G_2$. We conclude that $\zeta_{G_2}^{G_2,\s}(s)=1$.

\subsubsection{Local factors of $\zeta_{G}^{N,\s}(s)$:}\label{G_6, zeta_G^N, normal}
If $p\neq 2$, arguing as in the calculation of $\zeta_{G_p}^{G_p,\s}(s)$, we find that $\zeta_{G_p}^{N_p,\n}(s)=\zeta_p(s)^3$. 

Assume now that $p=2$. By Corollary \ref{integral expression for zeta_G^N}, we have
\begin{align*}
	\zeta_{G_2}^{N_2,\n}(s)=\frac{1}{(1-2^{-1})^3}\int_{\mc{T}} |t_{11}|_2^{s-1}|t_{22}|_2^{s-2}|t_{33}|_2^{s-3}d\mu, 
\end{align*}
where $\mc{T}\subset \on{T}_3^+(\mb{Z}_2)$ is defined by the following equivalent conditions:
\begin{align*}
	&	{}^\alpha (\x^{\t_1})\x^{-\t_1},\ ^\alpha(\x^{\t_2}),\ ^\alpha(\x^{\t_3}),\ ^\beta(\x^{\t_1})\x^{\t_1},\ ^\beta(\x^{\t_2})\x^{-\t_2},\ ^\beta(\x^{\t_3})\in \overline{\langle \x^{\t_1},\x^{\t_2},\x^{\t_3}\rangle}\\
	\Leftrightarrow\quad & 	x_2^{-2t_{12}}x_3^{-2t_{13}},\ \x^{-\t_2},\ \x^{-\t_3},\ x_2^{2t_{12}},\ x_3^{-2t_{23}},\ \x^{-\t_3}\in \overline{\langle \x^{\t_1},\x^{\t_2},\x^{\t_3}\rangle},\\
	\Leftrightarrow\quad  & x_3^{-2t_{13}},\ x_2^{2t_{12}},\ x_3^{-2t_{23}} \in \overline{\langle \x^{\t_1},\x^{\t_2},\x^{\t_3}\rangle}.
\end{align*}
According to Lemma \ref{obtaining cone conditions abelian cas}, these conditions are equivalent to
\begin{align*}
	\begin{array}{llll}
		(1)\quad t_{22}|2t_{12},& (2)\quad t_{33}|-\frac{2t_{12}}{t_{22}}t_{23},& (3)\quad t_{33}|2t_{13}, & (4)\quad t_{33}|2t_{23}. 
	\end{array}
\end{align*}
We express $\mc{T}$ as a disjoint union $\mc{T}^a\cup\mc{T}^b$ according to the cases (4a) $t_{33}|t_{23}$ and (4b) $|t_{33}|_2=|2t_{23}|_2$. 

In the first case, (2) is redundant and hence $\mc{T}^a$ is defined by the conditions $
	t_{22}|2t_{12},\ t_{33}|2t_{13},\ t_{33}|t_{23}$. Thus,
\begin{align*}
\frac{1}{(1-2^{-1})^3}\int_{\mc{T}^a}|t_{11}|_2^{s-1}|t_{22}|_2^{s-2}|t_{33}|_2^{s-3}d\mu=\frac{\zeta_2(s)}{(1-2^{-1})^2}\int_{\substack{t_{22}|2t_{12}\\ t_{33}|2t_{13}}} |t_{22}|^{s-2}|t_{33}|_2^{s-2}d\mu=\zeta_2(s)(1+2^{1-s}\zeta_2(s))^2.
\end{align*} 

Assume now (4b). Note that (2) becomes equivalent to $t_{22}|t_{12}$, whence (1) is redundant. Note also that (3) is equivalent to $t_{23}|t_{13}$. Therefore, the conditions defining $\mc{T}^b$ are $
	t_{22}|t_{12}$, $t_{23}|t_{13}$ and $|t_{33}|_2=|2t_{23}|_2$,
and the integral over $\mc{T}^b$ becomes
\begin{align*}
\frac{1}{(1-2^{-1})^3}\int_{\mc{T}^b}|t_{11}|_2^{s-1}|t_{22}|_2^{s-2} 2^{3-s}|t_{23}|_2^{s-3}d\mu=\frac{\zeta_2(s)2^{1-s}}{(1-2^{-1})^2}\int |t_{22}|_2^{s-1}|t_{23}|_2^{s-1}d\mu=2^{1-s}\zeta_2(s)^3,	
\end{align*}
where in the second equality we used Proposition \ref{proceso de simplifiacion de integrales} with the pivots $t_{12},t_{13},t_{33}$.

We conclude that
	$$\zeta_{G_2}^{N_2,\n}(s)=\zeta_2(s)\left((1+2^{1-s}\zeta_2(s))^2+2^{1-s}\zeta_2(s)^2\right)=(1+4\cdot 2^{-s}+4^{-s})\zeta_2(s)^3.$$

\subsubsection{Local factors of $\zeta_G^{G,\n}(s)$} If $p\neq 2$, then $[N_p,G_p]=\overline{\langle x_1^2,x_2^2,x_3^2\rangle}=\overline{\langle x_1,x_2,x_3\rangle}$, whence $\zeta_{G_p}^{G_p,\n}(s)=1$ by Proposition \ref{simplification for normal extension zeta function}. 
If $p=2$, then $\zeta_{G_2}^{G_2,\n}(s)=1$ since $\zeta_{G_2}^{G_2,\s}(s)=1$, as shown in \ref{G_6, zeta_G^G, subgroup}.

\subsubsection{Local factors of $\zeta_{G}^{H,\n}(s)$ for  $H=\langle \alpha,x_1,x_2,x_3\rangle$} Note that $H\cong\ms{G}_2$. Arguing as in \ref{G_4, zeta_G^H, normal}, we find that $\zeta_{G_p}^{H_p,\lhd}(s)=\zeta_{p}(s)$, and that $\zeta_{G_2}^{H_2,\n}(s)=\zeta_{G_2/\overline{\langle x_1\rangle}}^{H_2/\overline{\langle x_1\rangle},\n}(s)$. One easily checks that $G/{\langle x_1\rangle}\cong \textbf{p2gg}$; therefore, the latter series is equal to $1+2\cdot 2^{-s}$, according to the results of \cite[6.7]{Mc}.

\subsubsection{Calculation of $\zeta_{G}^{K,\n}(s)$ and $\zeta_{G}^{L,\n}(s)$ for $K=\langle \beta,x_1,x_2,x_3\rangle$ and $L=\langle \alpha\beta,x_1,x_2,x_3\rangle$:}
It is easy to check that $(\alpha,\beta,x_1,x_2,x_3)\mapsto (\beta, \alpha\beta,x_2,x_3,x_1^{-1})$ defines an automorphism of $G$ that sends $H$ onto $K$, and $K$ onto $L$. Thus, $\zeta_{G}^{K,\n}(s)=\zeta_{G}^{L,\n}(s)=\zeta_{G}^{H,\n}(s)$. This series was computed in the previous paragraph. 
\subsection{Computing the zeta functions of $G=\ms{B}_2$.}
  We call $y_1=x_3$, $y_2=x_1$ and $y_3=x_1x_2$, so that now the relations are $\varepsilon^2=y_2$, $ ^\varepsilon y_3=y_3$ and $ ^\varepsilon y_1=y_1^{-1}y_3$.
 We will use this presentation. Note that $Z=\langle y_2,y_3\rangle$, and hence $T\cong\mb{Z}$.

\subsubsection{Local factors of $\zeta_G^{G,\s}(s)$} If $p\neq 2$, then by (\ref{general remark}) we have $\zeta_{G_p}^{G_p,\s}(s)=\zeta_{Z_p}(s)\zeta_p(s-1)=\zeta_p(s)\zeta_p(s-1)^2$. 

Assume now that $p=2$. By Proposition \ref{integral expression for zeta_H^H and zeta_G^H} with $H=G$,  we have 
	$$\zeta_{G_2}^{G_2,\s}(s)=\frac{1}{(1-2^{-1})^3}\int_\mc{T}|t_{11}|_2^{s-2}|t_{22}|_2^{s-3}|t_{33}|_2^{s-4}d\mu,$$ 
	where $\mc{T}\subset  \on{T}_3^+(\mb{Z}_2)\times\mb{Z}_2^3$ is the set of pairs $(\t,\v)$ satisfying the following equivalent conditions:
	 \begin{align*}   &{}^{\varepsilon}(\y^{\t_1}) \y^{\t_1},\  ^{\varepsilon}(\y^{\t_2}),\ ^{\varepsilon}(\y^{\t_3}),  (\varepsilon\y^\v)^2\in\overline{\langle\y^{\t_1},\y^{\t_2},\y^{\t_3}\rangle}\\
	 \Leftrightarrow\quad & y_2^{2t_{12}}y_3^{t_{11}+2t_{13}}, y_2^{2v_2+1}y_3^{2v_3+v_1}\in \overline{\langle\y^{\t_1},\y^{\t_2},\y^{\t_3}\rangle}.
	 \end{align*}
  According to Lemma \ref{obtaining cone conditions abelian cas}, these conditions are equivalent to
	\begin{align*}
		(1)\ t_{22}|2t_{12},&& (2)\ t_{33}|-\frac{2t_{12}}{t_{22}}t_{23}+t_{11}+2t_{13},&& (3)\ t_{22}|2v_2+1,&& (4)\ t_{33}|-\frac{2v_1+1}{t_{22}}t_{23}+2v_3+v_1.
	\end{align*}
	Note that (3) is equivalent to $t_{22}\in \mb{Z}_2^*$, whence (1) is redundant. We now express $\mc{T}$ as a disjoint union $\mc{T}^a\cup\mc{T}^b$ according to the cases: a) $t_{33}\in\mb{Z}_2^*$ and b) $t_{33}\in 2\mb{Z}_2$. 
	
	Note that $\mc{T}^a$ is defined by $t_{22},t_{33}\in\mb{Z}_p^*$, and hence
	\begin{align*}
		\frac{1}{(1-2^{-1})^3}\int_{\mc{T}^a}|t_{11}|_2^{s-2}|t_{22}|_2^{s-3}|t_{33}|_2^{s-4}d\mu=\zeta_2(s-1).
	\end{align*}
	
	In case b), (2) can be written as $\frac{t_{33}}{2}|-\frac{t_{12}}{t_{22}}t_{23}+\frac{t_{11}}{2}+t_{13}$, and it implies that $t_{11}\in 2\mb{Z}_2$. By using  Proposition \ref{proceso de simplifiacion de integrales} with the pivots $(t_{13},v_1)$ the integral over $\mc{T}^b$ becomes
	\begin{align*}
		\frac{1}{(1-2^{-1})^3}\int_{\mc{T}^b}|t_{11}|_2^{s-2}|t_{33}|_2^{s-4}d\mu=\frac{1}{(1-2^{-1})^2}\int_{t_{11},t_{33}\in 2\mb{Z}_2}|t_{11}|_2^{s-2}2|t_{33}|_2^{s-2}d\mu=2^{3-2s}\zeta_2(s-1)^2.
	\end{align*}
	
	We conclude that 
	$
	\zeta_{G_2}^{G_2,\s}(s)=\zeta_2(s-1)(1+2^{3-2s}\zeta_2(s-1))=\zeta_2(s-1)^2(1-2^{-s}+8\cdot 4^{-s}).$

\subsubsection{Local factors of $\zeta_{G}^{N,\n}(s)$}
 If $p\neq 2$, then by (\ref{general remark}) we have $\zeta_{G_p}^{N_p,\n}(s)=\zeta_{Z_p}^\s(s)\zeta_p(s)=\zeta_p(s)^2\zeta_p(s-1)$. 

Assume now that $p=2$. By Corollary \ref{integral expression for zeta_G^N}, we have
$$\zeta_{G_2}^{N_2,\n}(s)=\frac{1}{(1-2^{-1})^3}\int_\mc{T}|t_{11}|_2^{s-1}|t_{22}|_2^{s-2}|t_{33}|_2^{s-3}d\mu,$$
 where 
 $$\mc{T}=\{\t\in \on{T}_3^+(\mb{Z}_2): {}^{\varepsilon}(\y^{\t_1})\y^{\t_1},\ ^{\varepsilon}(\y^{\t_2}),\ ^{\varepsilon}(\y^{\t_3})\in \overline{\langle\y^{\t_1},\y^{\t_2},\y^{\t_3}\rangle}\}=\{\t\in \on{T}_3^+(\mb{Z}_2): y_2^{2t_{12}}y_3^{t_{11}+2t_{13}}\in \overline{\langle\y^{\t_1},\y^{\t_2},\y^{\t_3}\rangle}\}.$$
  By Lemma \ref{obtaining cone conditions abelian cas}, the conditions defining $\mc{T}$ are equivalent to
	\begin{align*}
		(1)\quad t_{22}|2t_{12}, && 
		(2)\quad t_{33}|-\frac{2t_{12}}{t_{22}}t_{23}+t_{11}+2t_{13}.
	\end{align*}
 We can express $\mc{T}$ as a disjoint union $\mc{T}=\mc{T}^a\cup\mc{T}^b$ according to the cases: (1a) $|t_{22}|_2=|2t_{12}|_2$ and (1b) $t_{22}|t_{12}$.
 
  In the first case the coefficient of $t_{23}$ in (2) is a unit, so the integral over $\mc{T}^a$ becomes
	\begin{align*}
		\frac{2^{2-s}}{(1-2^{-1})^3}\int_{\mc{T}^a}|t_{11}|_2^{s-1}|t_{12}|_2^{s-2}|t_{33}|_2^{s-3}d\mu=\frac{2^{2-s}}{(1-2^{-1})^2}\int |t_{11}|_2^{s-1} 2^{-2}|t_{12}|_2^{s-1}|t_{33}|_2^{s-2}d\mu=2^{-s}\zeta_2(s)^2\zeta_2(s-1), 
	\end{align*}
where in the first equality we used Proposition \ref{proceso de simplifiacion de integrales} with the pivots $t_{22},t_{23}$.

	Assume now (1b). We express $\mc{T}^b$ as a disjoint union $\mc{T}^b=\mc{T}^{b,i}\cup\mc{T}^{b,ii}$ according to the sub-cases (i) $t_{33}\in\mb{Z}_2^*$ and (ii) $t_{33}\in 2\mb{Z}_2$. Note that $\mc{T}^{b,i}$ is just defined by $t_{22}|t_{12}$ and $t_{33}\in\mb{Z}_2^*$, so the integral over $\mc{T}^{b,i}$ becomes 
	\begin{align*}
		\frac{1}{(1-2^{-1})^3}\int_{\mc{T}^{b,i}}|t_{11}|_2^{s-1}|t_{22}|_2^{s-2}d\mu=\frac{1}{(1-2^{-1})^2}\int |t_{11}|_2^{s-1}|t_{22}|_2^{s-1}d\mu=\zeta_2(s)^2.
	\end{align*}
	In sub-case (ii), condition (2) implies that $t_{11}\in 2\mb{Z}_2$, whence $\mc{T}^{b,ii}$ is defined by the conditions $2|t_{11}$, $2|t_{33}$, $t_{22}|t_{12}$
	 and $\frac{t_{33}}{2}|-\frac{t_{12}}{t_{22}}t_{23}+\frac{t_{11}}{2}+t_{13}$. By using Proposition \ref{proceso de simplifiacion de integrales} with the pivot $t_{13}$, the integral over $\mc{T}^{b,ii}$ becomes
	\begin{align*}
	\frac{1}{(1-2^{-1})^3}\int_{\mc{T}^{b,ii}}|t_{11}|_2^{s-1}|t_{22}|_2^{s-2}|t_{33}|_2^{s-3}d\mu=\frac{1}{(1-2^{-1})^3}\int_{t_{11},t_{33}\in 2\mb{Z}_2}|t_{11}|_2^{s-1}|t_{22}|_2^{s-1}2|t_{33}|_2^{s-2}d\mu=2^{2-2s}\zeta_2(s)^2\zeta_2(s-1).
	\end{align*}

	We conclude that $
	\zeta_{G_2}^{N_2,\n}(s)=\zeta_2(s)^2(2^{-s}\zeta_2(s-1)+1+2^{2-2s}\zeta_2(s-1))= \zeta_2(s)^2\zeta_2(s-1)(1-2^{-s}+4\cdot 4^{-s})
	$.

\subsubsection{Calculation of $\zeta_{G}^{G,\n}(s)$} Note that $[N,G]= y_1^{-2}y_3$, and $G/[N,G]\cong \mb{Z}^2$. By Proposition \ref{simplification for normal extension zeta function}, we have $\zeta_{G}^{G,\n}(s)=\zeta_{G/[N,G]}^{G/[N,G],\n}(s)$ (with respect to $N/[N,G]$). Applying Lemma \ref{combinatorial lemma} with $Z_1=G/[N,G]$ and $Z_0=N/[N,G]$, we obtain that $\zeta_{G}^{G,\n}(s)=\zeta(s)\zeta(s-1)(1-2^{-s})$.

\medskip 
	
\subsection{Computing the zeta functions of $G=\ms{B}_3$.}
Note that $Z=\langle x_1\rangle$. One can easily check that  $E=G/Z$ is the plane crystallographic group $\mathbf{p2mg}$.

\subsubsection{Local factors of $\zeta_{G}^{G,\s}(s)$}\label{B_3, zeta_G^G, subgroup} 
If $p\neq 2$, then by (\ref{general remark}) we have $\zeta_{G_p}^{G_p,\s}(s)=\zeta_p(s)\zeta_{E_2}^{E_2,\s}(s)$, which is equal to $\zeta_p(s)\zeta_p(s-1)^2$ according to the results of \cite[5.6]{Mc}. 

We now compute $\zeta_{G_2}^{G_2,\s}(s)$. An open subgroup $A\lhd G_2$ such that $AN_2=G_2$ must contain  $\alpha \x^\v$ for some $\v\in\mb{Z}_2^3$, and hence it contains $(\alpha\x^\v)^2=x_1^{2v_1+1}$. Since $2v_1+1\in\mb{Z}_2^*$, we deduce that $Z_2\subset A$. Therefore, $\zeta_{G_2}^{G_2,\s}(s)=\zeta_{E_2}^{E_2,\s}(s)$, which is equal to  $(1+2^{1-s})\zeta_2(s-1)$ according to  \cite[5.6]{Mc}.

\subsubsection{Local factors of $\zeta_G^{N,\n}(s)$}\label{B_3, zeta_G^N, normal} If $p\neq 2$, then by (\ref{general remark}) we have $\zeta_{G_p}^{N_p,\n}(s)=\zeta_p(s)\zeta_{E_2}^{T_2,\n}(s)$, which is equal to $\zeta_p(s)^3$ according to \cite[6.6]{Mc}. 

Assume now that $p=2$. By Corollary \ref{integral expression for zeta_G^N},
	$$\zeta_{G_2}^{N_2,\n}(s)=\frac{1}{(1-2^{-1})^3}\int_{\mc{T}}|t_{11}|_2^{s-1}|t_{22}|_2^{s-2}|t_{33}|_2^{s-3}d\mu,$$
	where $\mc{T}\subset\on{T}_3(\mb{Z}_2)$ is defined by the following equivalent conditions:
	\begin{align*}
	&  ^\alpha(\x^{\t_1})\x^{-\t_1},\ ^\alpha(\x^{\t_2}),\ ^\alpha(\x^{\t_3}),\  ^\varepsilon(\x^{\t_1})\x^{-\t_1},\ ^\varepsilon(\x^{\t_2})\x^{-\t_2},\ ^\varepsilon(\x^{\t_3}) \in\overline{\langle\x^{\t_1},\x^{\t_2},\x^{\t_3}\rangle}\\
	\Leftrightarrow\quad 	& x_2^{2t_{12}}x_3^{2t_{13}},\ \x^{-\t_2},\ \x^{-\t_3},\  x_3^{2t_{13}}, x_3^{2t_{23}},\ \x^{-\t_3} \in \overline{\langle\x^{\t_1},\x^{\t_2},\x^{\t_3}\rangle}\\
	\Leftrightarrow\quad  & x_2^{2t_{12}}, x_3^{2t_{13}},\ x_3^{2t_{23}} \in \overline{\langle\x^{\t_1},\x^{\t_2},\x^{\t_3}\rangle}. 
	\end{align*}
This is the same integral that showed up in \ref{G_6, zeta_G^N, normal} when $p=2$, so the result is $(1+4\cdot 2^{-s}+4^{-s})\zeta_2(s)^3$.

\subsubsection{Local factors of $\zeta_{G}^{G,\n}(s)$}
 If $p\neq 2$, then $[N_p,G_p]=\overline{\langle x_2^2,x_3^2\rangle}=\overline{\langle x_2,x_3\rangle}$, whence $N_p/[N_p,G_p]\cong\mb{Z}_p$. Thus, by
 Proposition \ref{simplification for normal extension zeta function},  $\zeta_{G_p}^{G_p,\n}(s)=\zeta_p(s)$.
 
 Assume now that $p=2$. By the same argument used in the calculation of $\zeta_{G_2}^{G_2,\s}(s)$ in \ref{B_3, zeta_G^G, subgroup}, we have $\zeta_{G_2}^{G_2,\n}(s)=\zeta_{E_2}^{E_2,\n}(s)$, which is equal to $1+4\cdot 2^{-s}$ according to \cite[6.6]{Mc}.

\subsubsection{Local factors of $\zeta_{G}^{H,\n}(s)$ for  $H=\langle \alpha,x_1,x_2,x_3\rangle$} Note that $H=\ms{G}_2$. Arguing as in \ref{G_4, zeta_G^H, normal}, we find that $\zeta_{G_p}^{H_p,\n}(s)$ if $p\neq 2$, and that $\zeta_{G_2}^{H_2,\n}(s)=\zeta_{E_2}^{H_2/Z_2,\n}(s)$. The latter series is equal to $1+2^{1-s}$ according to \cite[6.6]{Mc}.

\subsubsection{Local factors of $\zeta_{G}^{K,\n}(s)$ for  $K=\langle \varepsilon,x_1,x_2,x_3\rangle$}
Note that $K=\ms{B}_2$. If $p\neq 2$, then by (\ref{general remark}) we have $\zeta_{G_p}^{K_p,\n}(s)=\zeta_{Z_p}^\s(s)\zeta_{E_p}^{K_p/Z_p,\n}(s)$, which is equal to $\zeta_p(s)\zeta_p(s)$ according to \cite[6.6]{Mc}. 

Assume now that $p=2$. We claim that if $A\lhd G_2$ is open and satisfies $AN_2=K_2$, then it must include $\overline{\langle x_2\rangle}$. Indeed, $A$ includes $[K_2,K_2]=\overline{\langle x_3^2\rangle}$ by Lemma \ref{Gamma_c(G) is included in A}. It also contains $\varepsilon \x^\v$ for some $\v\in\mb{Z}_2^3$, whence it contains $[\alpha,\varepsilon \x^\v]=x_2^{2v_2+1}x_3^{2v_3}$. It follows that $x_2^{2v_2+1}\in A$ and hence $x_2\in A$. Therefore, $\zeta_{G_2}^{K_2,\n}(s)=\zeta_{G_2/\overline{\langle x_2\rangle}}^{K_2/\overline{\langle x_2\rangle},\n}(s)$. We claim that this series is equal to $\zeta_{K_2/\overline{\langle x_2\rangle}}^{K_2/\overline{\langle x_2\rangle},\n}(s)$. Indeed, $\alpha$ and $\varepsilon$ commute in $G_2/\overline{\langle x_2\rangle}$ and their actions on $N_2/\overline{\langle x_2\rangle}$ are the same; thus, a normal subgroup of $K_2/\overline{\langle x_2\rangle}$ is already normal in $G_2/\overline{\langle x_2\rangle}$. Finally, observe that $K/\langle x_2\rangle$ is isomorphic to the plane crystallographic group $\textbf{pm}$. According to the results of \cite[6.3]{Mc}, we obtain that $\zeta_{G_2}^{K_2,\n}(s)=(1+5\cdot 2^{-s}+2\cdot 4^{-s})\zeta_2(s)$.

\subsubsection{Local factors of  $\zeta_{G}^{L,\n}(s)$ for $L=\langle \varepsilon\alpha ,x_1,x_2,x_3\rangle$:} Note that  $L=\ms{B}_2$.
If $p\neq 2$, then by (\ref{general remark}) we have $\zeta_{G_p}^{L_p,\n}(s)=\zeta_p(s)\zeta_{E_p}^{L_p/Z_p,\n}(s)$. This is equal to $\zeta_p(s)\zeta_p(s)$ according to the results of \cite[6.6]{Mc}.

Assume now that $p=2$. By Proposition \ref{integral expression for zeta_H^H and zeta_G^H} (with $H=L$), we have
	$$\ds \zeta_{G_2}^{L_2,\n}(s)=\frac{1}{(1-2^{-1})^3}\int_{\mc{T}}|t_{11}|_2^{s-2}|t_{22}|_2^{s-3}|t_{33}|_2^{s-4}d\mu,$$ 
	where $\mc{T}\subset\on{T}_3^+(\mb{Z}_2)\times\mb{Z}_2^3$ is the set of pairs $(\t,\v)$ satisfying the following equivalent conditions:
	\begin{align*}
&\begin{array}{c}  {}^{\alpha}(\x^{\t_1})\x^{-\t_1},\  {}^{\alpha}(\x^{\t_2}),\  {}^{\alpha}(\x^{\t_3}),\   {}^{\varepsilon}(\x^{\t_1})\x^{-\t_1},\  {}^{\varepsilon}(\x^{\t_2})\x^{-\t_2},\  {}^{\varepsilon}(\x^{\t_3})\\ (\varepsilon\alpha\x^{\v})^2,\  [\varepsilon\alpha,x_1],\  [\varepsilon\alpha,x_2],\  [\varepsilon\alpha,x_3],\  [\alpha,\varepsilon\alpha\x^\v] 
\end{array}
\in\overline{\langle\x^{\t_1},\x^{\t_2},\x^{\t_3}\rangle},\\
\Leftrightarrow\quad 	& x_2^{-2t_{12}}x_3^{-2t_{13}},\ \x^{-\t_2},\ \x^{-\t_3},\ x_3^{-2t_{13}},\ x_3^{-2t_{23}},\ \x^{-\t_3},\  x_1^{2v_1+1}x_3^{2v_3},\ x_2^2,\   x_2^{2v_2+1}x_3^{-2v_3}\in\overline{\langle\x^{\t_1},\x^{\t_2},\x^{\t_3}\rangle},\\
		\Leftrightarrow\quad &x_3^{2t_{13}},\  x_3^{2t_{23}},\ x_1^{2v_1+1}x_3^{2v_3},\  x_2^2,\ x_2x_3^{-2v_3}\in\overline{\langle\x^{\t_1},\x^{\t_2},\x^{\t_3}\rangle}.
	\end{align*}
By Lemma \ref{obtaining cone conditions abelian cas}, these conditions are equivalent to
	\begin{align*}
	&	(1)\	t_{33}|2t_{13},&&
		(2)\ t_{33}|2t_{23},&&
		(3)\ t_{11}|2v_1+1,&&
		(4)\ t_{22}|-\frac{2v_1+1}{t_{11}}t_{12},&&(5)\  t_{33}|\frac{2v_1+1}{t_{11}t_{22}}t_{12}t_{23}-\frac{2v_1+1}{t_{11}}t_{13}+2v_3,\\
	&	 (6)\ t_{22}|2,&& 
			(7)\ t_{33}|\frac{2}{t_{22}}t_{23},&&\ (8)\ t_{22}|1,&& (9)\ t_{33}|-\frac{1}{t_{22}}t_{23}-2v_3.&&
	\end{align*} 
Note that (3) and (8) are equivalent to $t_{11},t_{22}\in\mb{Z}_2^*$ and hence (4) and (6) are redundant. Notice also that (2) implies (7). To sum up, $\mc{T}$ is defined by (1), (2), (5), (9) and $t_{11},t_{22}\in\mb{Z}_2^*$.
We express $\mc{T}$ as a disjoint union $\mc{T}^a\cup\mc{T}^b$ according to the cases (2a) $t_{33}|t_{23}$ and (2b) $|t_{33}|_2=|2t_{23}|_2$.

		If we assume (2a), then (9) can be replaced by $t_{33}|2v_3$ and then (5) can be replaced by $t_{33}|t_{13}$, which makes (1) redundant. Thus, the conditions defining $\mc{T}^a$ are $t_{11},t_{22}\in\mb{Z}_2^*$, $t_{33}|t_{13}$, $t_{33}|t_{23}$ and $t_{33}|2v_3$, and the integral over $\mc{T}^a$ becomes
	$$
	\frac{1}{(1-2^{-1})^3}\int_{\mc{T}^a}|t_{33}|_2^{s-4}d\mu=\frac{1}{1-2^{-1}}\int_{ t_{33}|2v_3} |t_{33}|_2^{s-2}d\mu=1+2^{1-s}\zeta_2(s).
	$$

		If we assume (2b), then (9) implies $2|t_{23}$. Note also that (1) and (9) can be replaced by (1') $t_{23}|t_{13}$ and (9') $t_{23}|-\frac{t_{23}}{2t_{22}}-v_3$, and the latter implies that $t_{23}|2v_3$. This enables us to replace (5) by (5') $2|\frac{2v_1+1}{t_{11}t_{22}}t_{12}-\frac{2v_1+1}{t_{11}}\frac{t_{13}}{t_{23}}-\frac{2v_3}{t_{23}}$. The conditions defining $\mc{T}^b$ are, therefore, $t_{11},t_{22}\in\mb{Z}_2^*$, $2|t_{23}$, $|t_{33}|_2=|2t_{23}|_2$, (1'), (5'), and (9'), and the integral over $\mc{T}^b$ becomes
	$$
	\frac{1}{(1-2^{-1})^3}\int_{\mc{T}^b}|2t_{23}|_2^{s-4}d\mu=\frac{1}{1-2^{-1}}\int_{2|t_{23}}2^{4-s}|t_{23}|_2^{s-1}2^{-3}d\mu=2^{1-2s}\zeta_2(s),
	$$
	where in the first equality we first applied Proposition \ref{proceso de simplifiacion de integrales} with the pivots $(t_{11},t_{22},t_{13},t_{12},t_{33},v_3)$.

	Therefore, $\zeta_{G_2}^{L_2,\n}(s)=1+2^{1-s}\zeta_2(s)+2^{1-2s}\zeta_2(s)=\zeta_2(s)(1+2^{-s}+2\cdot 4^{-s})$. 

\medskip 

	\subsection{Computing the zeta functions of $G=\ms{B}_4$} Note that $Z=\langle x_1\rangle$. One can easily check that $E=G/Z$ is the plane crystallographic group $\textbf{p2gg}$.

	\subsubsection{Local factors of $\zeta_G^{G,\s}(s)$} 
	If $p\neq 2$, then by (\ref{general remark}) we have $\zeta_{G_p}^{G_p,\s}(s)=\zeta_p(s)\zeta_{E_p}^{E_p,\s}(s)$,
	and according to the results of \cite[5.7]{Mc} the latter is equal to $\zeta_p(s)\zeta_p(s-1)^2$. 
	
	Assume that $p=2$. An open subgroup $A\leq G_2$ such that $AN_2=G_2$ contains $\alpha\x^\v$ for some $\v\in\mb{Z}_2^3$, whence $(\alpha\x^\v)^2=x_1^{2v_1+1}\in A$. Since $2v_1+1\in\mb{Z}_2^*$, we conclude that $Z_2\subseteq A$. It follows that $\zeta_{G_2}^{G_2,\s}(s)=\zeta_{E_2}^{E_2,\s}(s)$, which is 1, according to \cite[5.7]{Mc}. 

\subsubsection{Calculus of  $\zeta_{G}^{N,\n}(s)$} Note that $N$, as $\mb{Z}[P]$-module, is isomorphic to the analogue module in the case $G=\ms{B}_3$. Therefore,
$\zeta_{G}^{N,\n}(s)=\zeta(s)^3(1+4\cdot 2^{-s}+4^{-s})$ according to the result of \ref{B_3, zeta_G^N, normal}.

\subsubsection{Local factors of $\zeta_G^{G,\n}(s)$}
 If $p\neq 2$, then $[N_p,G_p]=\overline{\langle x_2^2,x_3^2\rangle}=\overline{\langle x_2,x_3\rangle}$, whence $N_p/[N_p,G_p]\cong\mb{Z}_p$. Thus, by Proposition \ref{simplification for normal extension zeta function}, $\zeta_{G_p}^{G_p,\n}(s)=\zeta_p(s)$.
 
 Assume now that $p=2$. If $A\lhd G_2$ is open and satisfies $AN_2=G_2$, then $A$ contains $\alpha\x^\v$ and $\varepsilon\x^{\w}$ for some $\v,\w\in\mb{Z}_2^3$,
whence it contains both $(\alpha\x^\v)^2=x_1^{2v_1+1}$ and $(\beta\x^\w)^2=x_1^{2w_1}x_2^{2w_2+1}$. Since $2v_2+1, 2w_2+1\in\mb{Z}_2^*$, we conclude that $x_1,x_2\in A$. In addition, $[G_2,G_2]\subseteq A$ by Lemma \ref{Gamma_c(G) is included in A}; in particular $x_2x_3\in A$. It follows that $x_1,x_2,x_3\in A$ and hence $A=G_2$. Thus, $\zeta_{G_2}^{G_2,\n}(s)=1$.

\subsubsection{The local factors of $\zeta_{G}^{H,\n}(s)$ for $H=\langle \alpha,x_1,x_2,x_3\rangle$}
 Note that $H\cong \ms{G}_2$. Arguing as in \ref{G_4, zeta_G^H, normal}, we find that $\zeta_{G_p}^{H_p,\n}(s)=\zeta_p(s)$, and that $\zeta_{G_2}^{H_2,\n}(s)=\zeta_{E_2}^{H_2/Z_2,\n}(s)$. This series is equal to $1+2\cdot 2^{-s}$ according to the results of \cite[6.7]{Mc}.

\subsubsection{Local factors of $\zeta_{G}^{K,\n}(s)$ for $K=\langle\varepsilon,x_1,x_2,x_3\rangle$}
Note that $K=\ms{B}_2$. If $p\neq 2$, then by (\ref{general remark}) we have $\zeta_{G_p}^{K_p,\n}(s)=\zeta_p(s)\zeta_{E_p}^{K_p/Z_p,\n}(s)$, which is equal to $\zeta_p(s)\zeta_p(s)$ according to the results of \cite[6.7]{Mc}. 

Assume now that $p=2$. By Proposition \ref{integral expression for zeta_H^H and zeta_G^H} with $H=K$, we have
	$$\ds \zeta_{G_2}^{K_2,\n}(s)=\frac{1}{(1-2^{-1})^3}\int_\mc{T}|t_{11}|_2^{s-2}|t_{22}|_2^{s-3}|t_{33}|_2^{s-4}d\mu,$$ 
	where $\mc{T}\subset \on{T}_3^+(\mb{Z}_2)\times\mb{Z}_2^3$ is the set of pairs $(\t,\v)$ satisfying the following equivalent conditions:
	\begin{align*}
		&\begin{array}{c}  {}^{\alpha}(\x^{\t_1})\x^{-\t_1},\  {}^{\alpha}(\x^{\t_2}),\  {}^{\alpha}(\x^{\t_3}),\   {}^{\varepsilon}(\x^{\t_1})\x^{-\t_1},\  {}^{\varepsilon}(\x^{\t_2})\x^{-\t_2},\  {}^{\varepsilon}(\x^{\t_3})\\ (\varepsilon\x^{\v})^2,\  [\varepsilon,x_1],\  [\varepsilon,x_2],\  [\varepsilon,x_3],\  [\alpha,\varepsilon\x^\v] 
		\end{array}
		\in\overline{\langle\x^{\t_1},\x^{\t_2},\x^{\t_3}\rangle}\\
\Leftrightarrow	\quad	&\ x_2^{-2t_{12}}x_3^{-2t_{13}},\ \x^{-\t_2},\ \x^{-\t_3},  x_3^{-2t_{13}},\ x_3^{-2t_{23}}, \ \x^{-\t_3},\ x_1^{2v_1} x_2^{2v_2+1},\  x_3^2,\  x_2^{-1-2v_2}x_3^{-1+2v_3} \in\overline{\langle \x^{\t_1},\x^{\t_2},\x^{\t_3}\rangle}\\
\Leftrightarrow\quad	& x_2^{2t_{12}},\  x_1^{2v_1} x_3^{-1},\ x_3^2,\  x_2^{-1-2v_2}x_3^{-1} \in \overline{\langle \x^{\t_1},\x^{\t_2},\x^{\t_3}\rangle}.
	\end{align*}
By Lemma \ref{obtaining cone conditions abelian cas}, these conditions are equivalent to
	\begin{align*}
&		(1)\ t_{22}|2t_{12},&&
		(2)\ t_{33}|\frac{2t_{12}}{t_{22}}t_{23},&&
		(3)\ t_{11}|2v_1, &&
		(4)\ t_{22}|-\frac{2v_1}{t_{11}}t_{12},\\
		&		(5)\ t_{33}|\frac{2v_1}{t_{11}}(\frac{t_{12}}{t_{22}}t_{23}-t_{13})-1,&&
		(6)\ t_{33}|2,&&
		(7)\ t_{22}|-1-2v_2,&&
		(8)\ t_{33}|\frac{1+2v_2}{t_{22}}t_{23}-1.
	\end{align*}
		Note that (7) is equivalent to $t_{22}\in\mb{Z}_2^*$,  whence (1) and (4) are redundant and (6) implies  (2). 	We express $\mc{T}$ as a disjoint union $\mc{T}^a\cup\mc{T}^b$ according to the cases  (3a) $t_{11}|v_1$ and (3b) $|t_{11}|_2=|2 v_1|_2$.

	In the case (3a), (6) reduces (5) to the condition $t_{33}|-1$, that is $t_{33}\in\mb{Z}_2^*$; hence, $\mc{T}^a$ is defined by $t_{11}|v_1$, $t_{22},t_{33}\in\mb{Z}_2^*$. The integral over $\mc{T}^a$ becomes
	$$
	\frac{1}{(1-2^{-1})^3}\int_{\mc{T}^a}|t_{11}|_2^{s-2}d\mu=\frac{1}{1-2^{-1}}\int |t_{11}|_2^{s-1}d\mu=\zeta_2(s).
	$$
	
	We now assume (3b).	Note that $\mc{T}^b$ is defined by $|t_{11}|_2=|2v_1|_2$, $t_{22}\in\mb{Z}_2^*$,  (5), (6), and (8), and the integral over $\mc{T}^b$ becomes 
	\begin{align*}
		\frac{1}{(1-2^{-1})^3}\int_{\mc{T}^b}|2v_1|_2^{s-2}|t_{33}|_2^{s-4}d\mu=\frac{1}{(1-2^{-1})^2}\int_{ t_{33}|2}2^{2-s}|v|_2^{s-2} 2^{-2}|t_{33}|_2^{s-2}d\mu=(1+ 2\cdot 2^{-s})2^{-s}\zeta_2(s),
	\end{align*}
	where in the second equality we applied Proposition \ref{proceso de simplifiacion de integrales} with the pivots $(t_{11},t_{22},t_{13},t_{23})$.

		We conclude that $\zeta_{G_2}^{K_2,\n}(s)=\zeta_2(s)+(1+ 2^{-(s-1)})2^{-s}\zeta_2(s)=(1+2^{-s}+2\cdot 4^{-s})\zeta_2(s)$.

\subsubsection{Local factors of $\zeta_{G}^{L,\n}(s)$ for  $L=\langle \varepsilon\alpha,x_1,x_2,x_3\rangle$} 
Note that $L\cong \ms{B}_2$. If $p\neq 2$, then by (\ref{general remark}) we have $\zeta_{G_p}^{L_p,\n}(s)=\zeta_p(s)\zeta_{E_p}^{L_p/Z_p,\n}(s)$, which is equal to $\zeta_p(s)\zeta_p(s)$ according to the results of \cite[6.7]{Mc}. 

Assume now that $p=2$.
We claim that any $A\lhd G_2$ that is open and satisfies $AN_2=L_2$ must include $\overline{\langle x_1x_2,x_2x_3,x_2^2\rangle}$.
Indeed, $A$ includes $[G_2,N_2]=\overline{\langle x_2^2 \rangle}$ by Lemma \ref{Gamma_c(G) is included in A}. It also contains $\varepsilon \alpha\x^\v$ for some $\v\in\mb{Z}_2^3$, whence it contains $(\varepsilon\alpha\x^\v)^2=x_1^{2v_1+1}x_3^{2v_3+1}$ and $[\alpha,\varepsilon\alpha\x^\v]=x_2^{2v_2-1}x_3^{-2v_3-1}$. We deduce that $x_1^{2v_1+1}x_2^{2v_2-1}\in A$, and using that $x_2^2\in A$ we obtain that $(x_1x_2)^{2v_1+1}\in A$. It follows that $x_1x_2\in A$ since $2v_1+1$ is a unit. With a similar argument, we find that $x_2x_3\in A$. This proves the claim. Since $[\varepsilon,\alpha]=x_2x_3$, the quotient group $G':=G_2/\overline{\langle x_1x_2,x_2x_3,x_2^2\rangle}$ is abelian, and therefore, $\zeta_G^{L,\n}(s)=\zeta_{L'}^{L',\s}(s)$, where $L'$ is the image of $L_2$ at $G'$, and the partial zeta function on the right is computed with respect to $N'$, the image of $N_2$ at $G'$. Note that $H'\cong C_2\times C_2$ and $N'$ has index 2. Thus, $\zeta_{L'}^{L',\s}(s)=1+2\cdot 2^{-s}$.

\section{Computing the zeta functions of the 3-dimensional AB-groups}\label{Section: proof for AB-groups}
\noindent In this final section, we prove the formulae presented in \ref{3-dimensional AB groups and their zeta functions}.
We always denote by $G$ the AB-group under consideration and by $N$ the Fitting subgroup (which is generated by $x_1,x_2,x_3$). We keep the notation introduced at the beginning of Section \ref{Section: Local factors at good primes}. In particular,
$$Z=Z(N), \quad T=N/Z,\quad E=G/Z,\quad P=G/N.
$$
The formulae for the zeta functions of $G=N=\langle x_1,x_2,x_3: [x_2,x_1]=x_3^k,\ [x_1,x_3]=[x_2,x_3]=1\rangle$, where $k>0$, are already computed; see Table \ref{Table with zeta_G^G} and \ref{Table with zeta_G^N} with $F=C_1$.
We focus here on the cases $G\neq N$.  If $N\leq H\leq G$ is an intermediate subgroup, then $H$ is also a 3-dimensional AB-group with Fitting subgroup $N$. 
Therefore, in order to compute $\zeta_G^{G,\s}(s)$, it will be enough to calculate $\zeta_G^{G,\s}(s)$ since the other partial zeta functions will have been computed in previous calculations. The isomorphism classes of the intermediate subgroups $N\subsetneqq H\subsetneqq G$ (if there is any) will be identified when computing the partial zeta functions $\zeta_G^{H,\n}(s)$.

When dealing with local factors at primes $p|[G:N]$, we will sometimes use the method of $p$-adic integration described in Section \ref{Section: method of p-adic integration}; specifically  Corollary \ref{integral expression for zeta_G^N} and Proposition \ref{integral expression for zeta_H^H and zeta_G^H}. The following lemma collects some information about the $\mf{T}$-group $N$.

\begin{lemma}\label{obtaining cone conditions} 
Let $N=\langle x_1,x_2,x_3: [x_2,x_1]=x_3^k,\ [x_1,x_3]=[x_2,x_3]=1\rangle$, with $k=[Z:[N,N]]>0$. 
\begin{enumerate}
	\item For $\a,\b\in\mb{Z}_p^3$ and $r\in\mb{Z}_p$, it holds that 
	\begin{align*}
		\x^\a\cdot \x^\b=x_1^{a_1+b_1}x_2^{a_2+b_2}x_3^{a_3+b_3+kb_1a_2}, & & (\x^\a)^r=x_1^{ra_1}x_2^{ra_2}x_3^{ra_3+\binom{r}{2}ka_1a_2}, & & [\x^\a,\x^\b] &=x_3^{k(b_1a_2-a_1b_2)}.	
	\end{align*}
\item 	For $\t\in\on{T}_3^+(\mb{Z}_p)$, the following holds. 
\begin{enumerate}
	\item $\t\in\mc{M}_{N_p}^\s$ if and only if $t_{33}|kt_{11}t_{22}$. 
	\item  $\t\in\mc{M}_{N_p}^\n$  if and only if $t_{33}|kt_{11}$, $t_{33}|kt_{12}$ and $t_{33}|kt_{22}$. 
\end{enumerate}
\item 	Fix  $\t\in\mc{M}_{N_p}^\s$. Given $\a=(a_1,a_2,a_3)\in\mb{Z}_p^3$, it holds that $\x^{\a}\in \overline{\langle\x^{\t_1},\x^{\t_2},\x^{\t_3}\rangle}$ if and only if
\begin{enumerate}
	\item $t_{11}|a_1$.
	\item $t_{22}|-\frac{a_1}{t_{11}}t_{12}+a_2$.
	\item $t_{33}|-\frac{-\frac{a_1}{t_{11}}t_{12}+a_2}{t_{22}}t_{23}-\frac{a_1}{t_{11}}t_{13}-\frac{1}{2}k\frac{a_1}{t_{11}}(\frac{a_1}{t_{11}}-1)t_{11}t_{12}+a_3$.
\end{enumerate}
In the case that $\t\in\mc{N}_{N_p}^\n$, (c) can be replaced by: $t_{33}|-\frac{-\frac{a_1}{t_{11}}t_{12}+a_2}{t_{22}}t_{23}-\frac{a_1}{t_{11}}t_{13}+a_3$.
\end{enumerate}
\end{lemma} 
\begin{proof}
	The verification of (1) is straightforward, and (2) is a special case of (\ref{condition to represent a good basis}). As for (3),
	given that $\{\x^{\t_1},\x^{\t_2},\x^{\t_3}\}$ is a good basis for $\overline{\langle \x^{\t_1},\x^{\t_2},\x^{\t_3}\rangle}$, the relation $\x^{\a}\in \overline{\langle \x^{\t_1},\x^{\t_2},\x^{\t_3}\rangle}$ holds if and only  $\x^{\a}=(\x^{\t_1})^{\lambda_1}(\x^{\t_2})^{\lambda_2}(\x^{\t_3})^{\lambda_3}$ for some $\lambda_1,\lambda_2,\lambda_3\in\mb{Z}_p$. Using (1), we find that
	\begin{align*}
		(\x^{\t_1})^{\lambda_1}(\x^{\t_2})^{\lambda_2}(\x^{\t_3})^{\lambda_3}=x_1^{t_{11}\lambda_1}x_2^{t_{12}\lambda_1+t_{22}\lambda_2}x_3^{t_{13}\lambda_1+t_{23}\lambda_2+k\binom{\lambda_1}{2}t_{11}t_{12}+t_{33}\lambda_3}.
	\end{align*}
	The first part of (3) is now clear. In the case that $\t\in\mc{M}_{N_p}^\n$ (and hence $\overline{\langle \x^{\t_1},\x^{\t_2},\x^{\t_3}\rangle}$ is normal), (2) implies that $t_{33}|kt_{11}$. Therefore, the term $\frac{1}{2}k\frac{a_1}{t_{11}}(\frac{a_1}{t_{11}}-1)t_{11}t_{12}=\binom{\frac{a_1}{t_{11}}}{2}kt_{11}t_{12}$ can be deleted from (c). This completes the proof.
\end{proof}

Finally, Table \ref{List of integrals} below collects the results of some $p$-adic integrals that will appear in our calculations. The verification of these formulae is straightforward. The notation is as follows: $k$ denotes an integer; the letters $s$, $t$, and $u$ are complex variables; the variables of integration are  $p$-adic integers; and $\theta_p$ denotes
$\frac{1}{1-p^{-1}}$.
\begin{table}[h]
	\resizebox{\columnwidth}{!}{\begin{tabular}{ll}
			\hline\noalign{\smallskip}
			$p$-adic integral & result \\
			\noalign{\smallskip}\hline\noalign{\smallskip}
			$A_{p,k}(s,t):=\theta_p^2 \int_{x|ky}|x|_p^{s-1}|y|_p^{t-1}d\mu$ &$\zeta_p(t)\left(\frac{1-|k|_p^s}{1-p^{-s}}+\frac{|k|_p^s}{1-p^{-(s+t)}}\right) $ \\
			$B_{p,k}(s,t):=\theta_p^2\int_{p|x|ky}|x|_p^{s-1}|y|_p^{t-1}d\mu$ &$\zeta_p(t)\left(\frac{p^{-s}(1-|k|_p^s)}{1-p^{-s}}+\frac{p^{-(s+t)}|k|_p^s}{1-p^{-(s+t)}}\right) $ \\
			$C_{p,k}(s,t,u):=\theta_p^3\int_{\substack{z|kx\\ z|ky}}|x|_p^{s-1}|y|_p^{t-1}|z|_p^{u-1}d\mu$ & $\zeta_p(s)\zeta_p(t)\left(\frac{1-|k|_p^u}{1-p^{-u}}+\frac{|k|_p^u}{1-p^{-(s+t+u)}}\right) $
			\\
			$D_{p,k}(s,t,u):=\theta_p^3\int_{\substack{z|kx\\ z|pky}}|x|_p^{s-1}|y|_p^{t-1}|z|_p^{u-1}d\mu$& $\zeta_p(s)\zeta_p(t)\left(\frac{1-|k|_p^u}{1-p^{-u}}+\frac{(1+p^{-(u+s)}-p^{-(u+s+t)})|k|_p^u}{1-p^{-(s+t+u)}}\right) $
			\\
			$E_{p,k}(s,t,u):=\theta_p^3\int_{\substack{z|kx\\ z|pky\\ z|pw}}|x|_p^{s-1}|y|_p^{t-1}|z|_p^{u-1}d\mu$ & $\zeta_p(s)\zeta_p(t)\left(\frac{(1+p^{-u}-p^{-(u+1)})(1-|k|_p^{u+1})}{1-p^{-(u+1)}}+\frac{(1+p^{-(u+s)}-p^{-(s+t+u+1)})|k|_p^{u+1}}{1-p^{-(s+t+u+1)}} \right)$
			\\
			\noalign{\smallskip}\hline
	\end{tabular}}
	\caption{List of $p$-adic integrals}
\label{List of integrals}  
\end{table}

		\subsection{Computing the zeta functions of $G=G_{\mathbf{p2},2k}$, $k\in\mb{N}$}\label{zeta functions of Q=p2}\hfill
	
\subsubsection{Local factors of $\zeta_{G}^{G,\s}(s)$}\label{G_p2,2k, zeta_G^G,subgroup}
If $p\neq 2$, then $\zeta_{G_p}^{G_p,\s}(s)$ was given in Table \ref{Table with zeta_G^G} (the case $F=C_2$ and $F\subset SL(T)$).

	Assume now that $p=2$. If $A\leq G_2$ is open and satisfies $AN_2=G_2$, then $A$ contains  $\alpha \x^\v$ for some $\v\in\mb{Z}_2$, whence $(\alpha \x^\v)^2=x_3^{2v_3+1}\in A$. Since $2v_3+1\in \mb{Z}_2^*$, we obtain that $Z_p\subset  A$. This implies that $\zeta_{G_2}^{G_2,\s}(s) =\zeta_{E_2}^{E_2,\s}(s)$, which in turn is equal to $\zeta_2(s-1)\zeta_2(s-2)$, according to the results of \cite[5.2]{Mc}.

\subsubsection{Local factors of $\zeta_{G}^{N,\n}(s)$} If $p\neq 2$, then $\zeta_{G_p}^{N_p,\s}(s)$ was given in Table \ref{Table with zeta_G^N} (the case $F=C_2$ and $F\subset SL(T)$). 

Assume now that $p=2$.
We show that the assumptions of Proposition \ref{formula for the partial normal zeta function} are satisfied. Firstly, since the class of $\alpha$ at $P=G_2/N_2$ acts on $T_2$ as $-\on{id}_{T_2}$, any finite index subgroup of $T_2$ is a $\mb{Z}_2[P]$-submodule, and they are all isomorphic. Fix $U\leq T_2$  and $V\leq Z_2$ of finite index such that $[U,T_2]\subseteq V$. Let $\tilde{U}$ be the pre-image of $U$ in $N_2$. Note that $\tilde{U}$ has a good basis of the form $\{ x_1^{t_{11}}x_2^{t_{12}},x_2^{t_{22}},x_3\}$, and that $V=\overline{\langle x_3^{t_{33}}\rangle}$ for some non-zero $t_{33}$. It is easy to check that the condition $[U,T_2]\subseteq V$ is translated into $t_{33}|2kt_{11}$, $t_{33}|2kt_{12}$, $t_{33}|2kt_{22}$. Therefore, by Lemma \ref{obtaining cone conditions}(2), $\{x_1^{t_{11}}x_2^{t_{12}}, x_2^{t_{22}}, x_3^{t_{33}}\}$ represents a good basis for a normal subgroup $B\lhd N_2$. Clearly,  $(BZ_2)/Z=U$ and $B\cap Z_2=V$. We now show that $B$ is normal in $G_2$. Indeed, $^\alpha(x_1^{t_{11}}x_2^{t_{12}})=x_1^{-t_{11}}x_2^{-t_{12}}=(x_1^{t_{11}}x_2^{t_{12}})^{-1}x_3^{-2kt_{11}t_{12}}\in B$ and $^\alpha(x_2^{t_{22}})=x_2^{-t_{22}}\in B$.

We are now in position to apply Proposition \ref{formula for the partial normal zeta function}. For the computation of $[T_2:X(V)]$ we use Lemma \ref{index of the centralizer of V}, which tells us that $[T_2:X(V)]=1$ if $[Z_2:V]<2^{v_2(2k)}$ and $[T_2:X(V)]=|2k|_2^2[Z_2:V]^2$ if $V\leq \overline{\langle x_3^{2k}\rangle}$. For the computation of $|\on{Hom}_{\mb{Z}_2[P]}(T_2,Z_2/V)|$ we use Lemma \ref{Hom from C_2}, which tells us that $|\on{Hom}_{\mb{Z}_2[P]}(T_2,Z_2/V)|=1$ if $V=Z_2$ and is $2^2$ otherwise. Thus, by Proposition \ref{formula for the partial normal zeta function} and the results of \cite[6.2]{Mc},
\begin{align*}
	\zeta_{G_2}^{N_2,\n}(s)&=\zeta_{E_2}^{T_2,\n}(s)\left(1+\sum_{i=1}^{v_2(2k)-1} 2^{-is} 2^2+\sum_{V\leq 2k\mb{Z}_2}|2k|_2^{-2s}[\mb{Z}_2:V]^{-2s}[\mb{Z}_2:V]^{-s} 2^2 \right)
	\\
	&=\zeta_2(s)\zeta_2(s-1)\left( 1+2^{2-s}\frac{1-|k|^s}{1-2^{-s}}+|2k|^s 2^2\zeta_2(2s)\right). 
\end{align*}

\subsubsection{Local factors of $\zeta_{G}^{G,\n}(s)$}\label{G_p2,2k, zeta_G^G, normal}
 	If $p\neq 2$, then $\gamma_3(G_p,N_p)=\overline{\langle x_1^{4}, x_2^{4}, x_3^{4k}\rangle}=\overline{\langle x_1,x_2,x_3^k\rangle}$, whence $N_p/\gamma_{3}(G_p,N_p)\cong \mb{Z}_p/k\mb{Z}_p$. Thus, by Proposition \ref{simplification for normal extension zeta function}, $\zeta_{G_p}^{G_p}(s)=\frac{1-p^{-s}|k|_p^s}{1-p^{-s}}$. 

Assume now that $p=2$. As shown in the calculation of $\zeta_{G_2}^{G_2,\s}(s)$, any $A\lhd G_2$ that is open and satisfies $AN_2=G_2$ must include $Z_2$. Thus, $\zeta_{G_2}^{G_2,\n}(s)=\zeta_{E_2}^{E_2,\n}(s)$, which is equal to $1+6\cdot 2^{-s}+4\cdot 4^{-s}$ according to the results of \cite[6.2]{Mc}.

	\subsection{Computing the zeta functions of $G=G_{\mathbf{pg},2k}$}\label{zeta functions of Q=pg}\hfill 
	 
	\subsubsection{Local factors of $\zeta_{G}^{G,\s}(s)$}
	 If $p\neq 2$, then $\zeta_{G_p}^{G_p,\s}(s)$ was given in Table \ref{Table with zeta_G^G} (the case $F=D_1$ and $F\not\subset SL(T)$). 
	 
	 Assume now  that $p=2$. By Proposition \ref{integral expression for zeta_H^H and zeta_G^H} and Lemma \ref{obtaining cone conditions}(2), we have
	 $$\zeta_{G_2}^{G_2,\s}(s)=\frac{1}{(1-2^{-1})^3}\int_\mc{T}|t_{11}|_2^{s-2}|t_{22}|_2^{s-3}|t_{33}|_2^{s-4}d\mu,$$
	  where $\mc{T}\subset \on{T}_3^+(\mb{Z}_2)\times\mb{Z}_2^3$ is the set of pairs $(\t,\v)$ satisfying $t_{33}|2kt_{11}t_{22}$ and the following equivalent conditions: 
		\begin{align*}
		 &\x^{\t_1}\cdot{}^{\beta\x^\v}(\x^{\t_1}),\ [\beta\x^\v,\x^{\t_2}],\ {}^{\beta\x^\v}(\x^{\t_3}),\ (\beta\x^\v)^2\in\overline{\langle\x^{\t_1},\x^{\t_2},\x^{\t_3}\rangle}\\
		 \Leftrightarrow\quad & x_2^{2t_{12}}x_3^{-kt_{11}(1+2v_2)+2kv_1t_{12}-2kt_{11}t_{12}},\ x_3^{-2(t_{23}+kt_{22}v_1)},\ x_2^{2v_2+1}x_3^{kv_1(2v_2+1)} \in\overline{\langle\x^{\t_1},\x^{\t_2},\x^{\t_3}\rangle}.
		\end{align*}
	 By Lemma \ref{obtaining cone conditions}(3), the conditions defining $\mc{T}$ are equivalent to
		$$
		\begin{array}{lllll}
			(1)\ t_{33}|2kt_{11}t_{22},&\ \  & (2)\ t_{22}|2t_{12},& \ \  & (3)\ t_{33}|\frac{2t_{12}}{t_{22}}(-t_{23}+kv_1t_{22})-kt_{11}(1+2v_2)-2kt_{11}t_{12},\\
			(4)\ t_{33}|2(-t_{23}+kt_{22}v_1), & &(5)\ t_{22}|2v_2+1, & &(6)\ t_{33}|\frac{2v_1+1}{t_{22}}(-t_{23}+kv_1t_{22}).
		\end{array}
		$$
	Note that (5) is equivalent to $t_{22}\in\mb{Z}_2^*$, so (2) is redundant. Next, (6) implies (4), and (6) and (1) reduce (3) to the simpler condition $t_{33}|kt_{11}$. This makes (1) redundant. To sum up, the conditions defining $\mc{T}$ are: $t_{22}\in\mb{Z}_2^*$, $t_{33}|kt_{11}$, $t_{33}|-t_{23}+kv_1t_{22}$. Thus,
		\begin{align*}
			\zeta_{G_2}^{G_2,\s}(s)=\frac{1}{(1-2^{-1})^3}\int_\mc{T}|t_{11}|_2^{s-2}|t_{33}|_2^{s-4}d\mu=\frac{1}{(1-2^{-1})^2}\int_{  t_{33}|kt_{11}}|t_{11}|_2^{s-2}|t_{33}|_2^{s-3}d\mu=A_{2,k}(s-2,s-1),
		\end{align*}
	where we used Proposition \ref{proceso de simplifiacion de integrales} with the pivots $t_{22},t_{23}$.

\subsubsection{Local factors of $\zeta_{G}^{N,\n}(s)$}
 If $p\neq 2$, then $\zeta_{G_p}^{N_p,\n}(s)$ was given in Table \ref{Table with zeta_G^N} (the case $F=D_1$ and $F\not\subset SL(T)$). 
 
 Assume now that $p=2$. By Corollary \ref{integral expression for zeta_G^N} and Lemma \ref{obtaining cone conditions}(2),
		$$\zeta_{G_2}^{N_2,\n}(s)=\frac{1}{(1-2^{ -1})^3}\int_\mc{T}|t_{11}|_2^{s-1}|t_{22}|_2^{s-2}|t_{33}|_2^{s-3}d\mu,$$ where $\mc{T}\subset\on{T}_3^+(\mb{Z}_2)$ is defined by $t_{33}|2kt_{11}, t_{33}|2kt_{12}, t_{33}|2kt_{22}$ and the following equivalent conditions:
		\begin{align*}
		 &\x^{\t_1}\cdot{}^\beta(\x^{\t_1}),\ [\beta,\x^{\t_2}],\ \x^{\t_3}\cdot{}^\beta(\x^{\t_3})\in\overline{\langle \x^{\t_1},\x^{\t_2},\x^{\t_3}\rangle}\\
		 \Leftrightarrow\quad &  x_2^{2t_{12}}x_3^{-kt_{11}-2kt_{11}t_{12}}, x_3^{2t_{23}}\in\overline{\langle \x^{\t_1},\x^{\t_2},\x^{\t_3}\rangle}.
		\end{align*}
	Therefore, by Lemma \ref{obtaining cone conditions}(3), the conditions defining $\mc{T}$ are equivalent to
		\begin{align*}
			(1)\ t_{33}|2kt_{11},&& (2)\ t_{33}|2kt_{12},&& (3)\  t_{33}|2kt_{22},&& (4)\  t_{22}|2t_{12}&& (5)\ t_{33}|\frac{2t_{12}}{t_{22}}t_{23}+kt_{11},&&(6)\  t_{33}|2t_{23}.
		\end{align*}
	We express $\mc{T}$ as a disjoint union $\mc{T}^{a,a}\cup\mc{T}^{a,b}\cup\mc{T}^b$ according to the following three cases: (1a, 4a) $t_{33}|kt_{11}$ and $t_{22}|t_{12}$; (1a, 4b) $t_{33}|kt_{11}$ and $|t_{22}|_2=|2t_{12}|_2$; and (1b)  $|t_{33}|_2=|2kt_{11}|_2$.

	In the first case, conditions (2) and (5) are redundant, and hence the conditions defining $\mc{T}^{a,a}$ are $t_{33}|kt_{11}$, $t_{33}|2kt_{22}$, $t_{22}|t_{12}$ and $t_{33}|2t_{23}$, and the integral over $\mc{T}^{a,a}$ becomes 
	\begin{align*}
		\frac{1}{(1-2^{-1})^3}\int_{\mc{T}^{a,a}}|t_{11}|_2^{s-1}|t_{22}|_2^{s-2}|t_{33}|_2^{s-3}d\mu=\frac{1}{(1-2^{-1})^3}\int_{\substack{t_{33}|kt_{11}\\ t_{33}|2kt_{22}\\ t_{33}|2t_{23}}}|t_{11}|_2^{s-1}|t_{22}|_2^{s-1}|t_{33}|_2^{s-3}d\mu=E_{2,k}(s,s,s-2).
	\end{align*}

In the case  (1a, 4b), condition (5) can be replaced by $t_{33}|t_{23}$, and hence the conditions defining $\mc{T}^{a,b}$ are $t_{33}|kt_{11}$, $t_{33}|2kt_{12}$, $|t_{22}|_2=|2t_{12}|_2$, $t_{33}|t_{23}$. The integral over $\mc{T}^{a,b}$ becomes 
\begin{align*}
	\frac{1}{(1-2^{-1})^3}\int_{\mc{T}^{a,b}}|t_{11}|_2^{s-1}|2t_{12}|_2^{s-2}|t_{33}|_2^{s-3}d\mu=\frac{2^{-s}}{(1-2^{-1})^3}\int_{\substack{t_{33}|kt_{11}\\ t_{33}|2kt_{12}}}|t_{11}|_2^{s-1}|t_{12}|_2^{s-1}|t_{33}|_2^{s-2}d\mu=2^{-s}D_{2,k}(s,s,s-1),
\end{align*}
where in the second equality we used Proposition \ref{proceso de simplifiacion de integrales} with the pivots $t_{23},t_{22}$.

 We now assume (1b). Since in this case $t_{33}\nmid kt_{11}$, conditions (5) and (6) imply that $|t_{22}|_2=|2t_{12}|_2$. This, (1b) and (5) imply that $|t_{23}|_2=|kt_{11}|_2$. It is now easy to see that $\mc{T}^b$ is defined by the conditions  $|t_{33}|_2=|2kt_{11}|_2$, $|t_{22}|_2=|2t_{12}|$, $t_{33}|2kt_{12}$, $|t_{23}|_2=|kt_{11}|_2$, or equivalently, $|t_{33}|_2=|2kt_{11}|_2$, $|t_{22}|_2=|2t_{12}|$, $t_{11}|t_{12}$, $|t_{23}|_2=|kt_{11}|_2$. Thus, the integral over $\mc{T}^b$ becomes
 \begin{align*}
 	\frac{1}{(1-2^{-1})^3}\int_{\mc{T}^b}|t_{11}|_2^{s-1}|2t_{12}|_2^{s-2}|2kt_{11}|_2^{s-3}d\mu=\frac{2^{1-2s}|k|_2^{s-1}}{(1-2^{-1})^2}\int_{t_{11}|t_{12}} |t_{11}|_2^{s}|t_{12}|_2^{s-1}|t_{11}|_2^{s-1}d\mu=2^{1-2s}|k|_2^{s-1}A_{2,1}(2s-1,s), 	
 \end{align*}
where in the first equality we applied Proposition \ref{proceso de simplifiacion de integrales} with the pivots $(t_{22},t_{23},t_{33})$. 

It follows from Table \ref{List of integrals} that
$$
\zeta_{G_2}^{N_2,\n}(s)=\zeta_2(s)^2\left(
\zeta_2(s-1)(1+3\cdot 2^{-s})(1-|k|_2^{s-1})+\zeta_2(3s-1)(1+2^{-s}+3\cdot 2^{1-2s}-2^{1-3s}-2^{1-4s})|k|_2^{s-1} \right). 
$$

\subsubsection{Local factors of $\zeta_{G}^{G,\n}(s)$}\label{G_p2gg, zeta_G^G, normal} If $p\neq 2$, then  $\gamma_{3}(G_p,N_p)=\langle x_1^{4}, x_3^{4}\rangle=\overline{\langle x_1,x_3\rangle}$, whence $N_p/\gamma_3(G_p,N_p)\cong\mb{Z}_p$. Thus, by Proposition \ref{simplification for normal extension zeta function}, 
$\zeta_{G_p}^{G_p,\n}(s)=\zeta_p(s)$. 

For the computation of $\zeta_{G_2}^{G_2,\n}(s)$, we consider the cases $k$ even or odd separately.
		
		Assume that $k$ is even. We claim that if $A\lhd G_2$ is open and satisfies $AN_2=G_2$, then $A$  includes $\overline{\langle x_1^2,x_2,x_3^2\rangle}$. In fact, any such $A$ contains $\beta\x^\v$ for some $\v\in\mb{Z}_2^3$, whence $A$ contains $[\beta\x^\v,x_3]=[\beta,x_3]=x_3^{-2}$, $[\beta\x^\v,x_1]={}^\beta[\x^\v,x_1][\beta,x_1]=x_1^{-2} x_3^{-2kv_2}x_3^{-k}$ and $(\beta\x^\v)^2=x_2^{2v_2+1}x_3^{(2v_2+1)v_lk}$. From the first two  contentions, we deduce that $x_1^{-2}\in A$ since $k$ is even. From the first and third contentions and the fact that $k$ is even, we also deduce that $x_2^{2v_2+1}\in A$, and hence $x_2\in A$ since $2v_2+1\in\mb{Z}_2^*$. This proves the claim. It follows that
		$
		\zeta_{G_2}^{G_2,\n}(s)=\zeta_{G'}^{G',\n}(s)$, where $G'$ is the quotient of $G_2$ by $\overline{\langle x_1^2,x_2,x_3^2\rangle}$, and $\zeta_{G'}^{G',\n}(s)$ is computed with respect to $N'$, the image of $N_2$ at $G'$. Note that  $G':=(\mb{Z}/2\mb{Z})^3$, the classes of $\beta$, $x_1$ and $x_3$ being identified with $(1,0,0)$, $(0,1,0$ and $(0,0,1)$ respectively. The series $\zeta_{G'}^{G'}(s)$ enumerates the subgroups of $G'$ that are not included in $N'$. An easy calculation gives $\zeta_{G'}^{G',\n}(s)=1+6\cdot 2^{-s}+4\cdot 4^{-s}$.
		
		Assume now that  $k$ is odd. 
		With a similar analysis to the one in the previous case, we find that if $A\lhd G_2$ is open and satisfies $AN_2=G_2$, then $A$ includes $\overline{\langle x_1^2x_3, x_2^2, x_3^2\rangle}$.
	It follows that	$
		\zeta_G^{G,\n}=\zeta_{G'}^{G',\n}(s)$, where $G'$ is the quotient of $G$ by $\overline{\langle x_1^2x_3,x_2^2,x_3^2\rangle}$, and $\zeta_{G'}^{G',\n}(s)$ is computed with respect to $N'$, the image of $N$ at $G'$.
		Note that $G'\cong \mb{Z}/4\mb{Z}\times\mb{Z}/4\mb{Z}$, where the classes of $\beta$ and $x_1$ are identified with $(1,0)$ and $(0,1)$, respectively. Thus, the series $\zeta_{G'}^{G',\n}(s)$ enumerates the subgroups of $\mb{Z}/4\mb{Z}\times\mb{Z}/4\mb{Z}$ that are not included in  $2\mb{Z}/4\mb{Z}\times\mb{Z}/4\mb{Z}$. A simple inspection shows that the only such subgroups are $G'$, $\langle (1,0), (0,2)\rangle$, $\langle (1,1),(0,2)\rangle$, $\langle (1,0)\rangle$, $\langle (1,2)\rangle$, $\langle(1,3)\rangle$ and $\langle (1,4)\rangle$; thus, $\zeta_{G'}^{G',\n}(s)=1+2\cdot 2^{-s}+4\cdot 4^{-s}$.

	\subsection{Computing the zeta functions of  $G_{\textbf{p3},k,\epsilon}, G_{\textbf{p4},2k,\epsilon}, G_{\textbf{p6},2k,\epsilon}$}\label{G(k,p,r,d,e)} The general form of these groups is
		$$G_{q,\delta k,\epsilon}=\langle \gamma, x_1,x_2,x_3:\  [x_2,x_1]=x_3^{\delta k},\ [x_1,x_3]=[x_2,x_3]=1,\ \gamma^{q}=x_3^\epsilon,\ ^\gamma x_1=x_2,\ ^\gamma x_2=x_1^{-1}x_2^{-r} \rangle $$
	where $k\in\mb{N}$, $(d,\delta,r)\in \{(3,1,1), (4,2,0), (6,2,-1)\}$ and $\epsilon=\left\{\begin{array}{ll} 1\ \mbox{or}\ -1 &\mbox{if }k\equiv 0\mod \frac{q}{\delta}\\ -1&\mbox{if } k\equiv 1\mod \frac{q}{\delta}\\ 1&\mbox{otherwise} \end{array}\right.$. We let $G$ be one of these groups.

	\subsubsection{Local factors of $\zeta_G^{G,\s}(s)$}
	 If $p\nmid d$, then $\zeta_{G_p}^{G_p,\s}(s)$ was given in Table \ref{Table with zeta_G^G} (the case $F=C_d$; note that $\chi_6(p)=\chi_3(p)$ if $p\nmid 6$).
	 
	  Assume now that $p\mid d$. We claim that if $A\leq G$ is open and satisfies $AN_p=G_p$, then $Z_p\subset A_p$. Indeed, $A$ contains $\gamma\x^\v$ for some $\v\in\mb{Z}_p^3$, whence $(\gamma\x^\v)^{d}
	=  x_3^{d v_3+\epsilon-(r^2-r+1)\delta^2k\frac{(v_1+v_2)(v_1+v_2+r)}{2}}\in A$. The claim will follow from the assertion that $d v_3+\epsilon-(r^2-r+1)\delta^2k\frac{(v_1+v_2)(v_1+v_2+r)}{2}\in\mb{Z}_p^*$. If $d=4$, then $p=2$ and $\delta=2$, and the assertion is clear. The same argument works for $d=6$ and $p=2$. Assume that $d=6$ and $p=3$. Since $r=-1$, we have $r^2-r+1=3$; hence the assertion is true in this case too. We finally assume that $d=3$ and $p=3$. In this case $r=1$, whence $r^2-r+1=1$. Therefore, the assertion is that $\epsilon-k\frac{(v_1+v_2)(v_1+v_2+1)}{2}\in\mb{Z}_3^*$. If one of the elements $v_1+v_2, v_1+v_2+1$ is divisible by 3, the assertion is clear. Otherwise, necessarily $\frac{(v_1+v_2)(v_1+v_2+1)}{2}\equiv 1\mod 3$, so the assertion in this case is equivalent to $\epsilon-k\in\mb{Z}_3^*$, and this follows from the definition of $\epsilon$ in terms of $k$. This completes the proof of the claim. It follows that $\zeta_{G_p}^{G_p,\s}(s)=\zeta_{E_p}^{E_p,\s}(s)$. According to \cite[5.10, 5.13, 5.16]{Mc}, the latter series is equal to $\zeta_p(s-1)L_p(s-1,\chi_d)$ if $d\in\{3,4\}$ and to $\zeta_p(s-1)L_p(s-1,\chi_3)$ if $d=6$.

	\subsubsection{Local factors of $\zeta_{G}^{N,\n}(s)$}
	 If $p\nmid d$, then $\zeta_{G_p}^{N_p,\n}(s)$ was given in Table \ref{Table with zeta_G^N} (the case $F=C_d$). 
	 
	 Assume now that $p\mid d$.	By Corollary \ref{integral expression for zeta_G^N} and Lemma \ref{obtaining cone conditions}(2), we have
	$$\zeta_{G_p}^{N_p,\n}(s)=\frac{1}{(1-p^{-1})^3}\int_{\mc{T}}|t_{11}|_p^{s-1}|t_{22}|_p^{s-2}|t_{33}|_p^{s-3}d\mu,$$
	where $\mc{T}\subset\on{T}_3^+(\mb{Z}_p)$ is defined by the conditions $t_{33}|\delta k t_{11}, t_{33}|\delta k t_{12}, t_{33}|\delta k t_{22}$ and the following equivalent conditions:
	\begin{align*}
		& {}^{\gamma}(\x^{\t_1}), {}^{\gamma}(\x^{\t_2}), {}^{\gamma}(\x^{\t_3})\in\overline{\langle \x^{\t_1},\x^{\t_2},\x^{\t_3}\rangle}\\
		\Leftrightarrow \quad & x_1^{-t_{12}}x_2^{t_{11}-rt_{12}}x_3^{-k\delta t_{11}t_{12}+\binom{t_{12}}{2}\delta k r+t_{13}},\ x_1^{-t_{22}}x_2^{-rt_{22}}x_3^{\binom{t_{22}}{2}\delta k r+t_{23}}\in\overline{\langle \x^{\t_1},\x^{\t_2},\x^{\t_3}\rangle}
	\end{align*}
	By Lemma \ref{obtaining cone conditions}(3), the conditions defining $\mc{T}$ are equivalent to
	\vskip 0.1cm
	\begin{tabular}{lll}
		(1) $t_{33}|\delta k t_{11}$, $t_{33}|\delta k t_{12}$, $t_{33}|\delta k t_{22}$,&
		(2) $t_{11}|t_{12}$,&
		(3) $t_{22}|\frac{t_{12}^2}{t_{11}}+t_{11}-rt_{12}$,\\
		(4) $t_{33}|\frac{t_{12}}{t_{11}}t_{13}-\frac{\frac{t_{12}^2}{t_{11}}+t_{11}-rt_{12}}{t_{22}}t_{23}+\binom{t_{12}}{2}\delta k r+t_{13}$,&
		(5) $t_{11}|-t_{22}$,&
		(6) $t_{22}|\frac{t_{22}t_{12}}{t_{11}}-rt_{22}$,\\
		(7) $t_{33}|\frac{t_{22}}{t_{11}}t_{13}-(\frac{t_{12}}{t_{11}}-r)t_{23} +\binom{t_{22}}{2}\delta k r+t_{23}$.
	\end{tabular}
	\vskip 0.1cm
		\noindent  Condition (3) can be written as $\frac{t_{22}}{t_{11}}|(\frac{t_{12}}{t_{11}})^2-r\frac{t_{12}}{t_{11}}+1$.
		Condition (6) is implied by (2), and the second and third conditions in (1) are implied by the first one, (2) and (5).  
	Note that (4) and (7) have the form 
	\begin{enumerate}
		\item[(4)] $t_{33}|P_1 t_{13}+Q_1 t_{23}+\binom{t_{12}}{2}\delta k r $,
		\item[(7)] $t_{33}|P_2 t_{13}+Q_2 t_{23}+\binom{t_{22}}{2}\delta k r$,
	\end{enumerate}
	where $P_1$, $P_2$, $Q_1$ and $Q_2$ are rational functions that do not involve $t_{13}$ or $t_{23}$ and take values in $\mb{Z}_p$ when restricted to $\mc{T}$. Note also  that
	\begin{align}\label{determinant}
		\det\left(\begin{matrix}
			P_1&Q_1\\ P_2& Q_2
		\end{matrix}\right)&=-\frac{t_{12}+t_{11}}{t_{11}}(\frac{t_{12}}{t_{11}}-r-1)+\frac{t_{22}}{t_{11}}\frac{\frac{t_{12}^2}{t_{11}}+t_{11}-rt_{12}}{t_{22}}\\
		\nonumber &=-(\frac{t_{12}}{t_{11}})^2-\frac{t_{12}}{t_{11}}+r\frac{t_{12}}{t_{11}}+r+\frac{t_{12}}{t_{11}}+1 +(\frac{t_{12}}{t_{11}})^2+1-r\frac{t_{12}}{t_{11}}=2+r.
	\end{align}

		{\em We first assume that $(d,p)\in\{(3,3), (4,2), (6,3)\}$}. 
		The terms $\binom{t_{12}}{2}\delta k r $ and $\binom{t_{22}}{2}\delta k r$ in conditions (4) and (7) can be deleted. In fact, they are already zero when $d=4$ since in this case $r=0$, and if $d\in\{3,6\}$ and $p=3$, then (1) implies that these two terms are divisible by $t_{33}$.
		Now, one easily checks that the equation $T^2-rT+1\equiv 0\mod p$
		has $1-r-r^2$ as unique solution modulo $p$  in $\mb{Z}_p$, and that the equation $T^2-rT+1\equiv 0\mod p^2$ has no solutions in $\mb{Z}_p$. Therefore, by
		  (3), we can express $\mc{T}$ as a disjoint union $\mc{T}^a\cup\mc{T}^b$ according to the cases (3a) $|t_{11}|_{p}=|t_{22}|_{p}$, and (3b) $|p t_{11}|_{p}=|t_{22}|_{p}$ and  $pt_{11}|t_{12}-(1-r-r^2) t_{11}$. 
		  
				We describe $\mc{T}^a$. By (3a), the coefficient $P_2$ in (7) is a unit; therefore, by (\ref{determinant}), we can replace (4) by $t_{33}|(2+r)t_{23}$. We conclude that $\mc{T}^a$ is defined by the conditions $t_{33}|\delta k t_{22}$, $|t_{11}|_{p_0}=|t_{22}|_{p_0}$, $t_{22}|t_{12}$, $t_{33}|(2+r)t_{23}$ and (7), and the integral over $\mc{T}^a$ becomes
				\begin{align*}
					&\frac{1}{(1-p^{-1})^3}\int_{\mc{T}^a} |t_{22}|_p^{2s-3}|t_{33}|_p^{s-3}d\mu=\frac{1}{(1-p^{-1})^2}\int_{\substack{t_{33}|\delta k t_{22}\\ t_{33}|(2+r)t_{23}}}|t_{22}|_p^{2s-1}|t_{33}|_p^{s-2}d\mu\\
					&=\zeta_p(2s)+\frac{1}{(1-p^{-1})^2}\int_{\substack{p|t_{33}\\ t_{33}|\delta k t_{22}\\ \frac{t_{33}}{(2+r)}|t_{23}}}|t_{22}|_{p}^{2s-1}|t_{33}|_{p}^{s-2}d\mu=\zeta_p(2s)+\frac{p^{v_p(r+2)}}{(1-p^{-1})^2}\int_{\substack{p|t_{33}|\delta k t_{22}}}|t_{22}|_{p}^{2s-1}|t_{33}|_{p}^{s-1}d\mu,
								\end{align*}
							where in the first equality we used Proposition \ref{proceso de simplifiacion de integrales} with the pivots $(t_{12},t_{11},t_{13})$, and in the third one we used this proposition with the pivot $t_{23}$.
							
					We now describe $\mc{T}^b$. Condition (3b) implies that $|t_{11}|_{p}=|t_{12}|_{p}$. This makes condition (2) redundant. Next, note that $Q_1=-\frac{(\frac{t_{12}}{t_{11}})^2+1-r\frac{t_{12}}{t_{11}}}{\frac{t_{22}}{t_{11}}}$ is a unit of $\mb{Z}_{p}$ since we are assuming that $p|\frac{t_{22}}{t_{11}}$ and since $T^2-rT+1\equiv 0\mod p^2$ has no solutions in $\mb{Z}_{p}$. Therefore, by (\ref{determinant}), condition (7) can be replaced by the condition $t_{33}|(2+r)t_{13}$. We conclude that the conditions defining $\mc{T}^b$ are $t_{33}|\delta k t_{11}$,  $|p t_{11}|_{p}=|t_{22}|_{p}$, $p t_{11}|t_{12}-(1-r-r^2) t_{11}$, $t_{33}|(2+r)t_{13}$, and (4). The integral over $\mc{T}^b$ becomes
		\begin{align*}
		\frac{1}{(1-p^{-1})^3}\int_{\mc{T}^b}|t_{11}|_{p}^{s-1}|pt_{11}|^{s-2}|t_{33}|_{p}^{s-3}d\mu =\frac{1}{(1-p^{-1})^2}p^{-s}\int_{\substack{t_{33}|\delta k t_{11},\\ t_{33}|(2+r)t_{13} }}|t_{11}|_{p}^{2s-1}|t_{33}|_{p}^{s-2}d\mu,
	\end{align*}				
where in the first equality we used Proposition \ref{proceso de simplifiacion de integrales} with the pivots  $(t_{22}, t_{12}, t_{23})$. 

We conclude that
\begin{align*}
	\zeta_{G_p}^{N_p,\n}(s)&=\frac{1+p^{-s}}{(1-p^{-1})^2}\int_{\substack{t_{33}|\delta k t_{11},\\ t_{33}|(2+r)t_{13} }}|t_{11}|_{p}^{2s-1}|t_{33}|_{p}^{s-2}d\mu= \zeta_p(s)+\frac{1+p^{-s}}{(1-p^{-1})^2}\int_{\substack{p|t_{33}\\ t_{33}|\delta k t_{22}\\ \frac{t_{33}}{(2+r)}|t_{23}}}|t_{22}|_{p}^{2s-1}|t_{33}|_{p}^{s-2}d\mu\\
	&=\zeta_{p}(s)+|2+r|_{p}^{-1}\frac{1+p^{-s}}{(1-p^{-1})^2}\int_{p|t_{33}|\delta k t_{22}} |t_{22}|_p^{2s-1}|t_{33}|_p^{s-1}d\mu=\zeta_{p}(s)+|2+r|_{p}^{-1}(1+p^{-s})B_{p,\delta k}(s,2s).
\end{align*}

		 {\em We now assume that $(d,p)=(6,2)$.} The equation $T^2-r T+1\equiv 0\mod 2$ has no solutions in $\mb{Z}_2$, and hence  (3) and (5) can be replaced by $|t_{11}|_2=|t_{22}|_2$. It follows that the coefficient $P_2$ in condition (7) is a unit. Since the value of the determinant (\ref{determinant}) is $2+r=-1$, condition (4) can be replaced by one of the form $t_{33}|t_{23}+R_1$, where $R_1$ is a rational function that does not involve $t_{13}$ or $t_{23}$ and takes values in $\mb{Z}_2$ when restricted to $\mc{T}$. The conditions defining $\mc{T}$ are therefore $t_{33}|2 k t_{11}$, $t_{11}|t_{12}$, $|t_{11}|_2=|t_{22}|_2$, $t_{33}|t_{23}+R_1$, and (7). The integral over $\mc{T}$ becomes
		 	\begin{align*}
		 	\zeta_{G_2}^{N_2,\n}(s)=\frac{1}{(1-2^{-1})^3}\int_{\mc{T}}|t_{11}|_2^{2s-3}|t_{33}|_2^{s-3}d\mu=\frac{1}{(1-2^{-1})^2}\int_{t_{33}|2 k t_{11} }|t_{11}|_2^{2s-1}|t_{33}|_2^{s-1}d\mu=A_{2,2k}(s,2s),
		 \end{align*}
	 where in the first equality we used Proposition \ref{proceso de simplifiacion de integrales} with the pivots $(t_{12}, t_{22}, t_{23},t_{13})$.

\subsubsection{Local factors of $\zeta_{G}^{G,\n}(s)$}\label{G_pd, zeta_G^G, normal} If $p\nmid d$, then  $\gamma_{3}(G_p,N_p)=\overline{\langle x_1^{r+2},x_2^{r+2},x_3^{\delta k}\rangle}=\overline{\langle x_1,x_2,x_3^{ k}\rangle}$, where the last equality holds since $r+2\in\mb{Z}_p^*$. Thus, $N_p/\gamma_3(G_p,N_p)\cong\mb{Z}_p/ k\mb{Z}_p$, and by Proposition \ref{simplification for normal extension zeta function} we have $\zeta_{G_p}^{G_p,\n}(s)=\frac{1-p^{-s}|k|_p^s}{1-p^{-s}}$. 

Assume now that $p\mid d$.
Arguing as in the calculation of $\zeta_{G_p}^{G_p,\s}(s)$, we obtain that $\zeta_{G_p}^{G_p,\n}(s)=\zeta_{E_p}^{E_p,\n}(s)$. According to the results of \cite[6.10, 6.13, 6.16]{Mc}, this series is equal to $1+p\cdot p^{-s}$ if $d=3$ or $4$, and is $1$ if $d=6$.

\subsubsection{Local factors of $\zeta_{G}^{H,\n}(s)$ for $d=4$ and  $H=\langle \gamma^2,x_1,x_2,x_3\rangle$}
Note that  $H=G_{\textbf{p2},2k}$.  If $p\neq 2$, then $\gamma_3(H_p,N_p)=\overline{\langle x_1^{4}, x_2^{4}, x_3^{4k}\rangle}=\overline{\langle x_1,x_2,x_3^k\rangle}$, and hence $N_p/\gamma_3(H_p,N_p)\cong\mb{Z}_p/k\mb{Z}_p$. By Proposition \ref{simplification for normal extension zeta function}, we obtain $\zeta_{G_p}^{H_p,\n}(s)=\frac{1-p^{-s}|k|_p^s}{1-p^{-s}}$.

Assume now that $p=2$. Arguing as in the calculation of $\zeta_{G_{\textbf{p2},2k}}^{G_{\textbf{p2},2k},\n}(s)$ in \ref{G_p2,2k, zeta_G^G, normal}, we obtain that $\zeta_{G_2}^{H_2,\n}(s)=\zeta_{E_2}^{H_2/Z_2,\n}(s)$, which is equal to $1+2\cdot 2^{-s}+2\cdot 4^{-s}$ according to the results of \cite[6.10]{Mc}.

\subsubsection{Local factors of $\zeta_{G}^{H,\n}(s)$ for $d=6$ and  $H:=\langle \gamma^3,x_1,x_2,x_3\rangle$}
Note that  $H=G_{\textbf{p2},2k}$. Arguing as in the previous paragraph, we obtain that  $\zeta_{G_p}^{H_p,\n}(s)=\frac{1-p^{-s}|k|_p^s}{1-p^{-s}}$ if $p\neq 2$, and that $\zeta_{G_2}^{H_2,\n}(s)=\zeta_{E_2}^{H_2/Z_2,\n}(s)$. The latter is equal to $1+4^{-s}$ according to the results of \cite[6.16]{Mc}.

\subsubsection{Local factors of $\zeta_{G}^{K,\n}(s)$ for $d=6$ and  $K=\langle \gamma^2,x_1,x_2,x_3\rangle$}
Note that $K=G_{\textbf{p3},2k,\eta}$, where $\eta=\epsilon$ if $3\mid k$, and $\eta=-\epsilon$ if $3\nmid k$.  If $p\neq 3$, then arguing as in \ref{G_pd, zeta_G^G, normal} with $G=G_{\textbf{p3},2k,\eta}$, we find that  
$\zeta_{G_p}^{K_p,\n}(s)=\frac{1-p^{-s}|2k|_p^s}{1-p^{-s}}$, and that $\zeta_{G_3}^{K_3,\n}(s)=\zeta_{E_3}^{K_3/Z_3,\n}(s)$,  which is equal to $1+3^{-s}$ according to the results of \cite[6.16]{Mc}.

	\subsection{Computing the zeta functions of $G=G_{\textbf{p2gg},4k}$}

\subsubsection{Local factors of $\zeta_G^{G,\s}(s)$} If $p\neq 2$, then $\zeta_{G_p}^{G_p,\s}(s)$ was given in Table \ref{Table with zeta_G^G} (the case $F=D_2$). 

Assume now that $p=2$. Note that $\langle \alpha,x_1,x_2,x_3\rangle=G_{\mathbf{p2},4k}$. Arguing as in \ref{G_p2,2k, zeta_G^G,subgroup} for $G=G_{\textbf{p2},4k}$, we find that $\zeta_{G_2}^{G_2,\s}(s)=\zeta_{E_2}^{E_2,\s}(s)$. The latter is equal to 1 according to the results of \cite[5.7]{Mc}.

\subsubsection{Local factors of $\zeta_{G}^{N,\n}(s)$}
 If $p\neq 2$, then $\zeta_{G_p}^{N_p,\n}(s)$ was computed in Table \ref{Table with zeta_G^N} (the case $F=D_2$). 
 
 Assume now that $p=2$. By Corollary \ref{integral expression for zeta_G^N} and Lemma \ref{obtaining cone conditions}(2), we have
		$$
		\zeta_{G_2}^{N_2,\n}(s)=\frac{1}{(1-2^{-1})^3}\int_{\mc{T}}|t_{11}|_p^{s-1}|t_{22}|_p^{s-2}|t_{33}|_p^{s-3}d\mu,
		$$
		 where $\mc{T}\subset\on{T}_3^+(\mb{Z}_2)$ is defined by $t_{33}|4kt_{11},\ t_{33}|4kt_{12},\ t_{33}|4kt_{22}$ and the following equivalent conditions:
		\begin{align*}
			&
			^\alpha(\x^{\t_i})\cdot \x^{\t_i}, [\x^{\t_i},\beta]\in\overline{\langle\x^{\t_1},\x^{\t_2},\x^{\t_3}\rangle},\ i=1,2,3\\
			\Leftrightarrow\quad & x_3^{2t_{13}+2k(t_{11}-t_{12})-4kt_{11}t_{12}},\ x_3^{2t_{23}-4kt_{22}},\ x_2^{-2t_{12}}x_3^{-2t_{13}+2kt_{12}},\ x_3^{2kt_{22}}\in \overline{\langle\x^{\t_1},\x^{\t_2},\x^{\t_3}\rangle}.
		\end{align*}
		Therefore, by
		Lemma \ref{obtaining cone conditions}(3), the conditions defining $\mc{T}$ are equivalent to
		$$
		\begin{array}{lllllll} 
			(1)\ t_{33}|4kt_{11},\  &\ \ & (2)\ t_{33}|4kt_{12},\ &\ \ & (3)\ t_{33}|2t_{13}+2k(t_{11}-t_{12}),\ &\ \ &(4)\ t_{33}|2t_{23},\\
			(5)\ t_{22}|2t_{12},\ &\ \ &(6)\ t_{33}|\frac{2t_{12}}{t_{22}}t_{23}-2t_{13}+2kt_{12},\ &\ \ & (7)\ t_{33}|2kt_{22}.& &
		\end{array}
		$$
		Note that (2) follows from (5) and (7).
		Next, we express $\mc{T}$ as a disjoint union $\mc{T}^a\cup\mc{T}^b$ according to the following cases: a) $t_{33}\in\mb{Z}_2^*$ and b) $2|t_{33}$. 
		
		In case a), all the conditions but (5) are redundant, so
		\begin{align*}
			\frac{1}{(1-2^{-1})^3}\int_{\mc{T}^a}|t_{11}|_2^{s-1}|t_{22}|_2^{s-2}|t_{33}|_2^{s-3}d\mu=\frac{\zeta_2(s)}{1-2^{-1}}\int_{t_{22}|2t_{12}}|t_{22}|_2^{s-2}d\mu=\zeta_2(s)(1+2^{1-s}\zeta_2(s)).
		\end{align*}
		
		In case b), we split $\mc{T}^b$ as a union $\mc{T}^{b,a}\cup\mc{T}^{b,b}$ according to the following cases: (5a) $t_{22}|t_{12}$ and (5b) $|t_{22}|_2=|2t_{12}|_2$. Assume (5a). We sum the right-hand sides of (3) and (6) and obtain $t_{33}|\frac{2t_{12}}{t_{22}}t_{23}+2kt_{11}$. This and (4) imply that $t_{33}|2kt_{11}$. Thus, (1) can be replaced by $t_{33}|2kt_{11}$, and then (3) and (4)  imply (6) clearly. To sum up,  $\mc{T}^{b,a}$ is defined by the conditions $2|t_{33}$, $\frac{t_{33}}{2}|kt_{11}$, $\frac{t_{33}}{2}|kt_{22}$, $t_{22}|t_{12}$, $\frac{t_{33}}{2}|t_{13}-kt_{12}$, $\frac{t_{33}}{2}|t_{23}$. Therefore,
		\begin{align*}
			\frac{1}{(1-2^{-1})^3}\int_{\mc{T}^{b,a}}|t_{11}|_2^{s-1}|t_{22}|_2^{s-2}|t_{33}|_2^{s-3}d\mu=&\frac{2^{2-s}}{(1-2^{-1})^3}\int_{\substack{ u|kt_{11}\\ u|kt_{22}}}|t_{11}|_2^{s-1}|t_{22}|_2^{s-1} |u|_2^{s-1}d\mu=2^{2-s}C_{2,k}(s,s,s),
		\end{align*}
	 where in the first equality we used Proposition \ref{proceso de simplifiacion de integrales} with the pivots $(t_{12},t_{13},t_{23})$ and performed the change of variables $t_{33}=2u$.
		Assume now (5b).  Observe that (2) and (7) are equivalent. Note also that (6) implies that $t_{33}|2(\frac{2t_{12}}{t_{22}})t_{23}-4t_{13}+4kt_{12}$, and then that (1), (2), (3) imply that $t_{33}|4t_{13}$. It follows that (1), (2), (3), (6) imply (4). Summarizing, the conditions defining $\mc{T}^{b,b}$ are $2|t_{33}$, (1), (2), $\frac{t_{33}}{2}|t_{13}+2k(t_{11}-t_{12})$, $|t_{22}|_2=|2t_{12}|_2$ and (6). Therefore, the integral over $\mc{T}^{b,b}$ becomes 
			\begin{align*}
&				\frac{1}{(1-2^{-1})^3}\int_{\mc{T}^{b,b}}|t_{11}|_2^{s-1}|2t_{12}|_2^{s-2}|t_{33}|_2^{s-3}d\mu= \frac{2^{1-2s}}{(1-2^{-1})^3}\int_{\substack{u|2kt_{11}\\ u|2kt_{12}}}|t_{11}|_2^{s-1}|t_{12}|_2^{s-1}|u|_2^{s-1}d\mu\\
				&=2^{1-2s}(\zeta_2(s)^2+\frac{2^{-s}}{(1-2^{-1})^3}\int_{\substack{v|kt_{11}\\ v|kt_{12}}}|t_{11}|_2^{s-1}|t_{12}|_2^{s-1}|v|_2^{s-1}d\mu) =2^{1-2s}\left(\zeta_2(s)^2+2^{-s}C_{2,2k}(s,s,s)\right). 
			\end{align*}
		In the first equality we used  Proposition \ref{proceso de simplifiacion de integrales} with the pivots $(t_{22},t_{13}, t_{23})$ and applied the change of variables $t_{33}=2u$. 
	
		We conclude that  $\zeta_{G_2}^{N_2,\n}(s)=\zeta_2(s)(1+2^{1-s}\zeta_2(s))+2^{1-2s}\zeta_2(s)^2 +(2^{2-s}+2^{1-3s})C_{2,k}(s,s,s)$.

	\subsubsection{Local factors of $\zeta_{G}^{G,\n}(s)$} If $p\neq 2$, then one can easily check that $\gamma_3(G_p,N_p)=\overline{\langle x_1,x_2,x_3\rangle}$, and hence $\zeta_{G_p}^{G_p,\n}(s)=1$. 
		If $p=2$, then $\zeta_{G_2}^{G_2,\n}(s)=1$ since already $\zeta_{G_2}^{G_2,\s}(s)=1$.

	\subsubsection{Local factors of $\zeta_{G}^{H,\n}(s)$ for $H=\langle \alpha,x_1,x_2,x_3\rangle$} 
	Note that $H=G_{\textbf{p2},4k}$. Arguing as in \ref{G_p2,2k, zeta_G^G, normal} with $G=G_{\textbf{p2},4k}$,
	we find that  $\zeta_{G_p}^{H_p,\n}(s)=\frac{1-p^{-s}|k|_p^s}{1-p^{-s}}$ for $p\neq 2$, and that $\zeta_{G_2}^{H_2,\n}(s)=\zeta_{E_2}^{H_2/Z_2,\n}(s)$. The latter is equal  to $1+2\cdot 2^{-s}$ according to the results of \cite[6.7]{Mc}.

	\subsubsection{Local factors of $\zeta_G^{K,\n}(s)$ for $K=\langle \beta,x_1,x_2,x_3\rangle$} 
	Note that $K=G_{\textbf{pg},4k}$.
	If $p\neq 2$, then  $\gamma_{3}(H_p,N_p)=\langle x_1^{4}, x_3^{4}\rangle=\overline{\langle x_1,x_3\rangle}$, and hence $N_p/\gamma_3(G_p,N_p)\cong\mb{Z}_p$. Thus, by Proposition \ref{simplification for normal extension zeta function}, we have
	$\zeta_{G_p}^{K_p,\n}(s)=\zeta_p(s)$.
	
Assume now that $p=2$. 
 Arguing as in \ref{G_p2gg, zeta_G^G, normal} with $G=G_{\textbf{pg},4k}$ we find that
 $\zeta_{G_2}^{K_2,\n}(s)=\zeta_{G'}^{K',\n}(s)$, where $G'$ is the quotient of $G_2$ by  $\overline{\langle x_1,x_2^2,x_3^2\rangle}$, $K'$ is the image of $K$, and the new partial zeta function is computed with respect to $N'$, the image of $N$.
  Note that $K'\cong (\mb{Z}/2\mb{Z})^3$, the classes of $\beta$, $x_2$ and $x_3$  being identified with $(1,0,0)$, $(0,1,0)$ and $(0,0,1)$ respectively. 
  The action of the class of $\alpha$ on $K'$ is given by $(a,b,c)\mapsto (a,a+b,a+c)$. The series $\zeta_{G'}^{K',\n}(s)$ enumerates the subgroups of $G'$ that are not included in $N'$ and that are invariant under $(a,b,c)\mapsto (a,a+b,a+c)$. It is easy to check that $\zeta_{G'}^{K',\n}(s)=1+2\cdot 2^{-s}$.

\subsubsection{Local factors of $\zeta_G^{L,\n}(s)$ for $L=\langle \alpha\beta,x_1,x_2,x_3\rangle$} 
One can easily check that the assignment $\alpha\mapsto\alpha^{-1}$, $\beta\mapsto\beta\alpha$, $x_1\mapsto x_2$, $x_2\mapsto x_1$, $x_3\mapsto x_3^{-1}$ extends to an automorphism of $G$ mapping $K$ onto $L$. Therefore, $\zeta_G^{L,\n}(s)=\zeta_G^{K,\n}(s)$.

\end{small}

\bibliography{References}
\bibliographystyle{abbrv}
	\end{document}